
\ifx\CWIloaded\undefined \let\CWIloaded\relax \else   \fi

\ifx\topglue\undefined \def\topglue{\nointerlineskip\vglue-\topskip\vskip}\fi

\ifx\formloaded\undefined 
  \vsize=24.5 cm \advance \vsize -24pt 
  \hsize=16 cm
  \let\formloaded\relax 
\fi


\catcode`_=11 

\def\let_next{\let\next=}
{\obeylines\gdef\default_address%
 {Universit\'e de Poitiers, D\'epartement de Math\'ematiques,
  UFR Sciences SP2MI, T\'el\'eport~2, BP~30179, %
\ 86962 Futuroscope Chasseneuil Cedex, \ France}}

\newtoks\title \newtoks\MSCcodes \newtoks\CRcodes \newtoks\keywords
\newtoks\titlenote
\newif\if_started

 3

\font\seventeenssbf=cmssbx10 scaled \magstep 3

\font\twelvess=cmss12

\font\tenss=cmss10
\font\tenssit=cmssi10
\font\tenssbf=cmssbx10
\def\sans
{\def\rm{\fam0\tenss}\def\it{\fam\itfam\tenssit}\def\bf{\fam\bffam\tenssbf}\rm
}


\font\ninerm=cmr9 \font\nineit=cmti9 \font\ninebf=cmbx9 \font\ninesl=cmsl9
\font\ninett=cmtt9 \font\ninei=cmmi9 \font\ninesy=cmsy9
\font\niness=cmss9 \font\ninessit=cmssi9 \font\ninessbf=cmssbx10 at 9 pt
\skewchar\ninei='177 \skewchar\ninesy='60

\def\ninept
{\def\rm{\fam0\ninerm}\def\it{\fam\itfam\nineit}\def\bf{\fam\bffam\ninebf}%
 \def\sl{\fam\slfam\ninesl}\def\tt{\fam\ttfam\ninett}\rm
 \def\sans{\def\rm{\fam0\niness}\def\it{\fam\itfam\ninessit}%
           \def\sl{\fam\slfam\ninessit}\def\bf{\fam\bffam\ninessbf}\rm}%
 \textfont0\ninerm \textfont1\ninei \textfont2=\ninesy
 \textfont\itfam\nineit \textfont\slfam\ninesl
 \textfont\bffam\ninebf \textfont\ttfam\ninett
 \normalbaselineskip=11pt \normalbaselines
 \setbox\strutbox=\hbox{\vrule height 8 pt depth 3 pt width 0 pt}%
}


\long\def\centerpar#1{\begingroup
  \advance\leftskip by 0pt plus 1 fil minus 3pt \rightskip=\leftskip
	\parskip=0pt \parindent=0pt \parfillskip=0pt \linepenalty=100
  \noindent\ignorespaces#1\par\endgroup}

\def\good_break{\unskip\penalty-10\ \ignorespaces}

\catcode`@=11
\def\vfootnote#1%
{\insert\footins\bgroup \ninept
 \interlinepenalty=\interfootnotelinepenalty
 \splittopskip=\ht\strutbox 
 \splitmaxdepth=\dp\strutbox 
 \floatingpenalty=20000 
 \leftskip=1pc \rightskip=1pc \spaceskip=0pt \xspaceskip=0pt
 \parindent=0pt \parskip=3pt minus 1pt \parfillskip=0pt plus 1 fil
 \textindent{#1}\footstrut\futurelet\next\f@@t
}
\catcode`@=12

\let\_authors\empty 

%

\def\Author#1{\def\this_author{#1}%
  \ifx\this_author\_eq \let_next\Author 
  \else \let_next\let_next \afterassignment\_author 
  \fi \next}
\def\_author
 {\let\this_addr=\default_address
  \ifx\next\address \let_next\get_address
  \else \ifx\next\email \let_next\get_email
  \else \ifx\next\URL \let_next\get_URL
  \else \add_author
  \fi\fi\fi \next
}

\def\address#1{\errmessage
 {\string\address\space should only follow \string\Author}}
\def\email#1{\errmessage
 {\string\email\space should only follow \string\Author}}
\def\URL#1{\errmessage
 {\string\URL\space should only follow \string\Author}}

\def\get_address{\begingroup\obeylines\get_addr}
\def\get_addr#1%
{\endgroup\def\this_addr{#1}%
  \ifx\this_addr\_cr \let_next\get_address 
  \else\ifx\this_addr\_eq \let_next\get_address 
       \else \let_next\let_next \afterassignment\get_ad
  \fi\fi \next
}
\def\get_ad
{\ifx\next\email \let_next\get_email
 \else \ifx\next\URL \let_next\get_URL
 \else \add_author
 \fi\fi \next
}

\def\get_email#1%
{\def\this_email{#1}%
  \ifx\this_email\_eq \let_next\get_email
  \else \edef\this_addr{\this_addr\endgraf{\tt #1}}%
    \let_next\let_next \afterassignment\after_em
  \fi \next
}
\def\after_em{\ifx\next\URL \let_next\get_URL \else \add_author \fi \next}
\def\get_URL#1%
{\def\this_URL{#1}%
 \ifx\this_URL\_eq \let_next\get_URL
 \else 
  {\def~{\char "7E }%
   \xdef\this_addr{\this_addr\endgraf{\noexpand\tt URL: http://#1}}%
  }\let_next\add_author
 \fi \next
}

\def\_eq{=}{\obeylines\gdef\_cr{^^M}}

\def\add_author 
{\toks0=\expandafter{\_authors}%
 \xdef\_authors
 {\the\toks0 \noexpand\\{\this_author\noexpand\_addr{\this_addr}}}%
}


\def\start_art
{\vfil\eject \ifodd\pageno\else \shipout\vbox to \vsize{}\advancepageno\fi

 \topglue 10mm plus 3 mm minus 2mm
 \centerpar{\let\\=\good_break\seventeenssbf \baselineskip=22pt \the\title}%
 \ifx \_authors\empty \message{Warning: no authors specified.}\else
   \bigskip\centerpar
    {\twelvess
     \def\par{\ifhmode\endgraf\else\medskip\fi}\obeylines
     \def\_addr##1{\smallskip{\ninessit\let\\=\par\noindent##1\endgraf}}%
     \def\\##1{\bigskip##1}\_authors
    }\bigskip
 \fi
 \_startedtrue
}


\outer\def\abstract
{\if_started \errmessage{Abstract must be placed at start of article.}%
 \else \start_art \bigskip \begingroup \ninept \sans \narrower
   \leftline{\hskip\leftskip ABSTRACT}\nobreak\smallskip
	 \everypar{\everypar{}\setbox0=\lastbox}%
 \fi
}

\outer\def\maintext
{\if_started\medskip\endgroup\else\start_art\fi
 \if &\the\MSCcodes\the\CRcodes\the\keywords\the\titlenote&\else
  {\narrower\ninept\sans
   \if &\the\MSCcodes&\else
     \hangindent=4em \noindent
     {\it 1991 Mathematics Subject Classification:\/} \the\MSCcodes.
     \par\fi
   \if &\the\CRcodes&\else
     \hangindent=4em \noindent
     {\it 1991 Computing Reviews Subject Classification:\/} \the\CRcodes.
     \par\fi
   \if &\the\keywords&\else
     \hangindent=2em \noindent{\it Keywords and Phrases:\/} \the\keywords.
     \par\fi
   \if &\the\titlenote&\else
     {\def\\{\unskip\hfil\break\ignorespaces}%
      \hangindent=2em \noindent{\it Note:\/} \the\titlenote.
      \par}\fi
  }\bigskip
 \fi
}

\def\sticksection 
{\vbox\bgroup\parskip=0pt \everypar{\egroup\noindent}}

\def\newreport
{\vfil \eject \secno=0\subsecno=0\proclno=0 \eqnumber=0 \pageno=1 \mark{}
 \title={}\MSCcodes={}\CRcodes={}\keywords={}\titlenote={}
 \let\_authors\empty
 \_startedfalse
}

\catcode`_=8 

\ifx\formloaded\relax   \else \let\formloaded=\relax \fi

\hsize=16cm
\vsize=24cm
\normalbaselineskip=13pt\normallineskiplimit=.1pt\normallineskip.1pt
\normalbaselines
\setbox\strutbox=\hbox{\vrule height8.9pt depth4pt width0pt}
\parskip=\smallskipamount 
\multiply \smallskipamount by 4 \divide\smallskipamount by 3
\multiply \medskipamount by 4 \divide\medskipamount by 3
\multiply \bigskipamount by 4 \divide\bigskipamount by 3

\input marxmax.inc

\input epsf

\mathorddef\Part=P2 
\mathorddef\Sym=S\bffam
\mathorddef\Y=Y\bffam 

\def\set#1{{[#1]}}
\def\multiset#1{{\{\!\!\{#1\}\!\!\}}}
\def\multisetof#1:#2\endset{{\{\!\!\{\,#1\mid#2\,\}\!\!\}}}
\def\Multisetof#1:#2\endset
    {\left\{\!\!\left\{\,#1\,\,\vrule\,\,#2\,\right\}\!\!\right\}}

\mathorddef\n=n0
\def\tr{^{\rm t}}

\def\scal<#1|#2>{\left<\,#1\mid#2\,\right>}

\opdef\wt \opdef\Tab \opdef\spin \opdef\pos
\altopdef\hgt{ht} \opdef\form
\opdef\End

\def\r-{\hbox{$r$-}}
\def\precr{\prec_r}\def\succr{\succ_r}
\def\leqr{\leq_r}\def\geqr{\geq_r}
\let\leqh=\leftharpoondown
\let\geqh=\rightharpoonup
\let\leqv=\leftharpoonup
\let\geqv=\rightharpoondown
\def\leqhr{\leqh_r}
\def\geqhr{\geqh_r}
\def\leqvr{\leqv_r}

\def\rwt#1#2{\left|#1\if|#2|\else/#2\fi\right|_r}

\def\spinv{\spin\tr}

\let\meet=\wedge \let\join=\vee

\mathchardef\edgemap=\delta
\mathchardef\Edgemap=\Sigma
\def\edge#1{\edgemap(#1)}
\def\Edge#1{\Edgemap(#1)}
\def\Edgep#1{\Edgemap^+(#1)}
\def\Edgem#1{\Edgemap^-(#1)}

\def\multichoose(#1,#2){\left(\!{#1\choose#2}\!\right)}

\def\bKn{b^{\rm Kn}} \def\bBu{b^{\rm Bu}}
\def\bSW{b^{\rm SW}}\def\bWS{b^{\rm WS}}
\def\bBVG{b^{\rm BVG}}

\def\bA{b^{\rm A}}

\def\RH{{\bf RH}}

\mathorddef\n=n0 
\def\rev#1{#1^\diamondsuit}

\mathorddef\d=d0

\def\heads#1#2{I(#1/#2)}

\def\cramp#1#2#3 
  {\setbox#3\hbox{${}\atop#1'$}\def#2{\hss\raise3.5pt \box#3\hss}}

\hyphenation{semi-standard}
\ifx\Tableauxloaded\relax   \else \let\Tableauxloaded\relax \fi


\newcount\cols
{\catcode`,=\active\catcode`|=\active
 \gdef\Young(#1){\hbox{$\vcenter
 {\mathcode`,="8000\mathcode`|="8000
  \def,{\global\advance\cols by 1 &}%
  \def|{\cr
        \multispan{\the\cols}\hrulefill\cr
        &\global\cols=2 }%
  \offinterlineskip\everycr{}\tabskip=0pt
  \dimen0=\ht\strutbox \advance\dimen0 by \dp\strutbox
  \halign
   {\vrule height \ht\strutbox depth \dp\strutbox##
    &&\hbox to \dimen0{\hss$##$\hss}\vrule\cr
    \noalign{\hrule}&\global\cols=2 #1\crcr
    \multispan{\the\cols}\hrulefill\cr%
   }
 }$}}
 \gdef\Skew(#1:#2){\hbox{$\vcenter
 {\mathcode`,="8000\mathcode`|="8000
  \dimen0=\ht\strutbox \advance\dimen0 by \dp\strutbox
  \def\boxbeg{\vbox
    \bgroup\hrule\kern-0.4pt\hbox to\dimen0\bgroup\strut\vrule\hss$}%
  \def\boxend{$\hss\egroup\hrule\egroup}%
  \def,{\boxend\boxbeg}%
  \def|##1:{\boxend\vrule\egroup\nointerlineskip\kern-0.4pt
    \moveright##1\dimen0\hbox\bgroup\boxbeg}%
  \def\\##1\\##2:{\boxend\vrule\egroup\nointerlineskip\kern-0.4pt
    \kern ##1\dimen0\moveright##2\dimen0\hbox\bgroup\boxbeg}%
  \moveright#1\dimen0\hbox\bgroup\boxbeg#2\boxend\vrule\egroup
 }$}}
}


\font\sevenit=cmti7 \font\fiveit=cmti7 at 5 pt
\scriptfont\itfam=\sevenit \scriptscriptfont\itfam=\fiveit

\def\smallsquares
{\textfont0=\scriptfont0 \scriptfont0=\scriptscriptfont0
 \textfont1=\scriptfont1 \scriptfont1=\scriptscriptfont1
 \textfont\itfam=\scriptfont\itfam \scriptfont\itfam=\scriptscriptfont\itfam
 \textfont\bffam=\scriptfont\bffam \scriptfont\bffam=\scriptscriptfont\bffam
 \setbox0=\hbox{\raise0.4pt\hbox{$($}}
 \setbox\strutbox=\hbox{\vrule width 0pt height\ht0 depth\dp0 }
}

\def\bigsquares
{\setbox\strutbox=\hbox
 {\vrule width 0pt height1.3\ht\strutbox depth1.3\dp\strutbox}
}

\newbox\hpair \newbox\vpair
\def\installboxes
{\dimen0=\ht\strutbox \advance\dimen0 by \dp\strutbox
 \setbox0=\vbox{\hrule width \dimen0}
 \setbox1=\rlap{\strut\vrule}
 \setbox\hpair=\rlap{\raise\ht\strutbox\copy0 \lower\dp\strutbox\copy0}
 \setbox\vpair=\rlap{\raise\dimen0\copy1 \kern\dimen0 \copy1}
}

\newcount\ribbonlength
\newcount\joff \newcount\jmax \newcount\imid \newcount\jmid
\def\ribbon #1,#2:#3;#4 
{\setbox2=\hbox
 {\count0=0 \count1=0
  \copy1 \rlap{\lower\dp\strutbox\copy0}\if|#4|\else\ribs#4 \fi
  \raise\count1\dimen0\rlap
    {\kern\count0\dimen0 \raise\ht\strutbox\copy0 \copy1\raise0.4pt\copy1}%
  \global\joff=\count0
  \divide\dimen0 by 2 \raise\imid\dimen0\rlap{\kern\jmid\dimen0
    \hbox to 2\dimen0{\hfil$#3$\hfil}}%
 }%
 \lower#1\dimen0\rlap{\kern#2\dimen0\box2}%
 \advance\joff by #2\relax \ifnum \joff>\jmax \global\jmax=\joff \fi
}
\def\ribs#1#2
{\raise\count1\dimen0\rlap
   {\kern\count0\dimen0 \copy\ifnum #1=0 \hpair \else \vpair \fi}%
 \advance\ribbonlength by -2
 \ifnum \ribbonlength=0 \jmid=\count0 \imid=\count1 \fi
 \advance\count#1 by 1
 \ifnum \ribbonlength=1 \jmid=\count0 \imid=\count1 \fi
 \ifnum \ribbonlength>-1 \ifnum \ribbonlength<2
   \advance \jmid by \count0 \advance \imid by \count1
 \fi\fi
 \if |#2|\else \ribs#2 \fi 
}

\def\rtab#1
{\installboxes \ribbonlength=#1 \global\jmax=0
 \setbox0=\hbox{\aftergroup\endrtab \aftergroup}}
\def\endrtab{\advance\dimen0 by \jmax\dimen0 \wd0=\dimen0 \vcenter{\box0}}

\def\edgeseq#1,#2;#3 
{\setbox2=\hbox{\count0=0 \count1=0 \edges#3
  \global\joff=\count0
  }%
 \rlap{\kern#2\dimen0
       \dimen1=#1\dimen0 \advance\dimen1 by -\dimen0 \lower\dimen1\box2
      }%
 \advance\joff by #2\relax \ifnum \joff>\jmax \global\jmax=\joff \fi
}
\def\edges#1#2
{\rlap
 {\kern\count0\dimen0
  \ifnum #1=0
     \dimen1=\count1\dimen0 \advance\dimen1 by -\dp\strutbox
     \raise\dimen1\copy0
  \else \raise\count1\dimen0\copy1
  \fi}%
 \advance\count#1 by 1
 \if |#2|\else \edges#2 \fi 
}

\def\arrow#1,#2:#3
{\smash{\lower#1\dimen0\rlap{$\kern#2.5\dimen0
 \if r#3 \,\rightarrow \else
 \if l#3 \llap{$\leftarrow\,$} \else
 \if u#3 \setbox0=\hbox to 0pt{\hss$\uparrow$\hss}\ht0=0pt 
         \raise.5\dimen0 \box0 \else
 \if d#3 \setbox0=\hbox to 0pt{\hss$\downarrow$\hss}\dp0=-2pt
         \lower.5\dimen0 \box0 \else
 \if i#3 \setbox0=\hbox to 0pt{\hss$\nwarrow\,$}\ht0=0pt
         \raise.5\dimen0 \box0 \else
 \if o#3 \setbox0=\hbox to 0pt{$\,\searrow$\hss}\dp0=-2pt
         \lower.5\dimen0 \box0 \else
 \fi\fi\fi\fi\fi\fi
 $}}%
}


{\catcode`\"=\active\gdef\letqq{\let"}
 \catcode`;=\active\gdef\letsemi{\let;}
 \catcode`[=\active\gdef\letbrac{\let[}
 \catcode`\`=\active\gdef\letq{\let`} 
}


\newdimen\progindent \newif\ifbol
\mathdef\semic=Punct ;
\setbox0=\hbox{\bf do }\progindent=\wd0
\def\progcodes
{\mathcode`0=`0\mathcode`1=`1\mathcode`2=`2\mathcode`3=`3\mathcode`4=`4%
 \mathcode`5=`5\mathcode`6=`6\mathcode`7=`7\mathcode`8=`8\mathcode`9=`9%
 \mathcode`"="8000\mathcode`;="8000\mathcode``="8000\mathcode`!="8000%
}
\def\progQ{\ifbol \bolfalse \else \ \fi}
\def\progquote #1 {{}$\Q #1\/~$}
\def\progqquote #1 {{}$\nobreak\Q #1\/${}}
\def\progbq #1 {{}$\bolfalse#1\/~$}
\def\progbqq #1 {{}$\bolfalse#1\/${}}
\def\progactsemi {\semic\penalty\relpenalty\ }
\def\progtab #1{\hbox to \progindent{$#1 $ \hfil}}
\def\progcomma {\progtab,}
\def\progsemi {\progtab\semic}
\def\progcom #1{\ [\hbox{\rm#1\/}]}
\def\Qcr<#1>{\ifdim \hsize<15cm \<#1>\fi}
\def\Qbr<#1>{\Qcr<#1>$\bolfalse$}
\def\progendline{\ifhmode \errmessage{Improper line in program.}$\fi
   \ifmmode $\hskip\rightskip\egroup
            \ifdim\wd0>\hsize\line{\unhbox0}\else\box0\fi \fi}
\def\progskip{\progendline\needfil\smallskip}

\def\progdefs
{\let\Q\progQ \letq\progquote \letqq\progqquote \let\`\progbq \let\"\progbqq
 \letsemi\progactsemi \let\tab\progtab \let\;\progsemi 
 \let\-\,\let\,\progcomma \let\com\progcom \let\par\progendline
 \mathsurround 0pt\normalbaselines\boltrue\bf
}

\def\prog: 
 {\ifvmode\else\setbox0=\vtop\fi\bgroup \progdefs\progcodes
  \def\<##1>{\progendline \setbox0=\hbox\bgroup
                                   \hskip\leftskip\hskip##1\progindent$}
  \hbox{\strut}\vskip-\baselineskip}
\def\endprog{\progendline\egroup \ifvmode\else\box0\fi}

\def\progline #1\endprog
{\setbox0=\hbox{\progdefs\progcodes\def\<##1>{}$#1$}%
 \ifhmode\unhbox0\else\box0\fi}

\def\algol68
 {\prog:
  \everymath={\fam\itfam}\let\\=\progskip
  \def\com ##1{\hbox{\rm\everymath{} \# {##1\/} \#}}
  \def\plusab{\mathbin{+{:=}}}\def\minusab{\mathbin{-{:=}}}
  \mathdef\lbrack=Open [\mathcode`[="8000
  \def\next {\lbrack\futurelet\next\optspace} \letbrac\next
    \def\optspace {\if ]\next \- \else\if ,\next \- \fi\fi}
 }

\def\pascal:
 {\prog: \progindent=.7em
  \everymath={\fam\itfam}\def\str##1{\hbox{\rm'##1'}}
  \let\.=\ldotp 
  \def\com ##1{\ \{\hbox{\rm\everymath{}##1\/}\}}
 }

\def\MuPAD:
 {\prog:
  \everymath={\fam\itfam}\def\str##1{\hbox{\tt"##1"}}
  \def\for{\mathrel{\hbox{\bf\$}}}%
  \def\:{\ldotp\ldotp} 
  \def\com ##1{\kern1cm \hbox{\rm\everymath{}// ##1}}
 }

\def\ifabsent#1\intoks#2%
{\defnext##1#1##2##3\end {\ifx\ifabsent##2}%
 \expandafter\next\the#2#1\ifabsent\end
}
\let\intoks=\if 

\def\ifremoved#1\intoks#2{%
    \defnext##1#1##2#1##3\end{\ifx#1##3\global#2{##1##2}}%
    \expandafter\next\the#2#1#1\end}

\def\addtotoks#1#2{\ifabsent#2\intoks#1%
    \edef\next{\global#1{\the#1\noexpand#2}}\next\fi}

\newwrite\frwfile
\def\tofrw{\immediate\write\frwfile}
\def\writefrw{\immediate\openout\frwfile=\jobname.frw
              \global\let\writefrw=\tofrw \tofrw}

\newtoks\oldforward 
\newtoks\newforward 

\begingroup
  \openin15=\jobname.frw
  \ifeof15 \closein15
  \else
   \closein15
   \def\forward#1=#2{\gdef#1{#2}\addtotoks\oldforward#1}
   \input\jobname.frw
  \fi
\endgroup

\newif\ifforward
\def\deadref{{\bf?}}

\def\forward#1
{\ifx#1\undefined \global\let#1\deadref \forwardtrue
   \immediate\write16{Warning: unresolved forward reference to \string#1.}%
 \else \ifabsent#1\intoks\oldforward \forwardfalse \else\forwardtrue \fi
 \fi
 \ifforward \addtotoks\newforward#1\fi
 #1%
}

\outer\def\proclaim#1. #2#3\par
  {\secmedbreak \global\advance\proclno by 1 \reset\itemno
   \edef\lastlabel{\ifnum-1<\secno    \number\secno.\fi
                   \ifnum-1<\subsecno \number\subsecno.\fi
                   \number\proclno}
   {\interlinepenalty=300 \parskip=0pt\def\par{\endgraf\penalty200 }
    \noindent 
    \if\let\noexpand#2
     \ifx\undefined#2
      \global\let#2\lastlabel \deflab#2\letnext\ignorespaces
     \else
      \ifremoved#2\intoks\newforward 
       \ifx#2\lastlabel\else
          \ifx#2\deadref\else 
            \immediate\write16{Warning:
              forward reference to \string#2 is incorrect, run again!
              Changed from #2 to \lastlabel.}%
       \fi\fi
       \global\let#2\lastlabel
       \deflab#2\writefrw{\string\forward\string#2={#2}}%
       \letnext\ignorespaces
      \else 
       \ifremoved#2\intoks\oldforward 
        \global\let#2\lastlabel \deflab#2\letnext\ignorespaces
       \else\defnext{#2}
       \fi
      \fi
     \fi
    \else\defnext{#2}\fi
    \bf\lastlabel.\enspace#1.\enspace \sl\next#3\endgraf
    }%
   \ifdim \lastskip<\medskipamount \removelastskip\penalty55\medskip \fi
   }

\def\labelsec#1%
{\ifx\undefined#1
 \else
   \ifremoved#1\intoks\newforward 
     \ifx#1\lastlabel\else
        \ifx#1\deadref\else 
           \immediate\write16{Warning:
                forward reference to \string#1 is incorrect, run again!
                Changed from #1 to \lastlabel.}%
     \fi\fi
     \writefrw{\string\forward\string#1={\lastlabel}}%
   \else 
     \ifremoved#1\intoks\oldforward 
     \fi
   \fi
 \fi
 \global\let#1\lastlabel \deflab#1\nobreak 
}

\def\green#1{
#1
}

\itemindent=1.5pc

\readrefs
\Author={Marc A. A. van Leeuwen}
\email={maavl@math.univ-poitiers.fr}
\URL{www-math.univ-poitiers.fr/~maavl/}

\title={Spin-preserving Knuth correspondences for ribbon tableaux}
\MSCcodes={05E10}

\maketoc

\hyphenation{}

\abstract

The RSK~correspondence generalises the Robinson-Schensted correspondence by
replacing permutation matrices by matrices with entries in~$\N$, and standard
Young tableaux by semistandard ones. For $r\in\Np$, the Robinson-Schensted
correspondence can be trivially extended, using the \r-quotient map, to one
between coloured permutations and pairs of standard \r-ribbon tableaux built
on a fixed \r-core (the Stanton-White correspondence). Viewing coloured
permutations as matrices with entries in~$\N^r$ and total sum of coefficients
in each row or column equal to~$1$ (so the unique non-zero entry is one of the
$r$ a generators of~$\N^r$), this correspondence can also be generalised to
arbitrary matrices with entries in~$\N^r$ and pairs of semistandard \r-ribbon
tableaux built on a fixed \r-core; the generalisation is derived from the RSK
correspondence, again using the \r-quotient map. Shimozono and White recently
defined a more interesting generalisation of the Robinson-Schensted
correspondence to coloured permutations and standard \r-ribbon tableaux, one
that (unlike the Stanton-White correspondence) respects the spin statistic
(total height of ribbons) on standard \r-ribbon tableaux, relating it directly
to the colours of the coloured permutation. We define a construction
establishing a bijective correspondence between general matrices with entries
in~$\N^r$ and pairs of semistandard \r-ribbon tableaux built on a fixed
\r-core, which respects the spin statistic on those tableaux in a similar
manner, relating it directly to the matrix entries. We also define a similar
generalisation of the asymmetric RSK~correspondence, in which case the matrix
entries are taken from~$\{0,1\}^r$.

More surprising than the existence of such a correspondence is the fact that
its construction falsifies the conventional wisdom that Knuth correspondences
should be derived from Schensted correspondences via the method of
standardisation. Such a method does not work for general \r-ribbon tableaux,
since no ribbon Schensted insertion can preserve standardisations of
horizontal strips (with one notable exception in the case~$r=2$ of domino
tableaux). Instead, we use the analysis of Knuth correspondences by Fomin to
focus on the correspondence at the level of a single matrix entry and one pair
of ribbon strips. We define such a correspondence by a non-trivial
generalisation of the idea underlying the Shimozono-White correspondence,
which takes the form of an algorithm traversing the edge sequences of the
shapes involved. As a result of the particular way in which this traversal has
to be set up, our construction directly generalises neither the
Shimozono-White correspondence nor the RSK~correspondence: it specialises to
the transpose of the former, and to the variation of the latter called the
Burge correspondence.

Our constructions can be interpreted as bijective proofs of certain generating
series identities. With
$G^{(r)}_\\(q,X)
=\sum_Pq^{2\spin(P)}X^{\wt(P)}\in\Z[q][[X]]$, summed over
semistandard \r-ribbon tableaux~$P$ of shape~$\\$, the first such
identity is
$$\sum_{\lambda\geqr(0)}G^{(r)}_\lambda(q^{1\over2},X)
                        G^{(r)}_\lambda(q^{1\over2},Y)
 =\prod_{i,j\in\N}\prod_{k=0}^{r-1}{1\over1-q^kX_iY_j};
$$
this can be considered the $q$-analogue of an \r-fold Cauchy identity, since
for $q=1$ each $G^{(r)}_\\(q^{1\over2},X)$ factors into a product of
$r$~Schur functions. It is equivalent to a commutation relation for certain
operators acting on a $q$-deformed Fock~space, obtained by Kashiwara, Miwa
and~Stern. Our asymmetric correspondence bijectively proves
$$\sum_{\lambda\geqr(0)}G^{(r)}_\lambda(q^{1\over2},X)
             \widetilde G^{(r)}_\lambda(q^{1\over2},Y)
 =\prod_{i,j\in\N}\prod_{k=0}^{r-1}(1+q^kX_iY_j).
$$
where
$\widetilde G^{(r)}_\\(q,X)
=\sum_Pq^{2\spinv(P)}X^{\wt(P)}$ is the counterpart of
$G^{(r)}_\\(q^{1\over2},X)$ for transpose semistandard \r-ribbon tableaux,
their spin being defined using the standardisation appropriate for such
tableaux.

\maintext

\newsection Introduction.

The Robinson-Schensted correspondence has been generalised by various authors
in many different ways. Fomin has even described general schemes that allow
defining variants of this correspondence for any combinatorial structures that
satisfy certain basic relations. We shall apply such a scheme to find what can
be described as a Knuth correspondence extending the Schensted correspondence 
recently defined by Shimozono and White~\ref{Shimozono White}; it is based
however on a quite novel kind of basic construction. 

To indicate where our constructions fit in, we shall first need to review
various earlier generalisations of the Schensted algorithm, and Fomin's
general constructions. That however involves a rather long and abstract
discussion, which does not transmit very well the flavour of the operations we
shall actually be concerned with. Therefore we prefer, in order to whet the
reader's appetite, to first present some enumerative problems that do not
require much background, and yet are quite close to the questions related to
our main construction; indeed the enumerative claims we formulate below will
follow as special case from that construction. These problems allow us to
introduce in an informal manner several important ideas behind our
constructions. This discussion is for motivation only, so readers who wish to
skip this somewhat oversize \foreign{hors d'{\oe}uvre} can move on
to~\Sec\forward{\backgroundsec} without problem; the remainder of the paper
will explicitly provide any required notions and results where and when those
are needed.

\subsection Some enumerative problems.
\labelsec\enumsec

Fix an integer $r>0$ and an arbitrary bit string: a word~$w$ over the alphabet
$\{0,1\}$. To~$w$ we associate a lattice path: the successive bits of~$w$
determine the directions of the successive steps of the path, going a unit to
the right for each bit~$0$, and an unit upwards for each bit~$1$. We extend
this path indefinitely at both ends by steps in a fixed direction; for our
first problem we extend vertically at both ends (as if~$w$ floats in a sea of
bits~$1$). We shall count the number of ways to place a collection of
\r-ribbons below the path, according to the following rules. An \r-ribbon is a
polygon built up from sequence of $r$ squares arranged from bottom-left to
top-right, each following square being either directly above or directly to
the right of its predecessor. The first \r-ribbon placed must have all of its
top-left border along the (extended) path corresponding to~$w$ (one easily
sees that there are $r+1$ consecutive segments along the border of an
\r-ribbon that are either the left or top edge of one of its squares). For the
purpose of placing further \r-ribbons, the path is modified by replacing the
top-left border of this first ribbon by its bottom-right border, so that
further ribbons will be to the bottom-right of the first one. There is an
additional restriction however: the final (top-rightmost, and therefore
horizontal) segment of the top-left border of each ribbon placed must be part
of the \emph{original} path corresponding to~$w$, and further to the top-right
than (the segments on the border of) any previously placed ribbons. Since
every ribbon placed so ``uses'' at least one bit~$0$ of~$w$, it is clear that
the number of ribbons is bounded by the number of such bits. The final
restriction also ensures that no collection of ribbons can be constructed in
more than one way by ordering its ribbons differently. Here is a example of a
placement of five 4-ribbons below the path corresponding to
$w=0010110000010000$:
\bigdisplay
\smallsquares
{\rtab4
  \edgeseq9,0;11
  \ribbon6,0:;110 \ribbon5,1:;011 \ribbon4,3:;111 \ribbon2,4:;100
  \edgeseq1,7;0
  \ribbon0,8:;000
  \edgeseq0,12;11
}
$$
Rather than just counting the total number of placements possible, we refine
the counting by keeping track independently of two ``statistics'' for each
placement: the first is the number~$n$ of \r-ribbons placed, and the second is
their ``total height''~$t$, namely the number of vertical adjacencies between
squares of a same ribbon. In the example displayed one has $n=5$, and
$t=2+2+3+1+0=8$, where the sum shows the contributions of the individual
ribbons, in order of placement. By contributing a monomial $X^nY^t$ for each
placement, and taking the sum over all placements of collections of
\r-ribbons below the path corresponding to~$w$, one obtains a two-variable
generating polynomial $R_{1,1}(w,r)\in\Z[X,Y]$ (the two indices~`$1$' are
there to indicate that we have extended~$w$ by bits~$1$ at both ends). While
we do not know any more direct way of describing these polynomials, we do
remark the following property.

\proclaim Claim. \alpineclaim
If $\widetilde w$ is the word obtained by reversing the bits of~$w$, then
$R_{1,1}(w,r)=R_{1,1}(\widetilde w,r)$.

Although the polynomials~$R_{1,1}(w,r)$ are rather laborious to compute by
hand, their computation can be quite easily programmed. The basic observation
is that after having prefixed $w$ by $r$ bits~$1$ (more are not necessary),
each possible placement of a first \r-ribbon is characterised by the
simultaneous occurrence of a bit~$1$ and a bit~$0$ exactly $r$ places to
its right, and that the modification of the path due to placing the ribbon
corresponds to changing those two bits and nothing else. The possibilities of
adding further ribbons can be computed recursively if one takes care to ensure
that they can only be placed further to the right. This can be achieved by
removing in the recursive call the initial part of the word that may no longer
be altered, i.e., the part up to and including the first bit that changed
(from~$1$ to~$0$; since the bit disappears anyway there is no need to actually
perform this change). To illustrate the simplicity of the algorithm, we
present the complete code in the language of the MuPAD computer algebra
system. We hope that this is readable even to those not familiar with MuPAD;
comments are given after the symbol~``//''. The only technical point is the
procedure {\it R11} which prefixes $r$ bits~$1$ to the word before entering
the recursion, and makes sure the result is expanded and presented as a
polynomial in~$\Z[Y][X]$.
\bigdisplay
\eqalign{
R &:=
\MuPAD:
\<0>"proc (w,r) \com{$w$ is a list $[\li(w1..l)]$ with $w_i\in\{0,1\}$}
\<1>   `local l,result,i,j,ww;
\<0>`begin l:=nops(w); result:=1 \com{count the solution with no ribbons}
\<0>\; `for i `from 1 `to l-r \com{when $l\leq r$, the loop is skipped}
\<1>  \tab{"do }
          `if w[i]=1 `and w[i+r]=0
\<2>      `then ww:=[op(w,i+1\:l)]
           \com{copy the sub-list $[\li(wi+1..l)]$ of~$w$}
\<2>      \; ww[r]:=1 \com{change the bit that was copied from $w_{i+r}$}
\<2>      \; result:=result+ R(ww,r) *
                         X * Y\uparrow\_plus(ww[j]\for j=1\:r-1) 
\<2>      \com{the exponent of~$Y$ is $\sum_{j=1}^{r-1}ww_j
				      =\sum_{j=i+1}^{i+r-1}w_j$}
\<2>    "end\_if
\<1>  "end\_for
\<0>\; result \com{return the polynomial computed}
\<0>"end\_proc ;
\endprog
\cr
{\it R11} &:= \MuPAD:
\<0>    "proc (w,r) `begin poly(poly(R([1\for r,op(w)],r),[X,Y]),[X])
        "end\_proc ;
\endprog
\cr
}
$$
Thus one may compute that for $w=[0,0,1,0,1,1,0,0,0,0,0,1,0,0,0,0]$ that
$R_{1,1}(w,4)$ equals
\bigdisplay
\displaylines
{1\cr
 +X(2+2Y+2Y^2+Y^3)\cr
 +X^2(1+4Y+7Y^2+5Y^3+5Y^4+2Y^5+Y^6)\cr
 +X^3(Y+6Y^2+10Y^3+15Y^4+12Y^5+8Y^6+5Y^7+2Y^8+Y^9)\cr
 +X^4
  (2Y^3+11Y^4+19Y^5+23Y^6+20Y^7+16Y^8+8Y^9+5Y^{10}+2Y^{11}+Y^{12})\cr
 +X^5(Y^5+10Y^6+21Y^7+32Y^8+29Y^9+24Y^{10}+16Y^{11}+8Y^{12}+5Y^{13}+
  2Y^{14}+Y^{15})\cr
 +X^6(3Y^8+12Y^9+28Y^{10}+34Y^{11}+33Y^{12}+24Y^{13}+16Y^{14}+8Y^{15}+
  5Y^{16}+2Y^{17}+Y^{18})\cr
 +X^7(Y^{11}+10Y^{12}+21Y^{13}+32Y^{14}+29Y^{15}+24Y^{16}+16Y^{17}+
  8Y^{18}+5Y^{19}+2Y^{20}+Y^{21})\cr
 +X^8(2Y^{15}+11Y^{16}+19Y^{17}+23Y^{18}+20Y^{19}+16Y^{20}+8Y^{21}+
  5Y^{22}+2Y^{23}+Y^{24})\cr
 +X^9(Y^{19}+6Y^{20}+10Y^{21}+15Y^{22}+12Y^{23}+8Y^{24}+ 5Y^{25}+
  2Y^{26}+Y^{27})\cr
 +X^{10}(Y^{24}+4Y^{25}+7Y^{26}+5Y^{27}+5Y^{28}+2Y^{29}+Y^{30})\cr
 +X^{11}(2Y^{30}+2Y^{31}+2Y^{32}+Y^{33})\cr
 +X^{12}Y^{36},\cr
}
$$
and so does $R_{1,1}(\widetilde w,4)$, where 
$\widetilde w=[0,0,0,0,1,0,0,0,0,0,1,1,0,1,0,0]$.

Our claim above can be interpreted in a geometric fashion. If a placement of
\r-ribbons below the path corresponding to~$\widetilde w$ is rotated a half
turn, one obtains a placement of \r-ribbons \emph{above} the path
corresponding to~$w$ according to similar rules as for the placement below
(but note that the order of placement is now from right to left). So the claim
can be reformulated as: for any path that is ultimately vertical at both ends,
and any specified number of ribbons and total height, there are as many
placements possible above the path as there are below the path. Even more than
the original formulation this form begs for a bijective proof, a simple rule
to move ribbons to the other side of the path, so as to define an invertible
map from placements of ribbons on one side to those on the other, preserving
the number of ribbons and the total height; one would expect the inverse map
to be given by the same rule after rotating the configuration a half turn. Yet
we have not been able to find such a rule. Our proof of the claim (given at
the end of the paper) will be based on a bijection, but one corresponding an
identity obtained by multiplying both polynomials by an identical power
series. Although a general method exists to deduce from this a bijection
corresponding directly to the claim, the result it is way too complicated to
qualify as a bijective proof.

This negative finding notwithstanding, there are quite a few observations we
can make about this problem that do involve simple bijective constructions.
For instance, when seeing polynomials $R_{1,1}(w,r)$ such as the one displayed
above, it is hard to miss a symmetry in the coefficients (although we admit
having done so for quite some time): for every~$i$, the polynomial in~$Y$ that
is the coefficient of $X^{12-i}$ equals the coefficient of $X^i$ multiplied by
$Y^{36-6i}$. In the general case the coefficient of $X^{k-i}$ equals that of
$X^i$ multiplied by $Y^{(r-1)(k-2i)}$, where $k$ is the number of bits~$0$
in~$w$. This suggests that each placement of~$i$ \r-ribbons should correspond
to a placement of $k-i$ such ribbons with the same total \emph{width} (i.e.,
the number of horizontal adjacencies of squares within a same ribbon). Indeed
it does, but the latter will be a placement of ribbons \emph{above} the path
of~$w$, making the observation about the symmetry of the coefficients
equivalent to our claim~\alpineclaim. To describe the correspondence, first
observe that a placement of \r-ribbons is completely determined by it top-left
and bottom-right borders, i.e., by the path corresponding to $w$ and that path
as modified by the placement of all ribbons. Then if one shifts up the latter
path by $r$~units, it will be the top-left border of a unique placement of
\r-ribbons above the path of~$w$. For instance, here is the placement obtained
from the one displayed earlier (both have total width~$7$):
\bigdisplay
\smallsquares
{\rtab4
  \edgeseq9,0;11
  \ribbon6,0:;010
  \edgeseq5,3;10
  \ribbon3,4:;101 \ribbon3,5:;011 \ribbon3,7:;111
  \ribbon2,8:;110 \ribbon2,9:;101 \ribbon2,10:;011
  \edgeseq0,12;11
}
$$
Providing the details to prove that this correspondence is well defined and
has the required properties is an instructive exercise that we encourage the
reader to solve; those who would like a hint can turn to
lemma~\forward{\ribbonmerge} below, which provides a closely related result.

We can also give bijections proving several sub-cases of our
claim~\alpineclaim. To begin with, the case $r=1$ does not pose much
difficulty. Now $1$-ribbons are just single squares, and the total height,
being zero always, plays no role. The placement rule forces squares to be
added from left to right below the path of~$w$, advancing at least one column
at each step, so that no column can receive more than one square. Conversely
any path that remains weakly below that of~$w$ and weakly above the same path
shifted one unit down (thus leaving room for at most one square in each
column) can be obtained for an appropriate placement of squares; the set of
squares of such a placement is known as a horizontal strip. Each horizontal
component of the path of~$w$ (a maximal portion without vertical steps) can be
treated in isolation, and can be used to place any number of squares from~$0$
up to and including its length, in a unique way; this is true both for
placement above and below~it. So one can bring the horizontal strips above and
below the path into bijective correspondence, by requiring that for the former
has as many squares directly above any horizontal component as the latter has
directly below it, in other words that the former has as many squares in any
given row as the latter has in the next row. There is an equivalent
algorithmic description that treats the squares one at a time: traversing the
squares of the horizontal strip above the path of~$w$ from left to right, each
square is moved one place down, thus crossing the path of~$w$, and then if
necessary slid to the left until it finds a place where it can stay, namely
where its left edge is part of the path of~$w$ as possibly modified by the
previous placement of squares below it.

This case can be extended to cover a small part of the claim for general~$r$,
namely the part concerning the leading terms (in~$Y$) of the coefficients of
the~$X^i$, in other words the placements where all \r-ribbons are completely
vertical. For such placements only vertical portions of the path of length at
least~$r$ effectively produce separate compartments where the ribbons can be
placed independently; shorter vertical portions within such a compartment have
no effect on the number of vertical ribbons that can be placed in it. So again
each compartment can accommodate any number of vertical ribbons from~$0$ up to
and including its width, and does so both above and below. This correspondence
too can be described by processing the ribbons above the path from left to
right: each one is moved down across the path, and then to the bottom-left up
to the first place where it will fit. Incidentally a similar procedure also
works when there are only horizontal ribbons, but these cases are even more
marginal than those involving only vertical ribbons, since generally only
relatively few purely horizontal \r-ribbons can be placed at all.

All this only scratches the surface of the general problem. It should be noted
that one cannot expect a correspondence for the general case where each ribbon
above the path gives rise to a ribbon of the same height below the path, for
the simple reason that the distributions of heights within the collections of
placements that should match are not always the same. This can be seen for the
example given (for the given~$w$, there are $5$~different placements of five
ribbons that, like the one displayed, produce a total height~$8$ as sum of the
multiset $\multiset{0,1,2,2,3}$, whereas there are only~$4$ such placements
for~$\widetilde w$), but a smaller example is more convincing: below the path
of $01000$ one can place a vertical $3$-ribbon and a horizontal one (i.e.,
$3$-ribbons of heights $2$~and~$0$, respectively), but such a pair of
$3$-ribbons cannot be placed above the same path. This means that it is
important that we count only by \emph{total} height, and that any
correspondence one would hope to find must have some mechanism for the
exchange of height between ribbons (or alternatively it might treat placements
of ribbons as a whole without even considering individual ribbons, as our
first bijection did).

There are two more bijective correspondences that are worth mentioning in this
context, as they provide interesting new points of view, even though they do
not tackle the difficulties just indicated. The first proves the
specialisation for $Y:=1$ of our claim, i.e., it treats all configurations but
ignores the heights of the ribbons. The second handles the equality of the
coefficients of~$X^1$, i.e., it treats configurations consisting of a single
ribbon. If one ignores heights, matters become simpler if one forgets the
geometric description, and views placement of ribbons simply as operations on
bit strings. As we saw, the question of whether an \r-ribbon can be placed,
and the effect of placing it, can both be expressed in terms of just one pair
of bits, at indices $i$ and~$i+r$. So placement of different \r-ribbons
becomes completely independent unless the indices $i,i'$ of the bits involved
are congruent modulo~$r$ (in the latter case we shall say the ribbons are in
the same position class). Thus the possibilities for placing \r-ribbons
decompose completely following the $r$ different position classes, and the
specialisation for $Y:=1$ of $R_{1,1}(w,r)$ decomposes as a product
$\prod_{i=0}^{r-1}R_{1,1}(w^{(i)},1)$ of polynomials in~$X$, where $w^{(i)}$
is the word extracted from~$w$, of its bits at indices congruent to~$i$
modulo~$r$. So for instance for our example, the specialisation $1+7X + 25X^2
+ 60X^3 + 107X^4 + 149X^5 + 166X^6 + 149X^7 + 107X^8 + 60X^9 + 25X^{10} +
7X^{11} + X^{12}$ factors as
$R_{1,1}(0100,1)R_{1,1}(0100,1)R_{1,1}(1000,1)R_{1,1}(0010,1)
=(1+2X+2X^2+X^3)\*(1+2X+2X^2+X^3)\*(1+X+X^2+X^3)\*(1+2X+2X^2+X^3)$. Thus we
are back at the case $r=1$ we know how to handle. We find the following
procedure to transform a placement of \r-ribbons above the path of~$w$ into
one below, defined by the final value of a modified copy~$w'$ of~$w$. Process
the ribbons of any position class from left to right; the relative ordering
between members of different classes is irrelevant. For a ribbon with initial
bit $w_i=0$, search in the current value of~$w'$ (which still has $w'_i=0$),
testing the bits $w'_{i-r},w'_{i-2r},w'_{i-3r},\ldots$ until finding the first
bit~$w'_{i-kr}=1$; one has $w'_{i-kr+r}=0$, and the bit string~$w'$ is
modified by setting $w'_{i-kr}:=0$ and $w'_{i-kr+r}:=1$. The modifications
to~$w'$ do not always occur in the right order to describe the ribbons of the
placement eventually found, so the independence of operations on different
position classes of ribbons is crucial for proving that the same procedure
rotated a half turn defines an inverse.

We have seen that preserving heights of individual ribbons is not possible in
general, but that ignoring heights altogether makes our problem trivial. The
following bijection for the case of single ribbons gives some insight in the
role played by height, without the complications of interaction between
ribbons; it is based on observations in~\ref{Shimozono White}. When an
\r-ribbon of height~$h$ can be placed below the path of~$w$ with initial bit
$w_i=1$, this means that $w_{i+r}=0$ and $\sum_{j=i+1}^{i+r-1}w_j=h$, which
can also be formulated as $\sum_{j=i}^{i+r-1}w_j=h+1$ and
$\sum_{j=i+1}^{i+r}w_j=h$. Thus the places where a ribbon of height~$h$ can be
placed below the path of~$w$ correspond exactly to the places where the sum of
$r$ consecutive bits drops from $h+1$ to~$h$, and similarly a ribbon of
height~$h$ can be placed above the path precisely in the places where the
value of such sums rises from $h$ to~$h+1$. Therefore, starting from a place
where such a ribbon can be placed below, one can always find a place further
to the top-right where a similar ribbon can be placed above (since the path
ultimately becomes vertical), and from the first such place, the point of
departure can be found back as the first place to its bottom left that will
accommodate an \r-ribbon of height~$h$ below the path. This establishes our
bijection for the case of single ribbons. One may visualise all possibilities
to place \r-ribbons of height~$h$, both below and above, as the points of
crossing between the path of~$w$ and an appropriately shifted copy of the same
path; see the illustrations after corollary~\forward{\interleavecorr} below.

We close our discussion of this problem with an indication of why we think its
has no simple bijective solution (although we would love to be proved wrong).
When one tries to extend the height preserving procedure for single ribbons to
multiple ribbons, the main difficulty is not so much the exchange of heights
that may be necessary, as the fact that the left to right order among ribbons
cannot be preserved. We believe we could describe a bijection for the case of
two ribbons, but it already gets horribly complicated: when the second ribbon
placed needs to move beyond the place where the first landed, exchange of
height must be taken into account, and it may be necessary to relocate the
first ribbon. But the hardest part is to show that one gets a bijection: the
original ribbons must be reconstructed from the pair produced without knowing
in which order those were placed, so there is no question of a step-by-step
inverse; a proof would involve piecing together all the scenarios that can
arise. Unless there is some easy way to read off the order in which ribbons
have been placed, it is hard to envisage a similar technique handling the case
of three or more ribbons.

If we have spent much time on a problem for which we know no solution, it is
because it is superficially simpler than a second problem, a variant of the
first, but one for which we do have a solution; indeed the solution is closely
related to the main result of this paper (and it will not be detailed in this
introduction). The variant is simply obtained by extending the path described
by the finite bit string~$w$ not vertically, but horizontally at both ends; in
other words that string is now considered to float in a sea of bits~$0$. The
conditions for placing collections of \r-ribbons remain exactly the same, as
are the two statistics on placements of ribbons (number of ribbons~$n$ and
total height~$t$); analogously to the definition of~$R_{1,1}(w,r)$, the sum
of~$X^nY^t$ over all possible placements below this differently extended path
will now be denoted by~$R_{0,0}(w,r)$. One still has symmetry between
placements of ribbons above and below the path.

\proclaim Claim. \polderclaim
If $\widetilde w$ is the word obtained by reversing the bits of~$w$, then
$R_{0,0}(w,r)=R_{0,0}(\widetilde w,r)$.

If, in keeping with the laws of gravity, we think primarily of placing ribbons
above the path, then the path in our first ribbon placement problem resembles
a ledge in an otherwise sheer rock-face, while the second problem more
resembles a Dutch landscape with a polder to the left, a dike described by the
string~$w$, and the sea to the right (the sea being as high as the dike is not
quite realistic, fortunately). We shall therefore refer to first ribbon
placement problem as the alpine problem, and to this second ribbon placement
problem as its polder variant.

The change of landscape modifies the character of our problem in several ways.
While ribbons can lean against the rock face, the sea and the space above sea
level are inaccessible (the top-rightmost vertical edge of each ribbon must
belong either to the dike or to another ribbon). On the other hand, the
requirement that the bottom-leftmost horizontal edge of each ribbon lie on
the original path does not put a bound on the number of ribbons, since the
polder provides an infinite supply of such edges. Indeed, provided $w$ has at
least one bit~$1$, arbitrarily many ribbons can be placed above the path, for
instance using only horizontal ribbons in the polder. Hence the identity of
our second claim is one of formal power series rather than of polynomials.
Considering $X$ to be the principal indeterminate, one has in fact
$R_{0,0}(w,r)\in\Z[Y][[X]]$: the coefficient of $X^n$ is a polynomial in~$Y$
of degree at most $n(r-1)$.

One cannot compute a complete power series~$R_{0,0}(w,r)$, but the recursive
procedure~$R$ above can be easily adapted to produce an initial part of such
power series (in finite time!), by adding as a parameter a limit to the degree
in~$X$ of the terms it should compute. Thus one verifies that for
$w=1011000001$ (our earlier example without the now superfluous bits~$0$ at
either end), both $R_{0,0}(w,4)$ and $R_{0,0}(\widetilde w,4)$ start as
\bigdisplay
\displaylines
{1\cr
+X(2+Y+Y^2)\cr
+X^2(2+2Y+4Y^2+Y^3+Y^4)\cr
+X^3(2+2Y+5Y^2+4Y^3+4Y^4+2Y^5+Y^6)\cr
+X^4(2+2Y+5Y^2+5Y^3+7Y^4+5Y^5+5Y^6+2Y^7+2Y^8)\cr
+X^5(2+2Y+5Y^2+5Y^3+8Y^4+8Y^5+8Y^6+6Y^7+6Y^8+3Y^9+2Y^{10}+Y^{11})\cr
\cdots\cr
}
$$
Like before, the configurations counted by~$R_{0,0}(w,r)$ and by
$R_{0,0}(\widetilde w,r)$ can be interpreted respectively as placements of
ribbons below and above the same path, and one would like to prove the claim
by means of a bijection between such placements, which preserves the number of
ribbons and their total height.

It is interesting to observe how the case $r=1$, that of the horizontal
strips, has changed. In terms of horizontal components of the path, we have
effectively gained one such component, with infinite capacity, whether placing
squares above or below the path (in more formal terms: assuming that $w$
neither begins nor ends with a bit~$0$, one has
$R_{0,0}(w,1)=R_{1,1}(w,1)\(\sum_{n\in\N}X^n\)$, which combinatorialists would
write as $R_{1,1}(w,1)\over1-X$). But the infinite horizontal component is not
the same one in both cases, so if one wants to maintain the bijection based on
this decomposition into horizontal components, one has to decree that, while
most squares from the strip above the path descend below it and then shift to
the left, those that were just above polder level ``wrap around at infinity''
and come back from the extreme right, just below sea level. There is nothing
against that as a bijection for $r=1$, but as point of departure for the
general case it is better to consider a different bijection, one that moves
all squares in the same direction; this must be left-to-right when going from
a horizontal strip above the path to one below. Doing so, squares arrive under
a different horizontal component than the one they belonged to, and since the
capacities of those are unrelated, the level at which squares will be placed
cannot be as neatly predictable as before.

Yet there is a simple method for placing the squares, which is essentially to
take the first available place to the right of their original position, taking
into account the other squares. For instance, for a horizontal strip above the
path that consists just of $n$ squares on the lowest possible (polder) level,
and just to the left of the first vertical edge of the path, the corresponding
horizontal strip below the path will occupy the first $n$ columns to the right
of that vertical edge, at whatever level is necessary to be just below the
path. This is possible for any $n\in\N$ because the squares below any
horizontal component of the path can be filled up from left to right, and a
horizontal component of infinite capacity is available at the end to absorb
whatever number of squares might not fit elsewhere. One wants the same
description to define the inverse procedure, which means in this example that
the horizontal strip above the path that occupies the $n$ columns directly to
the left of the \emph{last} vertical edge of the path should correspond to the
strip of $n$ squares to the right of that edge, just below sea level. To
obtain that result, one must declare columns that contain a square of the
original horizontal strip to be unavailable for placing squares, even if doing
so could produce a horizontal strip as output. Therefore the rule should be:
for each square, taken in order from left to right of their original position,
move it to the first column to its right that contains neither a square that
was already moved nor one that has yet to be moved. One may verify that a
square may indeed be placed in that column, just below the path, and that the
same rule rotated a half turn will bring back each square to its original
position.

Apart from this one-square-at-a-time procedure, there is a description of the
same correspondence that treats all squares at the same time. Imagine a bus
driving along the path from left to right, taking the squares with it as
passangers. Each horizontal edge is a bus stop, where either a square enters
the bus (from above), or a square leaves the bus (from below, but never at a
stop where a square entered), or finally, in case the bus is empty and there
is no square waiting at the bus stop, the the bus just drives on. The
importance of this alternative description is not so much its cuteness or
greater efficiency, as the fact that it treats the squares without regard to
their individual identities: while for the purpose of showing equivalence with
the earlier description, one may imagine that each square leaving the bus is
the one among the current passengers that has been aboard the longest, this
``choice'' has no effect on the result, and all that matters is keeping track
of the number of passengers at each moment.

If we consider the other correspondences that were established for the alpine
problem, we see that the one that was introduced in connection with the
symmetry in the coefficients of the polynomial~$R_{1,1}(w,r)$ has no
counterpart in the polder variant, because the power series $R_{0,0}(w,r)$ has
no such symmetry; the two other correspondences however (the one counting
placements of ribbons ignoring their heights, and the correspondence for
placements of single ribbons) can be adapted to the polder variant without
much difficulty. The first one of these involves the same reduction to the
case~$r=1$ as before, which case is modified as just discussed; the resulting
correspondence can be described by transportation of ribbon by means of a bus
with $r$ separate compartments, one for each of the position classes. The
second one is adapted by inverting the direction of search for a ribbon of
appropriate height, due to fact that sums of $r$ consecutive bits now
ultimately become~$0$ in both directions, rather than~$r$.

Unlike for the alpine problem, a bijection handling the general case and
preserving height can be given here; this fact is at the heart of our main
result. Details about the bijection will be given later, but here is a hint
for the impatient. The correspondence is obtained by the passage of an \r-deck
bus transporting ribbons, but instead of segregation by class (an idea we
could not endorse anyway) the level of entry and exit of ribbons is related to
height. This level is not identical to the height of the ribbon entering or
leaving however (that would force the output to have the same distribution of
heights as the input, which, like for the alpine problem, cannot always be
achieved). Rather it is the height a ribbon at that position would have in a
third placement, one that extends at both sides of the path of~$w$, and
occupies the union of the areas occupied but the input and output placement
(even if the latter is still under construction). Moreover, the bus operates
in a stack-like fashion: whenever a ribbon enters or leaves at some level, all
higher levels are empty.

Neither the alpine problem nor its polder variant quite reflect the original
problem that motivated our work. The problem deals with Young diagrams, whose
boundary is given by a path that starts vertically and ends horizontally. By
extending $w$ to the left with bits~$1$ and to the right with bits~$0$ one
obtains such a path, and it can be used to define a generating series
$R_{1,0}(w,r)$. For the purpose of counting placements of ribbons above rather
than below such a path one similarly defines~$R_{0,1}(\widetilde w,r)$. The
symmetry observed for the other problems does not exist for this case however;
indeed $R_{1,0}(w,r)$ is a proper power series, while $R_{0,1}(\widetilde
w,r)$ is a polynomial. The simplest case is obtained for the empty
word~$\epsilon$. Now no ribbons can be placed at all above the path, so that
$R_{0,1}(\epsilon,r)=1$. On the other hand, it can be seen that for any given
multiset of heights, there is exactly one way to place ribbons of those
heights below the path, which does so by weakly decreasing order of height.
Thus one deduces that $R_{1,0}(\epsilon,r)=\prod_{h=0}^{r-1}{1\over1-XY^h}$,
the generating series of multisets on $\set{r}=\{0,\ldots,r-1\}$. Instead of
the equalities expressed in our first two claims, we observe here that the
generating series for placements below the path is always obtained from the
generating polynomial for those above by multiplication by this fixed power
series $R_{1,0}(\epsilon,r)$.

\proclaim Claim. \partitionclaim
If $\widetilde w$ is the reverse word of~$w$, then
$R_{1,0}(w,r)=R_{0,1}(\widetilde w,r)\prod_{h=0}^{r-1}{1\over1-XY^h}$.

In spite of the different nature of the statement, a bijective proof of this
final claim will turn out to be deduced immediately from one for
claim~\polderclaim: each multiset contributing to the factor
$\prod_{h=0}^{r-1}{1\over1-XY^h}$ can be interpreted as describing an initial
state of the bus when it arrives (instead of all decks being empty initially),
the occurences of~$i\in\set{r}$ occupying deck~$i$. The bus will still leave
the scene empty in this case, but that too changes if the path ends vertically
as in the alpine problem; by allowing for a non-empty bus both at entry to and
at exit from the scene, one obtains a bijective proof not of the identity of
claim~\alpineclaim, but of the identity derived from it by multiplying both
sides by the power series~$R_{1,0}(\epsilon,r)$.

\subsection Background.
\labelsec\backgroundsec

The basic form of the Schensted algorithm constructs a bijection between
permutations of~$n$ and pairs of standard Young tableaux of equal shape and
size~$n$. The two tableaux shall be referred to as the $P$-symbol and the
$Q$-symbol, and this terminology will be extended to all generalisations of
the algorithm considered. Its first generalisation already appears in the
original paper~\ref{Schensted}; it operates on arbitrary sequences of~$n$
numbers (with equal entries allowed), and it constructs as $P$-symbol a
semistandard (or column-strict) tableau with the same multiset of entries as
the word, while the $Q$-symbol remains a standard tableau. The symmetry that
is lost here is restored in a further generalisation by Knuth \ref{Knuth PJM},
which operates on matrices with natural numbers as entries, and produces pairs
of tableaux which are both semistandard, the multiplicities of their entries
being given by the column sums (for the $P$-symbol) and the row sums (for the
$Q$-symbol) of the matrix. The basic Robinson-Schensted correspondence is
recovered when all row and column sums are equal to~$1$ (the case of
permutation matrices); the generalisation given by Schensted corresponds to
the case where the row sums are~$1$ but columns sums are arbitrary.

While this generalises the correspondence considerably, the algorithm itself
changes only marginally. The case of a matrix with multiple entries in the
same row or column, or entries exceeding~$1$, is handled by operating in the
same way the basic algorithm would for a permutation matrix derived from it by
splitting up rows (each following non-zero entry getting a fresh row below
that of the previous one), multiplexing columns similarly, and replacing
entries $m>1$ by $m\times m$ identity sub-matrices. The Schur function~$s_\\$
is the generating series of the semistandard tableaux of shape~$\\$, so
Knuth's correspondence provides a bijective proof of the Cauchy identity
$\prod_{i,j}{1\over1-X_iY_j}=\sum_\\ s_\\(X)s_\\(Y)$. Knuth also defines a
second correspondence that provides a bijective proof of a ``dual'' identity
$\prod_{i,j}(1+X_iY_j)=\sum_\\s_\\(X)s_{\\\tr}(Y)$. The truncation that has
occurred here of the factors of the left hand side to their terms of
degree~$\leq1$, means that matrix entries are now restricted to the
values~$\{0,1\}$; the transposition of~$\\$ in the second factor on the right
means either that $P$-~and $Q$-symbols have transpose shapes, or (since we
prefer pairs of equal shape) that one is semistandard and the other
transpose-semistandard (row-strict). Like the first one, this second
correspondence can be constructed by first transforming the given matrix to a
permutation matrix and then applying the Schensted algorithm; the only
difference is that columns are multiplexed in the opposite sense: each next
non-zero entry gets a fresh column to the \emph{left} of the previous ones. In
the first correspondence the $P$-~and~$Q$-symbols have symmetric roles, and
transposition of the matrix leads to exchanging them. The second
correspondence lacks such a symmetry, an we shall refer to it as Knuth's
asymmetric correspondence; when not explicitly calling a correspondences
asymmetric, it will be assumed to be a symmetric one.

Fomin has shown in a series of papers \ref{Fomin RSK}--\ref{Fomin Schur} that
by identifying the various tableaux with paths in a graded partially ordered
set or in a directed graph, these correspondences and many other ones can be
described in a general framework that links local correspondences in the poset
or graph to global correspondences involving pairs of paths. This also means
that new correspondences can automatically be defined as soon as a poset or
graph with the required local structure is found. As a consequence of the
generality of these constructions, the terms ``Schensted correspondence'' and
``Knuth correspondence'' now acquire a generic meaning, and the specific
correspondences mentioned above will referred to as the Robinson-Schensted
correspondence and the RSK~correspondence (use this term only for the
symmetric one). Our construction will be an instance of this general
framework, so we shall recall the necessary parts of Fomin's work in detail in
section~\forward{\FominSec}; here we just sketch the outlines. Although the
nature of the poset elements (or vertices of the graph) is not specified, we
shall refer to them as ``shapes''; they are Young diagrams in the cases of the
Robinson-Schensted and RSK~correspondences.

The general constructions build a two-dimensional array of shapes, from which
the $P$-symbol and $Q$-symbol can be read off in the two different directions.
To have an analogue of the Robinson-Schensted correspondence one needs a
graded poset whose connected components are ``$r$-differential''. The most
important requirement for this is that for all shapes~$\mu$ the number of
shapes covering~$\mu$ exceeds the number of shapes covered by~$\mu$ by a fixed
number $r>0$ (the precise requirement is that the commutator of the ``up'' and
``down'' operators for the poset be $r$ times the identity operator). The
Young lattice~$\Y$, consisting of Young diagrams ordered by inclusion, is well
known to be $1$-differential. The $r$-differential property can be ``made
bijective'' by means of an \r-correspondence, which defines for every
shape~$\mu$ a bijection between on one side the set of shapes covering~$\mu$,
and on the other side the union of the set of shapes covered by~$\mu$ and a
set of $r$~extra values. Given such an \r-correspondence, Fomin's construction
will produce a ``Schensted correspondence'' between \r-coloured permutations
of~$n$ (where each term has an additional attribute with $r$ possible values,
whence their number is~$n!r^n$) and pairs of saturated increasing paths of
length~$n$ in the poset with a fixed minimal element as starting point and a
common (but varying) end point. For the Young lattice such paths correspond to
standard Young tableaux. For that case there are two natural choices for a
$1$-correspondence, one of which leads to the usual Schensted correspondence
by row insertion, the other to its transposed variant (using column
insertion).

For ``Knuth correspondences'' the general scheme, which is described
in~\ref{Fomin Schur}, is more complicated. The graded set of shapes used has
more than a poset structure: it is equipped with a directed graph, where edges
may relate shapes any number of levels apart. The entries of the matrices that
form the input of the construction come from a graded but otherwise
unstructured set~$S$. For the Young lattice there is edge from $\mu$ to~$\\$
whenever $\mu/\\$ is a horizontal strip (so that paths correspond to
semistandard tableaux), while $S=\N$. The notion that replaces that of an
\r-correspondence is what we shall call a shape datum; it gives for every pair
$(\mu,\nu)$ of shapes a bijection from shapes~$\kappa$ with edges toward both
$\mu$ and~$\nu$, to pairs $(a,\\)$, where $\\$ is a shape with edges
\emph{from} both $\mu$ and~$\nu$, and~$a\in S$. It must satisfy a
compatibility with the gradings. This implies that when restricted to edges
between shapes at most one level apart, it reduces to an \r-correspondence,
where $r$ is the size of the rank~$1$ subset of~$S$.

The shape datum that matches Knuth's original construction can be found by
considering how that construction deals with a single matrix entry and one
horizontal strip each of the $P$ and $Q$ symbols. Since such a strip is
treated just like a skew standard tableaux corresponding to it (by way of
``standardisation''), the mentioned shape datum is defined by a localised case
of the original Schensted correspondence. Shape data in general however need
not be derived from any Schensted correspondence.

The lattice~$\Y^r$ is an example of an \r-differential poset with $r>1$, and
an \r-correspondence for~$\Y^r$ can be defined by fixing $1$-correspondences
on each of its $r$~factors. Similarly a graph on~$\Y^r$, and a shape datum for
it with graded set $S=\N^r$, can be derived in a component-wise fashion from
the graph on~$\Y$ defined by horizontal strips. These structures lead to
Schensted- and Knuth correspondences that factor into $r$~independent copies
of the original ones, which is not very interesting. Another example is given
by the \r-rim hook lattices, whose elements are (single) Young diagrams, but
whose covering relation, adding an \r-ribbon, relates shapes that are
$r$~levels apart in~$\Y$; the construction of~\ref{Stanton White} provides a
Schensted correspondence this example. However, these lattices are
isomorphic to~$\Y^r$ by means of the so-called \r-quotient map, and the
Stanton-White correspondence thus reduces to the example above.

Yet rim hook lattices can also give rise to Schensted correspondences that do
not decompose into independent ones, by choosing an \r-correspondence not
derived from the \r-quotient map. Indeed, in~\ref{Shimozono White} a more
interesting \r-correspondence for the \r-rim hook lattices is defined, which
unlike the previous one takes the shapes of the ribbons into account. To edges
in an \r-rim hook lattice one may assign a value $h\in\{0,\ldots,r-1\}$,
namely the height of the associated \r-ribbon; with the $r$~extra values
occurring in an \r-correspondence also labelled with that set of values,
Shimozono and White define an \r-correspondence preserving these heights. As a
consequence, the Schensted correspondence obtained respects the ``spin''
statistic on standard \r-ribbon tableaux that gives the sum of the heights of
the ribbons: the sum of the colours of the input permutation determines the
sum of the spins of the output tableaux.

In this paper we define a height respecting shape datum for \r-rim hook
lattices equipped with the graphs defined by the notion of horizontal
\r-ribbon strip (which is essentially a placement of ribbons of the previous
subsection; this notion underlies that of semistandard \r-ribbon tableaux),
with graded set $S=\N^r$. The Knuth correspondence derived from it shares the
``colour-to-spin'' property of the Schensted correspondence of~\ref{Shimozono
White}, for the natural definitions of the respective statistics on matrices
and semistandard \r-ribbon tableaux. Given the way the original Knuth
correspondence is derived from the Robinson-Schensted correspondence, it may
seem a straightforward process to obtain such a shape datum from the
\r-correspondence of Shimozono and White; in any case, this is what we thought
initially. It is not. In fact we were unable to find a shape datum that would
reduce, when appropriately restricted, to that \r-correspondence; this is
essentially for the same reason that we know no bijective proof of our
claim~\alpineclaim. Instead our shape datum reduces to the \emph{transpose}
\r-correspondence (so to speak its column insertion variant, although that
term is not very appropriate in the \r-ribbon case). Since horizontal
\r-ribbon strips (unlike \r-ribbons) lack transposition symmetry, this
distinction is significant.

Even more surprisingly, our correspondence does not reduce for $r=1$ to the
shape datum corresponding to the RSK~correspondence, but to one associated to
the so-called Burge correspondence. These shape data are fundamentally
different, even if the RSK~correspondence and the Burge correspondence are
related (the relation also involves the Sch\"utzenberger involution on
semistandard tableaux). The most crucial observation we had to make in order
to find the shape datum used in our construction, was that although the shape
datum for the Burge correspondence can be defined by iterating Schensted
column-insertion, it has an alternative single-pass description (much like the
bus transport in the previous subsection) that can easily be adapted to the
context of semistandard
\r-ribbon tableaux.

The remainder of this paper is organised as follows. In
section~\forward{\FominSec} we recall in detail Fomin's general framework to
define Schensted and Knuth correspondences. In section~\forward{\KnuthSec} we
describe the shape data for the RSK~correspondence, and for the Burge
correspondence, while also indicating how the global correspondences defined
by them are related. We close that section by giving the most trivial examples
of correspondences with $r>1$, namely those using $\Y^r$. In
section~\forward{\ribbonSec} we recall the definitions involving \r-rim hook
lattices and semistandard \r-ribbon tableaux, and the factoring of many
questions concerning \r-ribbons due to the \r-quotient map; we then discuss
the spin statistic and the \r-correspondence defined by Shimozono and White,
which do not factor in this way. In section~\forward{\mainSec} we present our
main result, the shape datum that leads to a spin preserving Knuth
correspondence from matrices with entries in~$\N^r$ to pairs of semistandard
\r-ribbon tableaux. In section~\forward{\asymSec} we similarly generalise
Knuth's asymmetric correspondence to one from matrices with entries
in~$\{0,1\}^r$ to pairs consisting of a semistandard and a transpose
semistandard \r-ribbon tableau.

Although the presentation of our new constructions is our central goal around
which the paper is organised, much of it can also be considered as an
expository one of the various constructions on which we build forth. For this
reason our pace will often be leisurely and our discussion informal. We have
tried to limit the formal definitions and notations used, and most of those
that are introduced serve a very localised purpose; giving them will therefore
be postponed to the moment they are actually used. We shall also take time to
discuss some matters that are not essential to our constructions, notably the
reasons why certain approaches we tried were not successful.

\subsecno=0

\newsection Review of Fomin's constructions.
\labelsec\FominSec

\subsection Schensted correspondences.
\labelsec\Schenstedsec

Let us recall the framework described in \ref{Fomin duality}--\ref{Fomin
algorithms} for defining Schensted correspondences, while simplifying it
slightly by omitting some generality that is not needed in the current paper.
For a treatment of just the case of the Robinson-Schensted correspondence, we
may refer to~\ref{van Leeuwen festschrift,~\Sec3}. We keep in mind that
special case throughout the discussion, and keep our notation close to what is
customary there.

One starts with graded set~$\Part$ of shapes, whose the rank function
$\Part\to\N$ is written $\\\to|\\|$ and is such that each set
$\Part_i=\setof\\\in\Part:|\\|=i\endset$ is finite. A structure of graded
graph on~$\Part$ is defined by giving a relation contained in
$\Union_{i\in\N}\Part_i\times\Part_{i+1}$, in other words a set of edges
$(\\,\mu)$ with $|\mu|=|\\|+1$. Although two different structures of graded
graph on the same set~$\Part$ are used in the theory (one for each of the
``up'' and ``down'' operators introduced below), these always coincide in the
cases we shall consider, so we use a single symbol~`$\prec$' to denote this
relation. Its reflexive and transitive closure will be denoted by~`$\leq$',
which makes~$\Part$ into a graded poset, from which one can retrieve the
graded graph as its Hasse diagram. The latter need not be connected, and we
shall encounter examples where this is not the case; however, since we shall
be mainly interested in paths in the graph, we loose no generality by
considering one connected component at the time. For $\\\leq\mu$ we define a
path of shape~$\mu/\\$ to be a monotonously rising path
$\\=\\^0\prec\foldwith{\\^#1}\prec\cdots(1..n)=\mu$ in the graph, and we shall
write $|\mu/\\|=n=|\mu|-|\\|$. Due to further requirements each connected
component of the graph will have a unique minimal element, which by a shift in
the grading on that component could be assumed to have rank~$0$. A path of
shape $\mu/\\$ where $\\$ is the minimal element of its component will be
called simply a path of shape~$\mu$.

Let $\Z\Part$ denote the free $\Z$-module on the set~$\Part$ of shapes; one
defines endomorphisms $U,D$ of $\Z\Part$ by their action on basis elements:
$U\:\\\mapsto\sum_{\mu\succ\\}\mu$ and $D\:\\\mapsto\sum_{\mu\prec\\}\mu$
(these are well defined by the finiteness of each~$\Part_i$). The basic
assumption we make of our graded graph is the commutation relation
$$
  D\after U=U\after D+r\Id
\label(\DUcomm)
$$
for some $r\in\Np$, where $\Id$ is the identity operator. It is clear that in
any case both members of this equation preserve the grading. So if one applies
the equation to some~$\\\in\Part_i$, it states two things: taking the
coefficient of~$\\$ itself, one gets
$$
  \Card\setof\mu\in\Part_{i+1}:\mu\succ\\\endset
 =\Card\setof\mu\in\Part_{i-1}:\mu\prec\\\endset+r,
\label(\diageq)
$$
and taking the coefficient of some $\\'\in\Part_i$ with $\\\neq\\'$ one gets
$$
  \Card\setof\mu\in\Part_{i+1}:\mu\succ\\;\mu\succ\\'\endset
 =\Card\setof\mu\in\Part_{i-1}:\mu\prec\\;\mu\prec\\'\endset.
\label(\offdiageq)
$$
The latter condition implies that if there is any path between two shapes,
then there is one that is composed of a monotonous descending path followed by
a monotonous rising path, so in particular each connected component of~$\Part$
has a unique minimal element. Therefore equation (\DUcomm) means that each
such component is an \r-differential poset as defined in~\ref{Stanley
differential posets,~definition~1.1}. One also sees easily that if
condition~(\offdiageq) holds for all $\\,\\'$, then the cardinalities in the
equation cannot exceed~$1$. In case a component of $\Part$ is a lattice, that
condition amounts to the lattice being modular (i.e.,
$|\\\meet\mu|+|\\\join\mu|=|\\|+|\mu|$).

The classical example of a $1$-differential poset is the Young lattice~$\Y$
consisting Young diagrams (finite order ideals of~$\N^2$) ordered by
inclusion. A partition $\\_0\geq\\_1\geq\cdots$ represents a Young diagram
$\\=\setof(i,j)\in\N^2:i<\\_j\endset$, whose elements are referred to (and
displayed) as squares. In~$\Y$, a path
$(0)=\foldwith{\\^#1}\prec\cdots(0..n)=\\$ of shape~$\\$ corresponds to a
standard Young tableau of shape~$\\$, which can be displayed by filling the
squares of~$\\$ with the numbers in $\set{n}=\{0,\ldots,n-1\}$ in such a way
that for each~$i$ the set of entries in the diagram~$\\^i$ is $\set{i}$. Thus
the standard Young tableau displayed as
$$
  {\smallsquares\Young(0,2,3,5|1,6|4)}
\ttext{corresponds to the path}
\circ\prec{\smallsquares\Young()\prec\Young(|)\prec\Young(,|)\prec\Young(,,|)
     \prec\Young(,,||)\prec\Young(,,,||)\prec\Young(,,,|,|)}.
$$
For our purposes we shall consider standard Young tableaux of shape~$\\$
just to \emph{be} paths of shape~$\\$.

Let us check that $\Y$ is a $1$-differential poset. Since $\Y$ is a modular
(even distributive) lattice (with $\\\meet\mu=\\\thru\mu$ and
$\\\join\mu=\\\union\mu$), it suffices to check condition~(\diageq), which
states that for any Young diagram there is one more square that can be added
to it than there are squares that can be removed from it. Although this is
quite easy to see, we give a formal argument, since variations of it will be
used later to prove less obvious statements. A Young diagram~$\\$ is
completely determined by its edge sequence~$\edge\\$, a doubly infinite bit
string describing (as in \Sec\enumsec) the segments of the path that forms the
edge of~$\\$, from bottom left to top right (for a formal definition
of~$\edge\\$ see~\ref{edge sequences}; it is essentially the (Com\'et) code
$C_\\$ of~$\\$ used in \ref{Stanley EC2,~Exercise~7.59}, but with all bits
complemented). The bits eventually become~$1$ to the left (for the rows of
length~$0$) and~$0$ to the right (for the columns of height~$0$). Occurrences
of a substring~``$10$'' correspond to squares that can be added to the
diagram, causing the substring to be replaced by~``$01$'', occurrences of
which therefore correspond to squares that can be removed. We illustrate
this for the diagram $\\$ of~$(6,4,3,3,1)$, and the bracketed substring
of $\edge\\=(\cdots1101001[10]100100\cdots)$.
\bigdisplay
\dimen2=5mm
\dimen0=-0.35\dimen2
\dimen1=0.075\dimen2
\def~#1{\advance\dimen#1 by 0.5\dimen2
        \raise\dimen1\rlap{\kern\dimen0$\scriptstyle#1$\hss}%
        \advance\dimen#1 by 0.5\dimen2
       }
\epsfxsize=8\dimen2
\vcenter{\hbox{~1~1~0~1~0~0~1~1~0~1~0~0~1~0~0%
               \advance\dimen0 by -4.65\dimen2
               \advance\dimen1 by -3.4\dimen2
		~0~1%
               \epsfbox{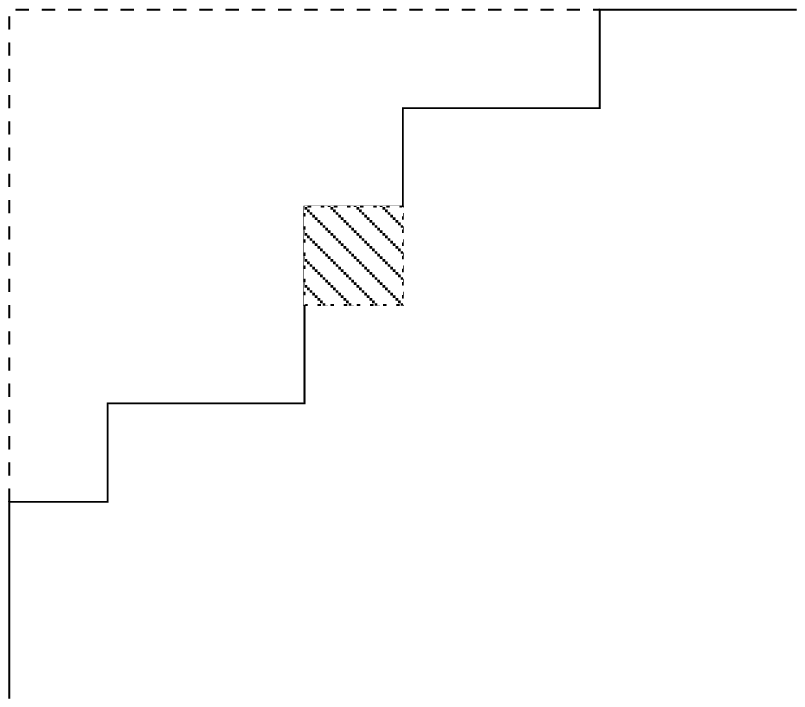}%
             }}
$$
In any word over the alphabet $\{0,1\}$ the occurrences of ``$10$'' and
``$01$'' are perfectly interleaved from left to right. The limiting behaviour
of letters in both directions ensures that the set of these occurrences is
finite and non-empty, and that the first and the last ones are occurrences of
``$10$''. Therefore the number of occurrences of ``$10$'' exceeds that of
occurrences of ``$01$'' by~one.

As is shown in \ref{Stanley differential posets}, many enumerative properties
can be derived uniformly for all \r-differential posets, i.e., they depend
only on the identity~(\DUcomm). In \ref{Fomin algorithms} bijections
corresponding to such identities are constructed, but this requires some
additional ingredients (as could be expected, since otherwise there is nothing
to discriminate between the elements of~$\Part$), namely a family of
bijections that correspond to local instances of~(\DUcomm), i.e., linking the
sets counted by the members of equations (\diageq) and~(\offdiageq). In fact
one needs not bother about~(\offdiageq), since we have observed that those
sets are either empty or singletons, which leaves no choice for a bijection.
It is not obvious which set should be counted by the final term~`$r$'
in~(\diageq); any $r$~element set disjoint from~$\Part$ would do. Recall that
we denoted by~$\set{r}$ the $r$~element order ideal
$\setof{i\in\N}:i<r\endset=\{0,\ldots,r-1\}$ of~$\N$, which we take as our
standard \r-element set. We could let the mentioned term~`$r$' correspond
to~$\set{r}$, but it will be clearer, and more flexible for generalising
later, to take a set $e_\set{r}=\setof e_i:i\in\set{r}\endset$ of $r$~symbols
that are reserved to serve as exceptional values. We then define an
\r-correspondence on~$\Part$ to consist of a family of bijections
$$
  b_\\\:\setof\mu\in\Part:\mu\succ\\\endset
     \to\setof\mu\in\Part:\mu\prec\\\endset\union e_\set{r}
\ttex{for $\\\in\Part$.}
\label(\rcorrdef)
$$
(In fact \r-correspondences are defined in \ref{Fomin algorithms} as
bijections between edges rather than between vertices of the graph, but in the
situation we are considering, the simpler notion above gives the same
information). Assuming that~(\DUcomm) holds, an \r-correspondence always
exists, although it might be hard to specify; usually however, the proof of
(\DUcomm) suggests a natural choice of an \r-correspondence. For $\Y$, we see
that there are in fact two natural choices: one may either associate to any
occurrence of~``$10$'' in~$\edge\\$ the occurrence of~``$01$'' to its left,
except that the leftmost occurrence of~``$10$'' is sent to the exceptional
value~$e_0$, or one can do the same thing with ``left'' replaced by ``right''.
Due to the symmetry of~$\Y$ with respect to transposition, there is no
fundamental difference between the two choices, but the choice that will lead
to the Schensted correspondence in its usual form (i.e., using row insertion)
is to move to the right: each square that can be added will then correspond to
the square that can be removed from the row directly above it, with the square
that can be added to the topmost row corresponding to~$e_0$.

Once an \r-correspondence has been fixed, it may be used to define bijective
correspondences that match identities derived from equation~(\DUcomm).
Although there are many of these, we focus on the following one:
$$
  \scal< \epsilon | D^n(U^n(\epsilon)) > = n!r^n
  \ttex{for a minimal element $\epsilon$ of~$\Part$, and $n\in\N$,}
\label(\RSid)
$$
where the scalar product is the canonical one in~$\Z\Part$, so that the left
hand side gives the coefficient of~$\epsilon$ in $D^n(U^n(\epsilon))$. The
identity can be derived by showing that when~(\DUcomm) is used repeatedly to
rewrite $D^n\after U^n$ as a sum of terms of the form $U^i\after D^i$ with
$0\leq i\leq n$, then $U^0\after D^0=\Id$ occurs with coefficient~$n!r^n$.
Since the operators $D$ and~$U$ are adjoint with respect to the scalar
product, the left hand side of (\RSid) can also be written as
$\scal<U^n(\epsilon)|U^n(\epsilon)>$. Every $\\\in\Part_{|\epsilon|+n}$ with
$\epsilon\leq\\$ occurs in $U^n(\epsilon)$, with as multiplicity the number of
paths of shape $\\$, so (\RSid) states that the sum of the squares of these
numbers equals~$n!r^n$. For $\Part=\Y$, paths of shape $\\$ are standard
tableaux of shape~$\\$, and for the value $r=1$ applicable to this case the
number $n!r^n=n!$ of course counts the permutations of~$n$, so here
equation~(\RSid) gives the enumerative identity that the Robinson-Schensted
correspondence proves. Any correspondence similarly proving~(\RSid) for
some $r$-differential poset, with the right hand side $n!r^n$ counting the
\r-coloured permutations of~$n$, can therefore be considered to be a
generalisation of the Robinson-Schensted correspondence.

Given an \r-correspondence on~$\Part$, one can obtain such a correspondence by
defining, for any minimal element~$\epsilon$ of~$\Part$ and $n\in\N$, an
intermediate set of ``growth diagrams'' that is in bijection both with the set
of pairs of paths of some common shape $\\\in\Part_{|\epsilon|+n}$ with
$\epsilon\leq\\$ and with the set of \r-coloured permutations of~$n$. These
bijections are just projections that extract partial information from a growth
diagram; their bijectivity therefore means that growth diagrams can always be
uniquely reconstructed from such information. We formally define \r-coloured
permutations of~$n$ to be matrices $(A_{i,j})_{i,j\in\set{n}}$ with entries
$A_{i,j}\in\{0\}\union e_\set r$ such that $|A_{i,j}|_{i,j\in\set{n}}$ is a
permutation matrix, where we set $|e_i|=1$ for all $i\in\set{r}$. Specifying
an~\r-coloured permutation of~$n$ is equivalent to giving the
permutation~$\sigma$ of~$n$ such that $|A_{i,j}|=1$ whenever $\sigma(i)=j$,
and the $n$~values $A_{i,\sigma(i)}\in e_\set{r}$ ($i\in\set{n}$), which
explains the name.

\proclaim Definition. \growthdef
Let $\Part$ be a graded poset equipped with an~\r-correspondence
$\setof{b_\\}:\\\in\Part\endset$ as defined by equation~(\rcorrdef), and let a
minimal element~$\epsilon$ of~$\Part$ and $n\in\N$ be fixed. A \newterm{growth
diagram}, or \newterm{Schensted-growth}, relative to these data consists of a
pair of maps, the first one mapping $\set{n+1}^2\to\Part$ and written
$(k,l)\mapsto\\^{(k,l)}$, the second one mapping $\set{n}^2\to\{0\}\union
e_\set r$ and written $(k,l)\mapsto A_{k,l}$, that satisfy the following
conditions for all~$k,l\in\set{n}$:
\itemno=-1
\statitem $\\^{(0,0)}=\epsilon$;
\statitem
   $\\^{(k,l)}\preceq\\^{(k,l+1)}$ and $\\^{(k,l)}\preceq\\^{(k+1,l)}$;
\statitem $\\^{(n,l)}\prec\\^{(n,l+1)}$ and $\\^{(k,n)}\prec\\^{(k+1,n)}$;
\statitem
   $|\\^{(k,l)}|=|\epsilon|+\sum_{(i,j)\in\set{k}\times\set{l}}|A_{i,j}|$;
\statitem
   If $\\^{(k,l+1)}=\\^{(k+1,l)}\prec\\^{(k+1,l+1)}$ then, putting
   $\mu=\\^{(k,l+1)}$, $\kappa=\\^{(k+1,l+1)}$ and $v=b_\mu(\kappa)$, one has
   either $\\^{(k,l)}=v$ (if $v\in\Part$), or $A_{k,l}=v$ (if $v\in
   e_\set{r}$).

The alternative term ``Schensted-growth'' for growth diagram will be used
later to distinguish this notion from a similar one used in relation to Knuth
correspondences. It follows from \statitemnr{2}~and~\statitemnr{3} that
$(A_{i,j})_{i,j\in\set{n}}$ is an~\r-coloured permutation of~$n$. One also has
$\\^{(i,0)}=\epsilon=\\^{(0,i)}$ for $i\in\set{n+1}$, so that
$P=(\\^{(n,0)}\prec\foldwith{\\^{(n,#1)}}\prec\cdots(1..n))$ and
$Q=(\\^{(0,n)}\prec\foldwith{\\^{(#1,n)}}\prec\cdots(1..n))$ are paths of
shape~$\\^{(n,n)}$; this defines the two mentioned projections. For the
permutation $\sigma$ of~$n$ corresponding to $|A_{i,j}|_{i,j\in\set{n}}$, one
has $\\^{(k,l)}=\\^{(k+1,l)}$ in condition~\statitemnr{1} if and only if
$\sigma(k)\geq l$, and $\\^{(k,l)}=\\^{(k,l+1)}$ if and only if
$\sigma\inv(l)\geq k$. In particular $|A_{k,l}|=1$ implies that
$\\^{(k,l)}=\\^{(k+1,l)}=\\^{(k,l+1)}$, so in condition~\statitemnr{4} one has
$A_{k,l}=0$ in the case $v\in\Part$ (since $b_\mu(\kappa)\neq\mu$), while
$\\^{(k,l)}=\\^{(k+1,l)}=\\^{(k,l+1)}$ in the case $v\in e_\set r$.

It is convenient to display a growth diagram by attaching the shape
$\\^{(k,l)}$ to the point $(k,l)$ of a grid (with $k$ increasing downwards and
$l$ increasing to the right), and writing the matrix entry $A_{k,l}$ into the
square of that grid with corners $(k,l)$, $(k,l+1)$, $(k+1,l)$,
and~$(k+1,l+1)$. Figure~1 thus illustrates a growth diagram for $\Part=\Y$
with the mentioned left-to-right $1$-correspondence; non-zero matrix entries,
which due to $r=1$ are necessarily equal to~$e_0$, are indicated by~`$\star$'.
\midinsert
$$\vcenter{
\catcode`@=11
\def\d{\vtop\bgroup\smallskip \futurelet\n\D}\let\f=\hfil
\def\D{\?\if *\n 
         \then \vbox to\z@{\vss\llap{$\star$\kern4pt}\kern 4pt}%
               \def\next*{\futurelet\n\E}\next
         \else \E
	 \fi}
\def\E{\?\ifx\n\f \then\smallskip\egroup \else\afterassignment\e\count@ \fi}
\def\e{\hbox{\loop\unhcopy\QEDbox\advance\count@\m@ne
             \ifnum\count@>\z@ \repeat}%
       \nointerlineskip\futurelet\n\E}
\def\g#1{\omit\hfil#1\hfil}
\def\h#1 #2{\vtop to 18pt
            {\smallskip\kern-\dimen@\kern#1\dimen@\hbox{#2}\vss}%
            \quad\vrule}
\def\z{\omit\vtop{\smallskip\hbox{$\circ$}}\hfil}
\offinterlineskip
\everycr{\noalign{\nobreak}}
\tabskip 1em plus 3pc
\halign{\h# &&\d#\f\cr
\omit\hfil\strut\vrule
    & \g0 & \g1 & \g2 & \g3\quad & \g4 &  \g5  &  \g6  &  \g7  &  \g8  \cr
\noalign{\hrule}
 1 0 & \z & \z  & \z  & \z  & \z  &  \z   &  \z   &  \z   &  \z   \cr
 1 1 & \z & \z  & \z  & \z  & \z  &*  1   &   1   &   1   &   1   \cr
 2 2 & \z & \z  &* 1  &  1  &  1  &  1 1  &  1 1  &  1 1  &  1 1  \cr
 2 3 & \z & \z  &  1  &  1  &  1  &  1 1  &  1 1  &* 2 1  &  2 1  \cr
 3 4 & \z &* 1  & 1 1 & 1 1 & 1 1 & 1 1 1 & 1 1 1 & 2 1 1 & 2 1 1 \cr
 3 5 & \z &  1  & 1 1 &*2 1 & 2 1 & 2 1 1 & 2 1 1 & 2 2 1 & 2 2 1 \cr
 3 6 & \z &  1  & 1 1 & 2 1 & 2 1 & 2 1 1 & 2 1 1 & 2 2 1 &*3 2 1 \cr
 3 7 & \z &  1  & 1 1 & 2 1 & 2 1 & 2 1 1 &*3 1 1 & 3 2 1 & 3 3 1 \cr
 3 8 & \z &  1  & 1 1 & 2 1 &*3 1 & 3 1 1 & 3 2 1 & 3 2 2 & 3 3 2 \cr}
}
$$
\centerline{{\bf Figure 1.} Schensted-growth for $\Y$ with
 $\sigma={0~1~2~3~4~5~6~7\choose4~1~6~0~2~7~5~3}$, 
 $P={\smallsquares\Young(0,2,3|1,5,7|4,6)}$, and
 $Q={\smallsquares\Young(0,2,5|1,4,6|3,7)}$.}
\endinsert

The pair $(P,Q)$ gives the values of $\\^{(i,j)}$ for
$(i,j)\in\set{n+1}^2\setminus\set{n}^2$, so reconstruction of a growth diagram
by decreasing values of $k$ and~$l$, from arbitrary $(P,Q)$, will be possible
if for any $k,l\in\set{n}$ one can uniquely determine the shape
$\\=\\^{(k,l)}$ and the matrix entry~$A_{k,l}$ once the shapes
$\mu=\\^{(k,l+1)}$, $\nu=\\^{(k+1,l)}$, and~$\kappa=\\^{(k+1,l+1)}$ are known.
And indeed this is the case: if $\kappa$ is equal to $\mu$ or to $\nu$, then
$\\^{(k,l)}$ will be equal to the other, and otherwise if $\mu$ and $\nu$ are
distinct, then $\\^{(k,l)}$ will be the unique element covered by both (which
exists because both are covered by~$\kappa$); in both these cases $A_{k,l}=0$.
The remaining case $\mu=\nu\neq\kappa$ is the one handled by
condition~\statitemnr{4}; there, as we saw, both $\\$ and~$A_{k,l}$ are always
determined. From this description one sees that the differential version
$|\kappa|-|\mu|-|\nu|+|\\|=|A_{k,l}|$ of condition~\statitemnr{3} holds in all
cases; the verification of the remaining requirements is now easy. A growth
diagram can be similarly reconstructed from the matrix
$(A_{i,j})_{i,j\in\set{n}}$, by increasing values of $k$ and~$l$: one
determines $\\^{(k+1,l+1)}$ once $\\^{(k,l)}$, $\\^{(k,l+1)}$, $\\^{(k+1,l)}$
(and of course $A_{k,l}$) are known; to prove this,
  similar cases as above can
be distinguished, with the case $|A_{k,l}|=1$ being singled out first.

For $\Part=\Y$ with the chosen $1$-correspondence, the resulting bijection
between permutations and pairs of standard Young tableaux is the classical
Robinson-Schensted correspondence. The Schensted insertion of the permutation
entry $\sigma(i)$ corresponds to the computation, given the line
$\foldwith{\\^{(i,#1)}}\preceq\cdots(0..n)$ of the growth diagram and the
entry $A_{i,\sigma(i)}=e_0$, of the next line
$\foldwith{\\^{(i+1,#1)}}\preceq\cdots(0..n)$.

The method of reconstructing a growth diagram from $(P,Q)$ or from
$(A_{i,j})_{i,j\in\set{n}}$ can be viewed as a data-flow network consisting of
$n^2$ copies of the \r-correspondence respectively of its inverse: for every
pair $k,l\in\set{n}$, a copy of the \r-correspondence links the values
$A_{k,l}$, $\\^{(k,l)}$, $\\^{(k+1,l)}$, $\\^{(k,l+1)}$, and~$\\^{(k+1,l+1)}$.
In fact we had to complement the \r-correspondence with several cases (where
no choice was involved) to determine all necessary values, but in the full
generality of \ref{Fomin algorithms} these cases are already incorporated into
the \r-correspondence itself. Moreover that \r-correspondence operates on the
edges rather than on the vertices of the square $\\^{(k,l)}$, $\\^{(k+1,l)}$,
$\\^{(k,l+1)}$, and~$\\^{(k+1,l+1)}$, with each of the two output edges being
one of the inputs to a unique other copy of the \r-correspondence (unless it
is part of the final output), so that the data-flow nature is even more
clearly visible there.

One may observe that apart from extraction of the pair $(P,Q)$ of paths
in~$\Part$, or of~$A$, the \r-coloured permutation of~$n$, other projections
allow unique reconstruction of a growth diagram as well: it suffices to know
$\\^{(i,j)}$ for all $(i,j)$ on some lattice path from $(n,0)$ to $(0,n)$, and
all matrix entries $A_{i,j}$ that lie below and to the right of that path.
Starting with the shapes along the bottom and right edges of the diagram as
given by~$(P,Q)$, one can therefore replace this information by the shapes
along paths that are gradually modified to eventually become the path along
the left at top edges, while retaining all matrix entries that the path has
moved across; then at each moment one has complete information to determine
the entire growth diagram. This procedure may be interpreted as describing the
process of rewriting the left hand side of~(\RSid) to its right hand side: it
provides a way to trace, and thus obtain a matching between, the contributions
to these two expressions and to all intermediate ones.

\subsection Knuth correspondences.
\labelsec\FominKnuthsec

We shall now consider the generalisation of this construction to ``Knuth
correspondences'', which is described in \ref{Fomin Schur}. The fundamental
difference with Schensted correspondences is that matrix entries will be
chosen independently, without the restriction that they should give rise to a
permutation matrix. Matrix entries will be chosen from some graded set~$S$,
whose rank function $S\to\N$ written $x\to|x|$ is such that each set
$S_i=\setof x\in S:|x|=i\endset$ is finite, and $S_0=\{0\}$. We shall define
growth diagrams similar to those defined in above, for which in particular
condition~\growthdef\statitemnr{4} will be required. From the independence of
matrix entries it then follows that one has to allow arbitrarily large
differences of rank between shapes on adjacent grid points (at least, growing
with~$n$), even if all non-zero matrix entries should have rank~$1$. Therefore
one needs to equip~$\Part$ with a relation that unlike `$\preceq$' can hold
between shapes any number of levels apart. We shall denote this relation,
which in general will not be transitive, by~`$\leqh$', with $\mu\geqh\\$
meaning the same as $\\\leqh\mu$. At this point we make no assumption about
this relation other than that $\\\leqh\mu$ and $|\\|=|\mu|$ imply $\\=\mu$.
When $\\\leqh\mu$ we shall as before write $|\mu/\\|=|\mu|-|\\|$.

The adaptation of the rules governing growth diagrams is quite
straightforward. Each matrix entry $a=A_{i,j}\in S$ is associated to a square
of the grid, the corners of which have shapes attached to them, which as
before we shall designate by $\\=\\^{(i,j)}$, $\mu=\\^{(i,j+1)}$,
$\nu=\\^{(i+1,j)}$, and $\kappa=\\^{(i+1,j+1)}$; these satisfy
$\\\leqh\mu\leqh\kappa$ and $\\\leqh\nu\leqh\kappa$. From
condition~\growthdef\statitemnr{4} it follows moreover that
$$
  |\kappa|-|\nu|-|\mu|+|\\|=|a|.
\label(\rankeq)
$$
As in the case of Schensted correspondences, we wish that given $\mu$, $\nu$,
and~$\kappa$, the shape~$\\$ and the entry~$a$ can be uniquely determined, and
that similarly given $\\$, $\mu$, $\nu$ and~$a$, the shape~$\kappa$ can be
determined. In other words, for every pair $\mu,\nu\in\Part$ there should be a
bijection between the sets
$\setof\kappa\in\Part:\kappa\geqh\mu;\kappa\geqh\nu\endset$ and
$\setof(a,\\)\in S\times\Part:\\\leqh\mu;\\\leqh\nu\endset$ such that
equation~(\rankeq) holds whenever $(a,\\)$ corresponds to~$\kappa$. Such a
family of bijections will be a central notion of this paper, so we state a
formal definition.

\proclaim Definition. \sddef
Let $\Part$ be a graded set equipped with a relation~`$\leqh$', and $S$ a
graded set. A \newterm{shape datum} for $(\Part,\leqh,S)$ consists of a family
$(b_{\mu,\nu})_{\mu,\nu\in\Part}$ of bijections
$$
  b_{\mu,\nu}\:
 \setof\kappa\in\Part:\mu\leqh\kappa\geqh\nu\endset
\to
 \setof(a,\\)\in S\times\Part:\mu\geqh\\\leqh\nu\endset,
\label(\bmunueq)
$$
such that $|\kappa|-|\nu|-|\mu|+|\\|=|a|$ holds whenever
$(a,\\)=b_{\mu,\nu}(\kappa)$.

In order for a shape datum to exist, $(\Part,\leqh,S)$ must satisfy an
enumerative identity that is most easily expressed using formal power series.
We shall use the ring $\End(\Z\Part)[[X,Y]]$ of formal power series in
indeterminates $X,Y$ over the ring of endomorphisms of~$\Z\Part$. Such a power
series can be applied to an element of~$\Z\Part$ to give a formal sum of
monomials in $X$ and~$Y$ with coefficients in~$\Z\Part$, i.e., an element of
the $\Z[[X,Y]]$-module $\Z\Part[[X,Y]]$. This action on elements of~$\Z\Part$
can be extended without problem to an action on power series, in other words
we can interpret elements of $\End(\Z\Part)[[X,Y]]$ as endomorphisms
of~$\Z\Part[[X,Y]]$. We define two such elements
$U_X,D_Y\in\End(\Z\Part)[[X,Y]]$ in terms of their action, by requiring
$$
  U_X(\\)=\sum_{\mu\geqh\\}X^{|\mu/\\|}\mu
\ttext{and}
  D_Y(\\)=\sum_{\mu\leqh\\}Y^{|\\/\mu|}\mu
\ttex{for all~$\\\in\Part$.}
\nn
$$
Let $F_S(T)\in\Z[[T]]$ be the rank generating series $\sum_{a\in
S}T^{|a|}=\sum_{i\in\N}c_iT^i$ of~$S$, where $c_i=\Card{S_i}$. Then the
enumerative requirement for the existence of a shape datum will be
$$
  D_Y\after U_X=(U_X\after D_Y)\(F_S(XY)\).
\label(\sdid)
$$
If one writes $U_X=\sum_{i\in\N}U_iX^i$ and $D_Y=\sum_{i\in\N}D_iY^i$, then
taking the coefficient of $X^iY^j$ in~(\sdid) gives
$$
  D_j\after U_i=\sum_{k\leq\min\{i,j\}}c_k(U_{i-k}\after D_{j-k}).
\label(\DjUiid)
$$
When both members of this equation are applied to~$\mu$, the values obtained
are of rank $|\mu|+i-j$, and if one takes the coefficient of some
$\nu\in\Part_{|\mu|+i-j}$, the resulting numbers count the subsets of the
domain and codomain of~$b_{\mu,\nu}$ for which $|\kappa|=|\mu|+i$ respectively
$|\\|=|\mu|-j+|a|$; these are subsets that should be in correspondence
under~$b_{\mu,\nu}$ if it is to satisfy~(\rankeq). Taking $i=j=1$ in~(\DjUiid)
gives $D_1\after U_1=U_1\after D_1+c_1\Id$, which shows that the current
situation extends the one considered earlier: if one defines $\\\prec\mu$ when
$\\\leqh\mu$ and $|\\|+1=|\mu|$, then (\diageq)~and~(\offdiageq) will hold,
with $r=c_1$. There is however no easy way to extend the structure
$(\Part,\prec)$ to $(\Part,\leqh,S)$; indeed it is not obvious why any
combinatorial structures satisfying~(\sdid) should exist at all. Yet several
examples are given in~\ref{Fomin Schur}, among which the example corresponding
to the RSK~correspondence (described in the next section), which is the
one that will concern us most.

Once one has defined a shape datum, the construction of a global ``Knuth''
correspondence from it is straightforward, by analogy to the construction of
Schensted correspondences. In fact, matters are slightly simpler since there
is no need to explicitly distinguish cases according to the differences of
rank between shapes associated to neighbouring points of the grid, as all
possible cases are already catered for by the shape datum itself. Calling
$\epsilon\in\Part$ minimal when $\\\leqh\epsilon$ implies $\\=\epsilon$, we
have the following notion of growth diagram.

\proclaim Definition. \Knuthgrowthdef
Let $(b_{\mu,\nu})_{\mu,\nu\in\Part}$ be a shape datum for $(\Part,\leqh,S)$
as defined in~\sddef, and let a minimal element~$\epsilon$ of~$\Part$, and
$n,m\in\N$ be fixed. A growth diagram, or Knuth-growth, relative to these data
consists of a pair of maps, the first one mapping
$\set{m+1}\times\set{n+1}\to\Part$ and written $(k,l)\mapsto\\^{(k,l)}$, the
second one mapping $\set{m}\times\set{n}\to S$ and written $(k,l)\mapsto
A_{k,l}$, that satisfy the following conditions, for all $k\in\set{m}$ and
$l\in\set{n}$:
\itemno=-1
\statitem $\\^{(0,n)}=\epsilon=\\^{(m,0)}$;
\statitem
   $\\^{(k,l+1)}\leqh\\^{(k+1,l+1)}$, and $\\^{(k+1,l)}\leqh\\^{(k+1,l+1)}$;
\statitem
   With $\\=\\^{(k,l)}$,  $\mu=\\^{(k,l+1)}$,  $\nu=\\^{(k+1,l)}$, and
   $\kappa=\\^{(k+1,l+1)}$, one has $(A_{k,l},\\)=b_{\mu,\nu}(\kappa)$.

Note that there is no longer a need to require explicitly that
$|\\^{(k,l)}|=|\epsilon|+\sum_{(i,j)\in\set{k}\times\set{l}}|A_{i,j}|$, as
this follows easily from the fact that $b_{\mu,\nu}$ respects
equation~(\rankeq). We shall again consider two projections from the set of
growth diagrams: one extracting the matrix
$A=(A_{k,l})_{k\in\set{m},l\in\set{n}}$, the other extracting the pair of
paths $P=(\\^{(m,0)}\leqh\foldwith{\\^{(m,#1)}}\leqh\cdots(1..n))$ and
$Q=(\\^{(0,n)}\leqh\foldwith{\\^{(#1,n)}}\leqh\cdots(1..m))$. In the latter,
$P$ and~$Q$ will be called paths of shape $\\^{(m,n)}$ (or more explicitly of
shape $\\^{(m,n)}/\epsilon$) in~$(\Part,\leqh)$. Besides this (final) shape, a
property of paths describing the successive ranks of the intermediate shapes
is important: we define the \newterm{weight} of a path
$p=(\foldwith{\\^#1}\leqh\cdots(0..n))$ to be the vector
$\wt(p)=(|\\^{i+1}/\\^i|)_{i\in\set{n}}$. One has
$\wt(P)=\(\sum_{i\in\set{m}}|A_{i,j}|\){}_{j\in\set{n}}$ and
$\wt(Q)=\(\sum_{j\in\set{n}}|A_{i,j}|\){}_{i\in\set{m}}$, in other words, the
weights of the paths $P,Q$ are the vectors of column sums of ranks and of row
sums of ranks, respectively, of the matrix~$A$. Both projections can be seen
to be bijections, respectively to the set of all matrices with coefficients
in~$S$ and to the set of all pairs of paths in~$(\Part,\leqh)$ of equal shape,
in the same way as for the Schensted-growths. The correspondence between the
matrices~$A$ and pairs $(P,Q)$ of paths, defined by composing one projection
with the inverse of the other, is the Knuth correspondence
for~$(\Part,\leqh,S)$ associated to the given shape datum.

Like Schensted correspondences, Knuth correspondences imply enumerative
identities, but since the sets linked by any bijection~$b_{\mu,\nu}$ are
infinite, these follow from the mentioned matching of the weights of~$(P,Q)$
and of~$A$. One uses two sets of indeterminates $X_\set{n}=\setof
X_j:j\in\set{n}\endset$ and $Y_\set{m}=\setof Y_i:i\in\set{m}\endset$ to
record the weights of $P$ and~$Q$ respectively: the contribution of each pair
$(P,Q)$ will be given by the monomial $X^{\wt(P)}Y^{\wt(Q)}$, where as usual
$X^{(\li(a0..n-1))}$ abbreviates $\foldwith{X_#1^{a_#1}}{}\cdots(0..n-1)$. By
what was observed above, that monomial can also be expressed in terms of the
matrix~$A$ corresponding to~$(P,Q)$ as
$\prod_{i\in\set{m},j\in\set{n}}(X_jY_i)^{|A_{i,j}|}$. Summation of that
monomial over all $m\times n$ matrices~$A$ with entries in~$S$ gives
$\prod_{i\in\set{n},j\in\set{m}}F_S(X_iY_j)$ (note that we disentangle the
indices, the mixing of which was caused by the unfortunate but conventional
choice of reading off~$P$ from the growth diagram by varying the \emph{second}
index). If one defines the generating series
$G_\\(X_\set{n})=\sum_PX^{\wt(P)}$, where the sum is over all paths
$P=(\epsilon=\foldwith{\\^#1}\leqh\cdots(0..n)=\\)$, then one finds the
identity
$$
  \prod_{i\in\set{n},j\in\set{m}}F_S(X_iY_j)
 =\sum_{\\\in\Part}G_\\(X_\set{n})G_\\(Y_\set{m})
\label(\CauchyKnuthid)
$$
Noting that
$G_\\(X_\set{n})=\scal<\\|\foldwith{U_{X_#1}}{}\cdots(n-1..0)(\epsilon)>$,
this equation can be seen to follow directly from~(\sdid):
$$
\eqalign
{ \smash{\sum_{\\\in\Part}G_\\(X_\set{n})G_\\(Y_\set{m})}
 &=\scal<\foldwith{U_{Y_#1}}{}\cdots(m-1..0)(\epsilon)
        |\foldwith{U_{X_#1}}{}\cdots(n-1..0)(\epsilon)>\cr
 &=\scal<\epsilon|\foldwith{D_{Y_#1}}{}\cdots(0..m-1)
                  \foldwith{U_{X_#1}}{}\cdots(n-1..0)(\epsilon)>\cr
 &=\scal<\epsilon|\foldwith{U_{X_#1}}{}\cdots(n-1..0)
                  \foldwith{D_{Y_#1}}{}\cdots(0..m-1)(\epsilon)>
   \smash{\prod_{i\in\set{n},j\in\set{m}}F_S(X_iY_j)}\cr
 &=\prod_{i\in\set{n},j\in\set{m}}F_S(X_iY_j),\cr
}
$$
where the last step follows from the fact that $\epsilon$ is minimal. This
derivation demonstrates (again) the straightforwardness of deriving the global
Knuth correspondence from the shape datum.

\newsection Examples of shape data.
\labelsec\KnuthSec

\subsection The shape datum for the RSK~correspondence.

Let us now consider the RSK~correspondence, and its shape datum. In this case
$S=\N$, $\Part=\Y$, and the relation $\\\leqh\mu$ means that $\mu/\\$ is a
horizontal strip: $\\\subset\mu$, and the skew diagram of~$\mu/\\$ (the
set-theoretic difference of diagrams $\mu\setminus\\$) has at most one square
in any column. When $\\\leqh\mu$, a specific path in~$(\Part,\prec)$ (such
paths are also called skew standard tableaux) of shape $\mu/\\$ is defined by
requiring that in the successive diagrams along the path from~$\\$ to~$\mu$,
the squares are added in left to right order, i.e., by (strictly) increasing
column number. By interpolating in this way each of its horizontal strips
(while eliminating any strips of size~$0$), a path in~$(\Part,\leqh)$ can be
transformed into a (skew) standard tableau of the same shape, called its
standardisation. The RSK~correspondence, in the ``extraction'' direction
(determining the matrix~$A$ from a pair $(P,Q)$ of equal shaped paths
in~$(\Part,\leqh)$), can now be informally described as follows. First the
standardisations of $P$ and~$Q$ are determined, and then the Schensted-growth
corresponding to this pair of standard Young tableaux is constructed. Due to
the initial interpolation, this growth diagram is defined on a grid that is
too extensive for the matrix~$A$ one wishes to find; therefore the grid is now
reduced by forgetting its points on horizontal or vertical lines that were
introduced by the interpolation, so that the paths remaining along the bottom
and right edges are the original $P$ and~$Q$, before standardisation. The
shapes associated to the remaining grid points will be the ones of a
Knuth-growth; the matrix entry to be associated with a grid square whose
corners carry shapes $\\~\mu\choose\nu~\kappa$ is $|\kappa|-|\nu|-|\mu|+|\\|$,
which is the sum over the corresponding rectangular area of the
Schensted-growth of the associated permutation matrix entries. Proving that
this description defines (a~Knuth-growth for) a shape datum for~$(\Y,\leqh)$
essentially amounts to showing that the shapes associated to adjacent grid
points after reduction differ by horizontal strips; this can be deduced easily
from a consideration of local portions of the Schensted-growth consisting of
just two horizontally or vertically adjacent squares of the grid for that
growth diagram.

Let us demonstrate the construction by a concrete example, which we take
from Knuth's original paper \ref{Knuth PJM}. The semistandard tableaux can be
readily interpreted as paths, so for instance
$$
 {\smallsquares\Young(1,1,1,2,4,7|2,3,3,5|3,4,6,6|6)}
 \text{represents}
\left(
 \circ\leqh{\smallsquares\smallsquares\Young(,,)}
      \leqh{\smallsquares\smallsquares\Young(,,,|)}
      \leqh{\smallsquares\smallsquares\Young(,,,|,,|)}
      \leqh{\smallsquares\smallsquares\Young(,,,,|,,|,)}
      \leqh{\smallsquares\smallsquares\Young(,,,,|,,,|,)}
      \leqh{\smallsquares\smallsquares\Young(,,,,|,,,|,,,|)}
      \leqh{\smallsquares\smallsquares\Young(,,,,,|,,,|,,,|)}
\right).
$$
The semistandard tableaux considered are
$$
  P={\smallsquares\Young(1,1,1,2,4,7|2,3,3,5|3,4,6,6|6)}
  \text{and}
  Q={\smallsquares\Young(1,2,2,3,3,6|3,3,3,4|4,5,5,5|5)}
  \text{whose standardisations are}
  {\let~\!\smallsquares\Young(0,1,2,4,9,1~4|3,6,7,1~0|5,8,1~2,1~3|1~1)}
  \text{and}
  {\let~\!\smallsquares\Young(0,1,2,6,7,1~4|3,4,5,9|8,1~1,1~2,1~3|1~0)}
$$
It would take too much space to draw the entire Schensted-growth for this
example, so figure~2 just displays the Knuth-growth derived from it, and
details the Schensted-growth for one of its grid squares. The square that has
been refined is indicated by brackets; matrix entries are indicated when
non-zero.

\midinsert
$$\vcenter{
\catcode`@=11
\def~{\unhcopy\QEDbox}
\setbox\z@=\hbox{~}\dimen@=.5\ht\z@ \advance\dimen@.5\dp\z@
\let\f=\hfil
\def\d#1 {\vtop\bgroup\smallskip
          \vbox to\z@{\vss\llap{$#1$ }\kern.77pt}%
          \futurelet\n\E}
\def\E{\?\ifx\n\f \then\smallskip\egroup \else\afterassignment\e\count@ \fi}
\def\e{\hbox{\loop~\advance\count@\m@ne \ifnum\count@>\z@ \repeat}%
       \nointerlineskip\futurelet\n\E}
\def\g#1{\omit\hfill\clap{\quad#1}\hfill\ignorespaces}
\def\h#1 #2{\vtop to 14pt
            {\medskip\kern-\dimen@\kern#1\dimen@\hbox{#2}\vss} \vrule}
\def\z{\omit\vtop{\smallskip\hbox{$\circ$}}\hfil}
\offinterlineskip
\everycr{\noalign{\nobreak}}
\tabskip=\z@
\halign{\h# \tabskip=0.9em plus 1pc minus 0.9em &&\d#\f\cr
\omit\hfil\strut\vrule
    & \g1 & \g2 & \g3 & \g4\quad & \g5 &  \g6  &  \g7  &  \g8 \cr
\noalign{\hrule}
 1 1 &\z& \z  & \z    & \z     	&  \z     &  \z     &    \z     &    \z    \cr
 1 2 &\z& \z  & \z    & 1  1   	&{}  1    &{} 1     &{} 1       &{}   1    \cr
 1 3 &\z& \z  & \z    &{}  1   	&{}  1    &{} 1 \f
\vtop to\z@{\kern0.23\wd\QEDbox
            \llap{$\left[\vrule depth2\wd\QEDbox width\z@
                   \right.$\kern2.8\wd\QEDbox}\vss}%
                                                    & 2 3 \f
\vtop to\z@{\kern0.23\wd\QEDbox
            \rlap{\kern0.8\wd\QEDbox
                  $\left.\vrule depth2\wd\QEDbox width\z@\right]$}\vss}%
                                                                &{}   3    \cr
 2 4 &\z& 1 1 & 1 2   & 1 3 1  	& 1 4 1   &{} 4 1   & 1 5 3     &{}  5 3   \cr
 3 5 &\z&{} 1 &{} 2   & 1 4 1  	&{} 4 2   & 1 5 2   &{} 5 4 1   &{} 5 4 1  \cr
 3 6 &\z& 2 3 & 1 4 1 &{} 4 3 1 & 1 5 3 2 &{} 5 4 2 &{} 5 4 4 1 &{} 5 4 4 1\cr
 3 7 &\z&{} 3 &{} 4 1 &{} 4 3 1 &{} 5 3 2 &{} 5 4 2 &{} 5 4 4 1 & 1 6 4 4 1\cr
}
}
\hskip 1cm minus 9mm
\left[\,
\vcenter{
\catcode`@=11
\def~{\unhcopy\QEDbox}
\dimen@=.5\ht\QEDbox \advance\dimen@.5\dp\QEDbox
\let\f=\hfil
\def\d#1 {\vtop\bgroup\smallskip
          \vbox to\z@{\vss\llap{$#1$\kern2pt}\kern.77pt}%
          \futurelet\n\E}
\def\E{\?\ifx\n\f \then\smallskip\egroup \else\afterassignment\e\count@ \fi}
\def\e{\hbox{\loop~\advance\count@\m@ne \ifnum\count@>\z@ \repeat}%
       \nointerlineskip\futurelet\n\E}
\offinterlineskip
\everycr{\noalign{\nobreak}}
\tabskip \z@
\halign{\d#\f\tabskip=1em &\d#\f&\d#\f&\d#\f\tabskip=\z@\cr
{} 1   & {} 2   & {} 3   & {} 3   \cr
{} 1 1 & {} 2 1 & {} 3 1 & {} 3 1 \cr
{} 2 1 & {} 2 2 & {} 3 2 & {} 3 2 \cr
{} 3 1 & {} 3 2 & {} 3 3 & {} 3 3 \cr
{} 4 1 & {} 4 2 & {} 4 3 & {} 4 3 \cr
{} 4 1 & {} 4 2 & {} 4 3 & 1  5 3 \cr
}
}
\,\right]
$$
\centerline{{\bf Figure 2.}
 A Knuth-growth for the RSK~correspondence,
 with detail of the associated Schensted-growth.}
\endinsert

For each rectangular area of the Schensted-growth for which the paths along
the bottom and right edges are standardisations of horizontal strips (like the
bracketed area above), any entries~$1$ in the matrix will be produced along
the end of the diagonal that ends in the bottom right corner of the rectangle.
Therefore the computation of the shape datum in the ``insertion'' direction,
determining the shape~$\kappa$ at the bottom right corner when the shapes
$\\$, $\mu$, $\nu$ at the other corners and the matrix entry~$a$ are known,
can be done as follows: first complete a Schensted-growth with the
standardisations of $\mu/\\$ and $\nu/\\$ along its top and left edges and
with all matrix entries equal to~$0$, then extend the grid by $a$ more lines
at the right and at the bottom, copying the skew standard tableaux found along
the right and bottom edges of the rectangle onto each of them, and finally
complete the empty $a\times a$ grid square so created (at the bottom right of
the grid) as a Schensted-growth containing an $a\times a$ identity matrix.
Note that if one combines these individual Schensted-growths to a global one
matching the full Knuth correspondence, these extensions of the grid apply to
an entire row or column; therefore the entries in a single row of~$A$ will
correspond to a collection of entries~$1$ in the final permutation matrix that
descend, when traversed from left to right, into consecutive rows; similarly
entries in a single column of~$A$ give rise to entries~$1$ that, when
traversed from top to bottom, move rightwards into consecutive columns.

Let us describe this the shape datum for the RSK~correspondence, in other
words the relation between the shapes $\\~\mu\choose\nu~\kappa$ at the corners
of a square of the Knuth-growth and the matrix entry~$a$ inside, in somewhat
different terms. Being defined by a part of a Schensted-growth, it can be
computed by a localised version of the Schensted algorithm, in fact by an
instance of the algorithm for skew tableaux of~\ref{Sagan Stanley}. In their
terminology, it can be described (in the insertion direction) as follows:
successively internally insert the squares of the standardisation of~$\nu/\\$
into a tableau of shape~$\mu/\\$ all of whose entries are equal, and then
externally insert $a$ more copies of the same entry; the shape of the
resulting tableau gives $\kappa/\nu$. The fact that $\mu/\\$ and~$\nu/\\$ are
horizontal strips severely restricts the insertion process, which allows us to
describe the shape datum more directly. Only the internal insertion steps can
involve ``bumping'' (corresponding in the Schensted-growth to the case
of~\growthdef\statitemnr{4}) and they can do so at most once. This happens
whenever the square inserted from~$\nu/\\$ also occurs in the diagram of~$\mu/\\$; it results in adding a new square on the first available
place in the next row. We note that for $\mu=\nu$ and~$a=0$ one finds
essentially the bijection given in the introduction to prove the case $r=1$ of
our claim~\alpineclaim.

A somewhat more formal direct description of this shape datum can be given
using a description of horizontal strips as multisets of rows: a horizontal
strip of the form $\\/\mu$ will be described by specifying for each row~$i$
the number $\\_i-\mu_i$ of its squares in that row (while this information
does not completely describe the strip, it suffices in our setting, where
always one of the shapes $\\,\mu$ involved is known beforehand). The operation
of multiset intersection (in which each element gets as multiplicity the
minimum is its multiplicities in the operands) provides a succinct way to
describe the common squares of two strips like $\mu/\\$ and~$\nu/\\$. We also
need an operation that shifts one row upwards: if $A$ is a multiset of rows
then $A^\uparrow$ denotes another such multiset
$\multisetof{i-1}:i\in{A},i>0\endset$, in other words the multiplicity of
row~$i$ in $A^\uparrow$ is the multiplicity of row~$i+1$ in $A$. Now our shape
datum $\bKn$ is then specified by $\bKn_{\mu,\nu}(\kappa)=(a,\\)$ where
$$
  \mu/\\=\kappa/\nu -(\kappa/\mu\thru\kappa/\nu)
                    +(\kappa/\mu\thru\kappa/\nu)^\uparrow
\text{and}
  \nu/\\=\kappa/\mu -(\kappa/\mu\thru\kappa/\nu)
                    +(\kappa/\mu\thru\kappa/\nu)^\uparrow,
\nn
$$
(the two equations are equivalent), and $a=|\kappa|-|\nu|-|\mu|+|\\|$, which
is the multiplicity in $\kappa/\mu\thru\kappa/\nu$ of row~$0$.
Equations~(\lastlabel) express the relation
$\mu/\\\thru\nu/\\=(\kappa/\mu\thru\kappa/\nu)^\uparrow$, which is the fact we
saw above that common squares in row~$i$ of $\mu/\\$ and~$\nu/\\$ are in
bijection with common squares in row~$i+1$ of $\kappa/\mu$ and~$\kappa/\nu$,
as well as the fact that the remaining squares of $\mu/\\$ and of~$\nu/\\$
respectively match squares in the same row in $\kappa/\nu$ and in~$\kappa/\mu$
(in fact the very same squares; this part is already implied by the inclusions
among the shapes $\\,\mu,\nu,\kappa$). These equations can be restated as
identities in terms of the individual parts of the shapes, which can then be
simplified to give the equivalent equations
$$ \eqalign
 {a&=\kappa_0-\max(\mu_0,\nu_0),
  \cr
  \\_i&=\min(\mu_i,\nu_i)+\max(\mu_{i+1},\nu_{i+1})-\kappa_{i+1}
  \ttex{for~$i\in\N$.}
 }
\nn
$$
Using the fact that $\\\leqh\mu$ means $\mu_{i+1}\leq\\_i\leq\mu_i$ for
all~$i\in\N$, one easily checks that with these equations, \
$\mu\leqh\kappa\geqh\nu$ is equivalent to $a\geq0$ and $\mu\geqh\\\leqh\nu$,
so that $\bKn$ is indeed a shape datum.

\subsection An alternative shape datum.
\labelsec\Burgesec

There is different shape datum for $(\Y,\leqh,\N)$ that can be obtained by
similar methods; it suffices to replace the (row insertion) $1$-correspondence
used above by its transpose (column insertion) $1$-correspondence. One may
proceed in exactly the same way to construct a Knuth-growth, interpolating
horizontal strips, then building a Schensted-growth, and finally reducing the
grid again. This works because in a Schensted-growth for the transpose
$1$-correspondence, it still holds that any rectangular area for which the
path along the right or bottom edge of is the standardisation of a horizontal
strip, has such a path along the opposite edge as well (ignoring trivial steps
that repeat the same shape). The Knuth-correspondence so obtained is known as
the Burge correspondence (although only the case where $P$-~and $Q$-symbol
coincide is actually used in~\ref{Burge}). It should not be confused with the
Knuth's asymmetric correspondence mentioned in the introduction; that
correspondence is not of the type we are currently considering, where the
$P$-~and $Q$-symbol must both be semistandard and of equal shape. We
illustrate the construction of the Burge correspondence in figure~3, in a
similar way as we did for the RSK~correspondence.

\midinsert
$$\vcenter{
\catcode`@=11
\def~{\unhcopy\QEDbox}
\setbox\z@=\hbox{~}\dimen@=.5\ht\z@ \advance\dimen@.5\dp\z@
\let\f=\hfil
\def\d#1 {\vtop\bgroup\smallskip
          \vbox to\z@{\vss\llap{$#1\,$}\kern.77pt}%
          \futurelet\n\E}
\def\E{\?\ifx\n\f \then\smallskip\egroup \else\afterassignment\e\count@ \fi}
\def\e{\hbox{\loop~\advance\count@\m@ne \ifnum\count@>\z@ \repeat}%
       \nointerlineskip\futurelet\n\E}
\def\g#1{\omit\hfill\clap{\quad#1}\hfill\ignorespaces}
\def\h#1 #2{\vtop to 14pt
            {\medskip\kern-\dimen@\kern#1\dimen@\hbox{#2}\vss} \vrule}
\def\z{\omit\vtop{\smallskip\hbox{$\circ$}}\hfil}
\offinterlineskip
\everycr{\noalign{\nobreak}}
\tabskip=\z@
\halign{\h# \tabskip=0.9em plus 1pc minus 0.9em &&\d#\f\cr
\omit\hfil\strut\vrule
    & \g1 & \g2 & \g3 & \g4\quad & \g5 &  \g6  &  \g7  &  \g8 \cr
\noalign{\hrule}
 1 1 &\z& \z  & \z    & \z     	&  \z     &  \z     &    \z     &    \z    \cr
 1 2 &\z& \z  & \z    & \z     	&  \z     &  \z     &    \z     & 1   1    \cr
 1 3 &\z& 2 2 & 1 3   &{}  3   	& 1  4    &{} 4     &{} 4       &{}   5    \cr
 2 4 &\z&{} 2 &{} 3   & 1 3 1  	&{} 4 1   & 1 4 2   &{} 4 2     &{}  5 2   \cr
 2 5 &\z& 1 3 & 1 4 1 & 1 4 3  	& 1 5 3 1 &{} 5 4 1 \f
\vtop to\z@{\kern0\wd\QEDbox
            \llap{$\left[\vrule depth4.1\wd\QEDbox width\z@
                   \right.$\kern4.5\wd\QEDbox}\vss}%
                                                    & 1 5 4 2 \f
\vtop to\z@{\kern0\wd\QEDbox
            \rlap{\kern-0.5\wd\QEDbox
                  $\left.\vrule depth4.1\wd\QEDbox width\z@\right]$}\vss}%
                                                                 &{} 6 4 2  \cr
 3 6 &\z&{} 3 &{} 4 1 &{} 4 3   &{} 5 3 1 &{} 5 4 1 & 2 5 4 3 1 &{} 6 4 3 1\cr
 3 7 &\z&{} 3 &{} 4 1 & 1 4 3 1 &{} 5 3 2 &{} 5 4 2 &{} 5 4 4 1 &{} 6 4 4 1\cr
}
}
\hskip 1cm minus 9mm
\left[\,
\vcenter{
\catcode`@=11
\def~{\unhcopy\QEDbox}
\setbox\z@=\hbox{~}\dimen@=.5\ht\z@ \advance\dimen@.5\dp\z@
\let\f=\hfil
\def\d#1 {\vtop\bgroup\smallskip
          \vbox to\z@{\vss\llap{$#1$\kern2pt}\kern.77pt}%
          \futurelet\n\E}
\def\E{\?\ifx\n\f \then\smallskip\egroup \else\afterassignment\e\count@ \fi}
\def\e{\hbox{\loop~\advance\count@\m@ne \ifnum\count@>\z@ \repeat}%
       \nointerlineskip\futurelet\n\E}
\offinterlineskip
\everycr{\noalign{\nobreak}}
\tabskip \z@
\halign{\d#\f\tabskip=0.8em &\d#\f&\d#\f&\d#\f\tabskip=\z@\cr
{} 5 4 1 & {} 5 4 1   & {} 5 4 1   & {} 5 4 2   \cr
{} 5 4 1 & {} 5 4 1   &  1 5 4 1 1 & {} 5 4 2 1 \cr
{} 5 4 1 &  1 5 4 1 1 & {} 5 4 2 1 & {} 5 4 3 1 \cr
}
}
\,\right]
$$
\centerline{{\bf Figure 3.}
  Knuth-growth for the Burge correspondence;
  detail of the associated Schensted-growth.}
\endinsert

To extract a shape datum $\bBu$ from the Burge correspondence, one may
again limit the Schensted-growth to a rectangular area like the bracketed one,
with standardisations of horizontal strips at the right and bottom. From those
standardisations, the rules for a Schensted-growth with the column insertion
$1$-correspondence determine the other shapes and the matrix entries; the
shape~$\\$ at the top left corner and the sum~$a$ of all matrix entries
produced define the shape datum $\bBu_{\mu,\nu}(\kappa)=(a,\\)$. Any
matrix entry~$1$ occurring in the Schensted-growth now adds a square in
column~$0$, which is the first in any horizontal strip in which it occurs, so
it can only be added to the shape~$\\$. Therefore the matrix entries~$1$ form
an anti-diagonal at the top left corner of the rectangle, and the computation
in the opposite (insertion) direction should proceed as follows: construct a
partial Schensted-growth with matrix entries~$1$ along the anti-diagonal of an
$a\times a$ square in the top left corner and~$0$ elsewhere, and shapes~$\\$
repeated along the top and left of edges that square, followed respectively by
standardisations of $\mu/\\$ and $\nu/\\$.

Considering the permutation matrix entries~$1$ contributing to one row or one
column of the final matrix, the above description implies a different
arrangement than for the RSK~correspondence: they are arranged in an
anti-diagonal sense (bottom-left to top-right) in a sequence of consecutive
rows respectively columns. In terms of the original ``bumping'' description,
this means that there are, in addition to using column insertion rather than
row insertion, two more differences of the Burge insertion algorithm with
respect to that of Knuth: the order of insertion is reversed among columns of
the two-line array with equal top index (making it weakly \emph{decreasing} by
bottom index), and among equal entries in the insertion tableau~$P$, the more
recently inserted ones are treated as smaller, so that bumping will replace an
entry of the same value if present in the column.

Although it is of no importance to our paper, we should mention that there
is a relation between the global RSK and Burge correspondences, which 
involves the Sch\"utzenberger involution. It is illustrated by the examples
shown in figures 2~and~3: their matrices are vertical mirror images (the same
rows appear in reverse order), their $P$-symbols are the same while their
$Q$-symbols are Sch\"utzenberger duals:
$$
  P={\Young(1,1,1,2,4,7|2,3,3,5|3,4,6,6|6)}=P'
,\quad
  Q={\Young(1,2,2,3,3,6|3,3,3,4|4,5,5,5|5)}
,\quad
  Q^*={\Young(1,2,2,2,2,4|3,3,4,4|4,4,5,6|5)}=Q'.
$$
This property is readily derived from the well known property of the
Robinson-Schensted correspondence that column inserting permutation entries in
reverse order gives the same $P$-symbol as ordinary row insertion, and the
Sch\"utzenberger dual $Q$-symbol. Note that reversing the order means
vertically reflecting the permutation matrix; to get the same relation between
the final matrices for the RSK and Burge correspondences, the rules for
multiplexing rows and columns must also be reversed.

In spite of this relation between the global correspondences, and the formal
similarity between their definitions, the shape data for the RSK and Burge
correspondences have quite different characteristics. Contrary to what we saw
for the former shape datum, the partial Schensted-growth that defines the
latter can have multiple ``bumping'' configurations \growthdef\statitemnr{4}
in the same row or column of the grid. In fact, when a horizontal strip
contains squares in some range of consecutive columns, any insertion or
extraction step that bumps one of these squares will go on bumping the other
squares until an unoccupied column is reached (or possibly the left edge of
the diagram in case of extraction). By studying the effect of successive
insertions one can deduce a description of the new shape datum in terms of
occupancy of columns; we shall omit the details of the reasoning since this
description will anyway serve merely as a motivation for or the construction
that is our final goal.

If horizontal strips $\mu/\\$ and $\nu/\\$ are given as well as a matrix
entry~$a$, the corresponding shape~$\kappa$ for the Burge shape datum can be
found as follows. We shall treat~$a$ as a modifiable variable, and traverse
the columns from left to right; after visiting column~$j$ we shall have
determined the length $\kappa\tr_j$ of that column in the diagram~$\kappa$.
When visiting column~$j$ we first set
$c=(\mu\tr_j-\\\tr_j)+(\nu\tr_j-\\\tr_j)\in\{0,1,2\}$, which counts the
occurrences in the two given horizontal strips of a square in column~$j$. If
$c=1$ we have $\kappa\tr_j=\\\tr_j+1=\max\{\mu\tr_j,\nu\tr_j\}$, and we
continue to the next column without further action. If $c=2$, then
$\kappa\tr_j=\\\tr_j+1=\mu\tr_j=\nu\tr_j$, and we increase the value of~$a$
by~$1$ before continuing to the next column. In the final case $c=0$ we have
$\\\tr_j=\mu\tr_j=\nu\tr_j$; in this case if $a>0$ then
$\kappa\tr_j=\\\tr_j+1$ and we decrease the value of~$a$ by~$1$, but if $a=0$
we put $\kappa\tr_j=\\\tr_j$ and we continue to the next column without
changing~$a$. This procedure can be continued indefinitely, but there is no
need to go on once a zero-length column of~$\kappa$ is found. Note that
$\kappa\tr_j-\\\tr_j\in\{0,1\}$ for all~$j$, which means that $\\\leqh\kappa$.

As an example, in the bracketed square of the Knuth-growth of figure~3 one
finds for~$c$ successively the values $0,1,0,0,\ldots$ as $j=0,1,2,3,\ldots$
while $a=2$ initially; it follows from the above description that
$\kappa\tr_j=\\\tr_j+1$ for $j=0,1,2$ while $\kappa\tr_j=\\\tr_j$ for all
other values of~$j$, which can be checked in the diagram. For a case where $a$
also increases, consider the square with
$\\={\smallsquares\smallsquares\Young(,)}$ in the row labelled~$4$ and the
column labelled~$2$, which also has
$\mu=\nu={\smallsquares\smallsquares\Young(,,)}$, and~$a=1$. Then for
$j=0,1,2,3$ one finds $c=0,0,2,0$, so that $a$~decreases to~$0$ for $j=0$,
stays so for $j=1$, raises again to~$1$ for $j=2$ and finally drops back
to~$0$ for~$j=3$; this causes $\kappa\tr_j=\\\tr_j+1$ for $j=0,2,3$ and
$\kappa\tr_j=\\\tr_j$ for other values of~$j$, whence
$\kappa={\smallsquares\smallsquares\Young(,,,|)}$.

It may be checked independently from the relation with Schensted insertion
that this description defines a shape datum. The variable~$a$ ultimately
becomes $a=0$, since $c=0$ ultimately, after which $a$ must decrease until it
is~$0$. Then one easily sees that (\rankeq) holds (for the initial value
of~$a$ of course). One never has $\kappa\tr_j>\kappa\tr_{j-1}$ for $j>0$ (a
column of~$\kappa$ made longer than its predecessor), since this requires
$\kappa\tr_j>\\\tr_j=\mu\tr_j=\nu\tr_j$, and therefore $a>0$ when column~$j$
is visited, but if visiting column~$j-1$ leaves $a>0$ one always has
$\kappa\tr_{j-1}=\\\tr_{j-1}+1$. Finally, an inverse operation (which is the
actual direction of the shape datum $\bBu$) is easily defined in the
same terms. The variable~$a$ will assume the same sequence of values as
before, but of course in the opposite order; in particular it starts at~$0$
and its final value gives the number~$a$ to be determined. For the number~$c$,
which is now determined by the horizontal strips $\kappa/\mu$
and~$\kappa/\nu$, the values are not the same as for the forward direction,
but they can be read off in all specified cases of the forward operation. They
can be checked to indeed always produce the inverse change of the variable~$a$
and reconstruct the column length~$\\\tr_j$, thus guaranteeing a step-by-step
inverse procedure.

Our specific interest in the shape datum for the Burge correspondence is due
to two properties: the first is that it admits a description that treats the
squares collectively rather than one at a time, the second is that this
description can be easily stated in terms of edge sequences. To demonstrate
the second point, whose importance will become clear when ribbons are
considered, we shall use (in addition to the variable~$a$) a modifiable doubly
infinite bit sequence~$w$. For the insertion direction that we shall describe,
it starts as $w=\edge\\$ and ends as $w=\edge\kappa$. We shall traverse $w$
from left to right, considering a pair of adjacent bits at the time; such a
pair $w_{k-1},w_k$ corresponds to edges crossing the diagonal
$d_k=\setof(i,j)\in\N^2:j-i=k\endset$, and as done above for columns, we count
the occurences in $\mu/\\$ and $\nu/\\$ of a square in the diagonal~$d_k$,
calling the resulting number $c_k\in\{0,1,2\}$.

For each index~$k$ we proceed as follows. If $(w_{k-1},w_k)\neq(1,0)$, we move
to $k+1$ without further action. If $(w_{k-1},w_k)=(1,0)$ and $c_k=1$, then we
set $(w_{k-1},w_k):=(0,1)$ before proceeding to~$k+1$. If
$(w_{k-1},w_k)=(1,0)$ and $c_k=2$, then we also set $(w_{k-1},w_k):=(0,1)$,
but in addition $a:=a+1$. If $(w_{k-1},w_k)=(1,0)$ and $c_k=0$ (the case where
a square of $d_k$ can be added in~$\kappa$ that is not already in $\mu$
or~$\nu$), the action depends on the current value of~$a$: if $a=0$ nothing
happens, but if $a>0$ we set $(w_{k-1},w_k):=(0,1)$ and $a:=a-1$. For this
procedure to work properly it must traverse a sufficiently large range of
diagonals; we may start with the smallest $k$ for which $(w_{k-1},w_k)=(1,0)$
initially (i.e., with $w=\edge\\$; in fact this gives $k=-\\\tr_0$), and we
may stop once $a=0$ is obtained after treating the last diagonal with $c_k>0$.
Here is the calculation for the two examples cited above
\bigdisplay
\def~{\copy\QEDbox}
\tabskip=0pt plus 1 fil
\halign to \hsize {\hfil#\tabskip=1em minus 0.3em &
  \hfil$\cdots#\cdots$\hfil \tabskip=0.5em minus 0.3em &
  \hfil$#$\hfil & \hfil$#$\hfil & \hfil$#$\hfil \ \tabskip=1em minus 0.3em  &
  #\hfil& \hfil$\cdots#\cdots$\hfil \tabskip=0.5em minus 0.3em &
  \hfil$#$\hfil & \hfil$#$\hfil \tabskip=0pt plus 1 fil\cr
\omit $k$\hfil& \omit \hfil\it $w$ before\hfil &
 \mu\!/\!\\ & \nu\!/\!\\ & c_k &
 \it action &
 \omit \hfil\it $w$ after\hfil &
  \kappa\!/\!\mu & \kappa\!/\!\nu \cr
\omit & \omit\hfil$\edge\\={}$\hfil & & & & $a=2$ \cr
-3 & 11[10]100010100 &   &   & 0 & $a:=a-1=1$, swap
   & 11[01]100010100 & ~ & ~ \cr
-2 & 110[11]00010100 &   &   & 0 & none
   & 110[11]00010100 \cr
-1 & 1101[10]0010100 &   & ~ & 1 & swap
   & 1101[01]0010100 & ~ &   \cr
 0 & 11010[10]010100 &   &   & 0 & $a:=a-1=0$, swap
   & 11010[01]010100 & ~ & ~ \cr
\omit & \omit         &   &  &   & terminate 
   & \omit\hfil${}=\edge\kappa$\hfil\cr
\noalign{\medskip}
\omit & \omit\hfil$\edge\\={}$\hfil & & & & $a=1$ \cr
-1 & 11[10]010000 &   &   & 0 & $a:=a-1=0$, swap
   & 11[01]010000 & ~ & ~ \cr
 0 & 110[10]10000 &   &   & 0 & none ($a=0$)
   & 110[10]10000 \cr
 1 & 1101[01]0000 &   &   & 0 & none
   & 1101[01]0000 \cr
 2 & 11010[10]000 & ~ & ~ & 2 & $a:=a+1=1$, swap
   & 11010[01]000 &   &   \cr
 3 & 110100[10]00 &   &   & 0 & $a:=a-1=0$, swap
   & 110100[01]00 & ~ & ~ \cr
 \omit & \omit   &    &   &  & terminate 
   & \omit\hfil${}=\edge\kappa$\hfil\cr
}
$$
To show that the inverse of this procedure defines a shape datum requires a
bit of effort, but is not too difficult. To see that
$\mu\union\nu\subset\kappa$ requires proving (by induction on~$k$) that
$c_k>0$ implies that $(w_{k-1},w_k)=(1,0)$ holds when this pair of bits is
considered. Termination with $a=0$, and equation~(\rankeq) for the initial
value of~$a$, are then proved like for the column description. Then one proves
(again by induction on~$k$) that cases with $w_{k-1}=0$ can be encountered
only when $a=0$ holds. This fact allows an inverse procedure to be formulated;
its description is in fact very similar to the one given above.

For the RSK~correspondence a similar description of the shape datum is
possible, but would be considerably messier than this one; this is essentially
due to the fact that squares are bumped in the direction opposite to the
traversal of the horizontal strips (namely towards the bottom left).

\subsection An \r-correspondence with $r>1$, and a corresponding shape datum. 
\labelsec\Yrcor 

So far we have only shown examples of $1$-correspondences, and related shape
data. Our goal however are Knuth correspondences whose shape data restrict to
\r-correspondences with $r>1$. Before discussing the ones that really interest
us, let us treat a construction that builds such correspondences in a fairly
trivial way. For these correspondences the natural replacement for the graded
set $S=\N$ with generating series $F_\N(T)=\sum_{i\in\N}X^i={1\over1-T}$ will
be $S=\N^r$, graded by $|a|=\sum_{i\in\set{r}}a_i$, with generating series
$$
  F_{\N^r}(T)=F_{\N}(T)^r
 ={1\over(1-T)^r}=\sum_{i\in\N}\multichoose(r,i)T^i
 \ttex{where $\multichoose(r,i)=\Card\setof a\in\N^r:|a|
             =i\endset={r+i-1\choose i}$.}
\label(\Nrgeneq)
$$

The simplest example of an \r-differential poset is~$\Y^r$, the set of
\r-tuples of Young diagrams, partially ordered by inclusion of each of the
$r$~components separately. Like~$\Y$ this is a distributive lattice, and to
check equation~(\diageq), one observes that any element covering a given
\r-tuple of shapes is obtained by adding a square to one of the component
shapes, while keeping the rest fixed; since for each component there is one
more shape covering it than there are covered by it, the number of elements
covering a given \r-tuple exceeds the number of elements covered by it by~$r$.
In fact it is easy to see that the Cartesian product of an \r-differential
poset and an $s$-differential poset always gives an $r+s$-differential poset.

It is not difficult either to define an \r-correspondence for this situation.
One starts by choosing a $1$-correspondence for~$(\Y,\prec)$ that is to be
used in the individual components; it will be denoted by $(b_\\)_{\\\in\Y}$
(it would be perfectly legal to make different choices for each component, but
such frivolity would only complicate notation). If
$\vec\\=(\\^i)_{i\in\set{r}}$ is an \r-tuple of shapes covered by another such
\r-tuple $\vec\mu=(\mu^i)_{i\in\set{r}}$, there is a unique index~$i$ for
which $\\^i\neq\mu^i$, and one has~$\\^i\prec\mu^i$. Then if
$b_{\\^i}(\mu^i)=e_0$ one defines $b_{\vec\\}(\vec\mu)=e_i\in e_\set{r}$,
while in other cases $b_{\vec\\}(\vec\mu)$ is defined by replacing $\mu^i$ in
$\vec\mu$ by $b_{\\^i}(\mu^i)\in\Y$.

If $\vec\\~\vec\mu\choose\vec\nu~\vec\kappa$ is a grid square of a
Schensted-growth for such a \r-correspondence, with matrix entry~$0$, then
either $\vec\\$ and $\vec\mu$ differ in the same component as $\vec\nu$ and
$\vec\kappa$, or there is equality in both cases; a similar relation holds for
$\vec\nu/\vec\\$ and $\vec\kappa/\vec\mu$. A grid square with matrix
entry~$e_i$ introduces a difference in component~$i$ across both its row and
column of the grid, and by the above property these differences will propagate
along row and column into steps of the $Q$-~and $P$-symbols that will still
involve a change in component~$i$. Consequently, the Schensted correspondence
may be computed in each of the components separately, as follows: for each
$i\in\set{r}$ the positions of the entries~$e_i$ determine on one hand the
sets of steps in the $P$-~and $Q$-symbols that will involve a change in
component~$i$, and on the other hand a permutation matrix (by removing all
other rows and columns, and replacing each $e_i$ by~$1$). For each~$i$ the
Schensted correspondence for~$\Y$ then defines a pair of paths~$(P,Q)$
in~$\Y$, which $r$ pairs can be spliced together into pairs of paths in~$\Y^r$
by taking at each step the next change in the appropriate component~$i$.

The same idea also works to create new Knuth correspondences. On $\Y^r$ one
defines $\vec\\\leqh\vec\mu$ to mean $\\^i\leqh\mu^i$ for all~$i\in\set{r}$,
and for $S$ one takes the graded set $\N^r$. Then a shape datum can be defined
by components: $b_{\vec\mu,\vec\nu}(\vec\kappa)=(\vec{a},\vec\\)$ where
$(a_i,\\^i)=b_{\mu^i,\nu^i}(\kappa^i)$ for all~$i\in\set{r}$; it is trivial to
verify the conditions of definition~\sddef. The Knuth correspondence for this
shape datum operates independently in each component even more evidently than
the Schensted correspondence above: matrices with entries in $S=\N^r$ can be
viewed as \r-tuples of matrices with entries in~$\N$, and paths in
$(\Y^r,\leqh)$ as \r-tuples of paths in $(\Y,\leqh)$ (i.e., of semistandard
tableaux); the Knuth-growth is defined with no interaction whatsoever between
different components, so one may compute Knuth-growths separately for each
component.

\newsection Ribbons, and the Shimozono-White \r-correspondence.
\labelsec\ribbonSec

\subsection Ribbons, edge sequences, \r-cores and \r-quotients.

Partitions, partially ordered by repeated removal of \r-ribbons (also called
\r-rim hooks), provide a less artificial example of a poset that satisfies the
commutation relation~(\DUcomm) for some~$r>1$. The relation $\\\precr\mu$
between Young diagrams is said to hold if $\\\ssubset\mu$, and if the diagram
of $\mu/\\$ consists of one square on each of $r$ consecutive diagonals (in
other words: that skew diagram is connected, has $r$ squares, and contains no
$2\times2$ blocks). In this case $\\$ is said to be obtained by removing an
\r-ribbon from~$\mu$, and $\mu$ by adding an \r-ribbon to~$\\$. The reflexive
transitive closure of~`$\precr$' defines a partial ordering~`$\leqr$' on~$\Y$.

In order to make $(\Y,\precr)$ into a graded graph as discussed in
\Sec\Schenstedsec, the grading on~$\Y$ must be adapted so that $\\$ and~$\mu$
have consecutive ranks when $\\\precr\mu$. This is easily done by defining a
grading $\rwt\\{}=\lfloor{|\\|\over r}\rfloor$, the quotient of the Euclidean
division of $|\\|$ by~$r$. The rank of individual shapes is of no importance,
since only differences in rank between comparable elements in $(\Y,\leqr)$ are
used; therefore we put $\rwt\mu\\=\rwt\mu{}-\rwt\\{}={|\\|-|\mu|\over r}$
whenever $\\\leqr\mu$. Paths of shape $\mu/\\$ in~$(\Y,\precr)$ are called
standard \r-ribbon tableaux of shape $\mu/\\$; such tableaux have $\rwt\mu\\$
individual \r-ribbons.

The operation of adding or removing an \r-ribbon is best understood in terms
of the edge sequences $\edge\\$ and~$\edge\mu$. We shall discuss these matters
here summarily, referring to~\ref{edge sequences} for an extensive discussion
and examples. When $\\\precr\mu$, the edge sequences $\edge\\$ and~$\edge\mu$
differ only in two places, which are at distance~$r$, and at those places
$\edge\\$ has $\cdots1\cdots0\cdots$ while $\edge\mu$ has
$\cdots0\cdots1\cdots$ (the dots represent unchanged bits). Hence addition or
removal of a single \r-ribbon affects bits of the edge sequence whose
positions are in the same congruence class modulo~$r$; we shall say two
\r-ribbons are in the same position class if the positions of the bits affected
by their addition or removal are congruent modulo~$r$. As operations on edge
sequences, the addition or removal of \r-ribbons in distinct position classes
always commute (although the the ribbons themselves, viewed as skew diagrams,
may change). The bits with positions in a fixed congruence class form, up to a
shift, the edge sequence of a unique Young diagram, which changes by a single
square for any modification by an \r-ribbon in the corresponding position
class. Thus the induced sub-poset of $(\Y,\leqr)$ on the subset of shapes
reachable from a given one by addition or removal of \r-ribbons in a single
position class is isomorphic to~$\Y=(\Y,\subset)$. It follows that any
connected component of $(\Y,\leqr)$ is isomorphic to $\Y^r$. In particular
each such component has a unique minimal element for~`$\leqr$', which is
called an \r-core; $\gamma$ is an \r-core if and only if each of the $r$
sequences extracted from $\edge\gamma$ by selecting a congruence class
modulo~$r$ of bit positions is of the form $\cdots111000\cdots$ (these
extracted sequences may be shifted with respect to each other). Any $\\\in\Y$
is uniquely determined by the minimal element in its connected component of
$(\Y,\leqr)$, called the \r-core of~$\\$, together with its image in $\Y^r$
under the isomorphism, an \r-tuple of partitions called the
\r-quotient of~$\\$.

\subsection Rim hook lattices and semistandard \r-ribbon tableaux.

Being isomorphic to~$\Y^r$, any connected component of $(\Y,\leqr)$ is clearly
a distributive lattice; it is called an \r-rim hook lattice. As all \r-rim
hook lattices are isomorphic as posets, one often considers only the one
containing the empty partition as \r-core, but there is no good reason to do
so: if one were only interested in the abstract poset structure, one could
study~$\Y^r$, and it would be useless to introduce ribbons in the first place.
Meanwhile the isomorphism, which is called the \r-quotient map, makes it
immediately clear that \r-rim hook lattices are \r-differential posets, and it
translates any choice of an \r-correspondence for~$\Y^r$ into such a choice
for any \r-rim hook lattice. In the case of the \r-correspondence for~$\Y^r$
of~\Sec\Yrcor, using in each component~$\Y$ the usual row-insertion
(respectively column-insertion) $1$-correspondence, the resulting
\r-correspondence for~$(\Y,\leqr)$ can be described directly as follows. Let
$\mu\precr\kappa$, then $b_\mu(\kappa)=\\\precr\mu$ is such that $\mu/\\$ is
the first \r-ribbon, if it exists, that is removable from~$\mu$, in the same
position class as $\kappa/\mu$ and to its top right (respectively to its
bottom left). If no such ribbon exists, $b_\mu(\kappa)=e_i$ where
$i\in\set{r}$ represents the position class of~$\kappa/\mu$: if the top right
square of the ribbon, which we shall called its head, lies on diagonal~$d_k$
then $i=k\bmod{r}$. The validity of this description rests on the fact that
the \r-ribbons in a single position class that can be added to respectively
removed from~$\mu$, are perfectly interleaved in bottom left to top right
order, with at both extremes ribbons that can be added; this can be seen
directly by the same argument as given in~\Sec\Schenstedsec\ for the
$1$-correspondence for~$\Y$, but applied to the sequence extracted
from~$\edge\mu$ of bits at positions congruent to~$i$ modulo~$r$. Note that
for these \r-correspondences in~$(\Y,\leqr)$ the names row- and
column-insertion would not be very appropriate, since the ribbons $\kappa/\mu$
and $\mu/\\$ might be separated by any number of rows and columns in either
case.

The \r-quotient maps can of course also be used to transport the
relation~`$\leqh$' and the shape data for it from~$\Y^r$ to~$(\Y,\leqr)$. We
shall define $\\\leqhr\mu$ when $\\\leqr\mu$ and moreover $\\^i\leqh\mu^i$ for
each pair of corresponding components $\\^i,\mu^i$ of their respective
\r-quotients. To better understand this relation, note that 
$\\\leqh\mu$ means that $\edge\\$ can be transformed into~$\edge\mu$ by a
sequence of replacements of a substring~`$10$' by~`$01$', proceeding from left
to right with overlap allowed (the bit~$1$ of the replacement may participate
in the next replacement). This fact, which we already saw implicitly in the
description of the Burge shape datum by edge sequences, follows from the
simple observation that the sequence of squares added in a skew standard
tableau has strictly increasing column numbers if and only if it has
increasing diagonal numbers (here strictness is for free). In this
description, the leftmost bit of a replacement string cannot be modified
afterwards, so while considering from left to right occurrences of~`$10$' for
replacement, the decision whether to or not replace is prescribed by the
target~$\edge\mu$ (if a possibility exists at all); therefore this description
gives a direct, backtrack-free, procedure to decide whether or not
$\\\leqh\mu$.

To decide whether $\\\leqhr\mu$ holds, we must make a similar test for each
pair of sequences of bits extracted from~$\\$ and~$\mu$ at positions in the
same conjugacy class modulo~$r$. Although the $r$ traversals required for
these tests are independent, we might as well combine them into a single left
to right pass over the edge sequences. Thus we find that $\\\leqhr\mu$ if and
only if $\edge\\$ can be transformed into~$\edge\mu$ by a left to right
sequence of replacements of a substring~`$1x0$' by~`$0x1$', where
$x\in\{0,1\}^{r-1}$ is any string of $r-1$ bits unaffected by the replacement,
and with overlap between successive replacements allowed. Each intermediate
bit string occurring during this transformation is edge sequence of a Young
diagram, and their sequence defines a standard \r-ribbon tableau of
shape~$\mu/\\$; the skew shape~$\mu/\\$ is called a \newterm{horizontal
\r-ribbon strip}, and the standard \r-ribbon tableau its standardisation. In
geometric terms, the standardisation of a horizontal \r-ribbon strip~$\mu/\\$
gives the unique decomposition of its diagram into a sequence of \r-ribbons
such that the head of each ribbon has its top edge on the inner border
of~$\mu/\\$ (i.e., on the boundary of~$\\$), or equivalently such that the
tail (bottom left square) of each ribbon has its bottom edge on the outer
border of~$\mu/\\$ (i.e., on the boundary of~$\mu$). A monotonically rising
path in~$(\Y,\leqhr)$ from $\\$~to~$\mu$ is called a \newterm{semistandard
\r-ribbon tableau} of shape~$\mu/\\$. Here are graphic representations of the
standardisation of a horizontal $5$-ribbon strip, and of a semistandard
$6$-ribbon tableau.
$$\smallsquares
{\rtab5
 \edgeseq8,0;1111111100000000000000
 \ribbon8,0:;0001
 \ribbon8,4:;1101
 \ribbon7,5:;0110
 \ribbon7,7:;1010
 \ribbon6,9:;0111
 \edgeseq3,11;010
 \ribbon1,13:;0100
}
\kern3cm
{\rtab6
 \edgeseq6,0;1111110000000
 \ribbon6,0:0;01101 \ribbon4,3:0;11100 \ribbon1,6:0;01000
 \ribbon5,2:1;00111 \ribbon1,8:1;00010
 \ribbon7,0:3;00100 \ribbon4,5:3;11000 \ribbon4,6:3;10001 \ribbon2,10:3;00101
 \ribbon7,3:4;00110 \ribbon2,13:4;01100
 \ribbon4,7:5;00010 \ribbon6,6:7;01000
}
$$
Although we used the \r-quotient map to define~`$\leqhr$', this is not evident
from the final description; usually horizontal \r-ribbon strips and
semistandard \r-ribbon tableaux are defined without using~\r-quotients.

\subsection Ribbon Schensted and Knuth correspondences that factor.
\labelsec\factorsec

Fomin's constructions reduce the question of defining Schensted and Knuth
correspondences for ribbon tableaux to the question of defining
\r-correspondences for~$(\Y,\leqr)$, respectively of defining shape data
for~$(\Y,\leqhr,\N^r)$. The \r-quotient map provides an easy way to do this,
as we already mentioned for \r-correspondences. We note that the Schensted
correspondence obtained from the \r-correspondence based on the row insertion
in each of the $r$~components was originally defined by a direct construction
in~\ref{Stanton White}, and was later found to factor via the \r-quotient map.
To define a shape datum $b^r$ for~$(\Y,\leqhr,\N^r)$ one may proceed
similarly, leading to the following description. Whenever
$(a,\\)=b^r_{\mu,\nu}(\kappa)$ is to be defined, each of $\mu,\nu,\kappa$ has
the same \r-core~$\gamma$; denoting by $\mu^i,\nu^i,\kappa^i$ the components
of their respective \r-quotients, $\\$ will be the partition with
\r-core~$\gamma$, and with the components of its \r-quotient defined, together
with the components of~$a$, by $(a_i,\\^i)=b_{\mu^i,\nu^i}(\kappa^i)$ for
all~$i\in\set{r}$, where one fixes either $b=\bKn$ or $b=\bBu$.

In this construction of shape data, the standardisations of the horizontal
\r-ribbon strips that arise are not used at all. We know these
standardisations exist, since within a \r-rim hook lattice $\\\leqhr\mu$ just
means that $\\^i\leqh\mu^i$ for all~$i$, but to find the standardisation of a
given horizontal \r-ribbon strip amounts to carefully ordering the individual
squares of the $r$ horizontal strips contributing to it, according to the
placements of the associated edges within the edge sequence merged together
from $r$ individual ones. Since the Knuth correspondences used in the
components were derived via standardisation from Schensted correspondences,
one might assume that these shape data, and the derived Knuth correspondences,
can also be obtained by combining the Stanton and White correspondence with
standardisation of horizontal \r-ribbon strips; indeed one can for instance
find an informal statement to this effect in \ref{Shimozono White,~\Sec5}.
Things do not work out that nicely however; in particular, one does not obtain
the proper standardisations directly from the Stanton and White algorithm,
even for the half-semistandard case (with the $Q$-symbol standard). To
illustrate this point, we shall briefly digress to study this Knuth
correspondence for ribbon tableaux (whose existence, based on arguments like
those above, has been mentioned occasionally, but which does not appear to
have been explicitly described anywhere) in some more detail.

The easiest way to define a correspondence based on standardisation is to
start with semistandard \r-ribbon tableaux $P,Q$, and (just like we did for
the RSK~correspondence) interpolate into standard tableaux, compute the
Schensted-growth for the Stanton and White \r-correspondence in the extraction
direction, and afterwards reduce the grid again, summing coloured permutation
matrix entries over rectangular areas that collapse to a square, to find the
matrix entries $a\in\N^r$ of the final result (addition is done for each
colour separately). This produces the same Knuth-growth and hence the same
matrix as the Knuth correspondence for the shape datum~$b^r$ defined above
(with $b=\bKn$), since both standardisation and application of the
\r-correspondence commute with the decomposition of shapes via the \r-quotient
map. More precisely, decomposing the shapes in the standardisation of a
horizontal \r-ribbon strip produces a sequence of \r-tuples of shapes, in
which only one component changes at each step; to extract from this an
\r-tuple of skew standard tableaux, one must eliminate in each component the
steps where the shape does not change, but then the resulting tableaux are
indeed the standardisations of the $r$ horizontal strips obtained from the
decomposition of the original horizontal \r-ribbon strip via the \r-quotient
map.

However, while the correct shapes are assigned to the grid points that remain
after reduction, the interpolations between them produced by the Schensted
growth are not the standardisations of the respective horizontal \r-ribbon
strips, and as a consequence the coloured permutation entries are not placed
in any predictable manner within their rectangle of the Schensted-growth; this
means that it is impossible to reconstruct the same Schensted-growth from the
matrix found by an inverse (insertion) procedure. The reason for this is that
the relative order of ribbons from different position classes contributing to
a horizontal \r-ribbon strip can freely change during the extraction process,
since the rules of the growth treat the $r$~classes of positions of ribbons
completely independently. The only thing that is guaranteed is that ribbons
within one position class keep their relative order, so that under the
\r-quotient map one gets standardisations of horizontal strips. Here are some
very simple examples of what can happen, with $r=2$ and without non-zero
matrix entries. The first one uses the Stanton and White
\r-correspondence, the second one uses the transpose \r-correspondence. Note
that the top rows are not the standardisations of their respective horizontal
\r-ribbon strips, while the bottom rows are such standardisations.
$$
\def\young(#1){\vtop{\kern0pt \smallsquares\Young(#1)}\hfill}
\def~{\vtop{\kern0pt \hbox{$\circ$}}}
\def\:{\lower10pt \hbox{$\longrightarrow$}}
\def\rt#1{\vtop{\kern0pt
                \hbox{$\smallsquares{\rtab2 #1 \edgeseq3,0;111000 }$}}\hfill}
\matrix { ~         & \young(,)  & \young(,|,)
  &\: &  \rt{\ribbon0,0:0;0 \ribbon1,0:1;0} \cr
   \noalign{\medskip}
  	  \young(,) & \young(,||)& \young(,|,|,)
  &\:  &  \rt{\ribbon2,0:0;1 \ribbon2,1:1;1} \cr}
\kern0.3\hsize
\matrix { ~         & \young(,)  & \young(,|,)
  &\: &  \rt{\ribbon0,0:0;0 \ribbon1,0:1;0} \cr \noalign{\medskip}
  	  \young(|) & \young(,|,)& \young(,,|,,)
  &\:  &  \rt{\ribbon1,1:0;1 \ribbon1,2:1;1}\cr}
$$
Note that in both cases one ribbon keeps its head on the diagonal, while the
other is bumped, moving two diagonals and thereby inverting the relative order
of the two ribbons. Such behaviour occurs during insertion as well as
extraction, so regardless of the rule one would choose to insert the different
contributions of a single matrix entry into the Schensted-growth, it is
impossible to achieve that the sequences interpolating horizontal \r-ribbon
strips remain their standardisations throughout the insertion process. All in
all, there is no reason that would justify using a Schensted-growth instead of
computing the Knuth correspondence directly via the \r-quotient map.

There is one more point that is worth noting about the shape datum $b^r$ of
this Knuth correspondence obtained from $b=\bBu$, the shape datum of the Burge
correspondence. In this case the description of~$\bBu$ in terms of edge
sequences given in~\Sec\Burgesec\ can be easily adapted to the ribbon setting,
merging $r$ independent traversals of the edge sequence into one in the same
manner as for the definition of horizontal \r-ribbon strips. Although it is
still a parallelised correspondence in disguise, we sketch it here because of
its similarity in form to the shape datum that we shall present later.

As before we consider insertion, so let shapes $\mu\geqhr\\\leqhr\nu$ and a
matrix entry~$a\in\N^r$ be given; also, $k$ shall be the index of a diagonal,
for which this time we count the number $c_k\in\{0,1,2\}$ of \r-ribbons in the
standardisations of $\mu/\\$ and of~$\nu/\\$ whose head lies on the
diagonal~$d_k$. A double infinite bit string~$w$ is again initialised
to~$\edge\\$, but instead of inspecting pairs $(w_{k-1},w_k)$ of adjacent
bits, one inspects the pairs $(w_{k-r},w_k)$ that would be affected by an
\r-ribbon with its head on~$d_k$. The action taken for~$k$ may swap these two
bits and may change the component~$a_i$ of~$a$ where $i=k\bmod{r}$, according
to the same case distinction as before; after this processing the value of~$k$
is increased by~$1$. The range of values that must be traversed by~$k$ is
determined in a similar manner as before, where we must now detail that all
components of~$a$ must have become~$0$ before termination can be decided.

A verification that (the inverse of) this procedure defines a shape datum
for~$(\Y,\leqhr,\N^r)$ can be given along the same lines as indicated at the
end of~\Sec\Burgesec. Since~$w_k$, the rightmost of the bits considered and
possibly modified, is put aside for some time before being reconsidered, the
condition that enables the formulation of an inverse procedure now read as
follows: after treating the pair $(w_{k-r},w_k)$, one has for all~$i$ with
$k-r<i\leq k$ that $w_i=0$ implies $a_{i\bmod{r}}=0$. Surprisingly this
algorithm, which is rather trivial because it factors into $r$ copies
simultaneously computing the shape datum for the Burge correspondence, can be
transformed into the spin preserving shape correspondence that is the main
subject of this paper by just one change: the selection $i=k\bmod{r}$ at each
step of the component of~$a$ to possibly modify is replaced by
$i=\sum_{j=1}^{r-1}w_{k-j}$. Proving the existence of an inverse will be a
bit harder though.

\subsection Form and height of ribbons.

From the above consideration we learn that from the enumerative point of view
${(\Y,\leqhr,\N^r)}$ allows the existence of Knuth correspondences, but the
ones found so far are not very interesting because they avoid really dealing
with ribbons by directly applying the \r-quotient map. The key to defining
interesting alternative correspondences is to focus on a property of
\r-ribbons that is not related to the \r-quotient maps. One such property is
their form, where two ribbons are considered to have the same form if their
diagrams are equal up to a translation. There are $2^{r-1}$ different forms
of~\r-ribbons. If $\edge\\$ can be transformed into $\edge\mu$ by the
replacement of a substring $1x0$ by $0x1$ with $x\in\{0,1\}^{r-1}$, then $x$
describes the form of the \r-ribbon $\mu/\\$ by telling for each following
square whether it is above or to the right of its predecessor; we therefore
define $\form(\mu/\\)=x$. Now an important statistic on ribbons is their
height, where $\hgt(\mu/\\)\in\set{r}$ is defined as the sum of the bits
in~$\form(\mu/\\)$. This is the number of vertical steps encountered when
going from the tail of~$\mu/\\$ to its head, or (since even ribbons without
vertical steps occupy one row) one less than the number of rows that meet the
diagram of~$\mu/\\$.

The \r-quotient map is not very well suited for studying the height of
ribbons, since the bits that contribute to the height of a ribbon are found
only in the components of the \r-quotient distinct from the one to which the
ribbon contributes a square. In fact, although \r-core and \r-quotient
together preserve complete information about a shape, the \r-quotients of
$\\$~and~$\mu$ alone do not suffice to determine $\hgt(\mu/\\)$, since it
cannot be determined \emph{which} bits of the other components should be added
(the edge sequences of the components of the \r-quotient are shifted by
amounts determined by the \r-core). Consequently, there is very little one can
say about the heights of the ribbons related by the \r-correspondences
for~$(\Y,\precr)$ and shape data for~${(\Y,\leqhr,\N^r)}$ discussed above. We
shall presently study correspondences that do respect the heights of ribbons;
it is these that will justify our interest in \r-rim hook lattices. Let us
start however by studying how the height of ribbons behaves in the cases
encountered in Schensted correspondences that are independent of any
particular choice of \r-correspondence.

Consider shapes $\\~\mu\choose\nu~\kappa$ with $\\\precr\mu,\nu\precr\kappa$
and $\mu\neq\nu$, in other words $\kappa$ is the unique element covering $\mu$
and~$\nu$ in~$(\Y,\precr)$, and $\\$ is the unique element covered by them.
Then the \r-ribbons $\mu/\\$ and~$\kappa/\nu$ will have their head on the same
diagonal (their heads may even coincide), as will $\nu/\\$ and~$\kappa/\mu$.
If these diagonals are respectively $d_s$ and~$d_t$, then $s-t$ cannot lie in
$\{-r,0,r\}$ (one can never successively add two ribbons with their heads on
the same diagonal, and $|s-t|=r$ would imply that there is only one shape
strictly between $\\$ and~$\kappa$ in~$(\Y,\leqr)$, contradicting
$\mu\neq\nu$). If $|s-t|>r$, then $\hgt(\mu/\\)=\hgt(\kappa/\nu)$ and
$\hgt(\nu/\\)=\hgt(\kappa/\mu)$, since the ribbons are separated by at least
one diagonal, and the indicated pairs of ribbons have identical diagrams. On
the other hand, if $|s-t|<r$, then both equalities of height fail, since the
fact of having added a ribbon with its head on diagonal~$s$ to~$\\$ before
forming a ribbon with its head on diagonal~$t$ will have changed one bit in
the form of the latter ribbon. In more detail: if $-r<s-t<0$ ($\kappa/\mu$
lies further to the top right than~$\mu/\\$) then
$\hgt(\kappa/\mu)=\hgt(\nu/\\)+1$ and $\hgt(\kappa/\nu)=\hgt(\mu/\\)-1$:
$$
 {\rtab7
  \ribbon3,0:\hss\mu/\\\hss;110010 \ribbon2,1:~\kappa/\mu\hss;001010
 }
\qquad
 {\rtab7
  \ribbon1,0:~\nu/\\\hss;001000 \ribbon3,0:\hss\kappa/\nu\hss;100010
 };
$$
similarly, if $0<s-t<r$ then $\hgt(\kappa/\mu)=\hgt(\nu/\\)-1$ and
$\hgt(\kappa/\nu)=\hgt(\mu/\\)+1$. Note that in all cases
$\hgt(\kappa/\mu)+\hgt(\kappa/\nu)=\hgt(\nu/\\)+\hgt(\mu/\\)$. In the
Schensted correspondence the ribbons at one side of this equation are computed
from those at the other side, so this amounts to a conservation of total
height.

\subsection The Shimozono-White \r-correspondence.

Now we consider \r-correspondences that respect the height of ribbons. For any
$\mu\in\Y$, we know that the set of shapes covering~$\mu$ in~$(\Y,\precr)$ has
$r$~more elements than the set of shapes covered by~$\mu$. A key observation,
made in~\ref{Shimozono White}, is that these sets, labelled by the height of
the \r-ribbon involved in the covering relation, have a very regular
structure. There exists for instance, for every possible height~$h\in\set{r}$,
at least one $\kappa\succr\mu$ with $\hgt(\kappa/\mu)=h$, a fact that does not
seem immediately obvious (except for $h=0$ and $h=r-1$). The heads of
\r-ribbons of the form $\mu/\\$ or $\kappa/\mu$ all lie on distinct diagonals,
so there is a natural total ordering on the set of those ribbons. For any
diagonal~$d_k$, the value $\sum_{j=1}^{r-1}\edge\mu_{k-j}$ gives the height of
the ribbon with its head on that diagonal \emph{if it exists}, and that sum
changes by steps at most~$1$ as~$k$ varies. The following somewhat surprising
proposition states that, even though many diagonals have no associated
ribbons, the heights of the ribbons that are present still change by steps at
most~$1$ as the diagonals are traversed, with even some additional
constraints.

\proclaim Proposition. \ribbonorderprop
Let $\mu\in\Y$, and let $\xi_0,\xi_1$ be \r-ribbons, both either of the form
$\mu/\\$ or~$\kappa/\mu$, such that the diagonals $d_s,d_t$ containing their
heads satisfy $s<t$, and no diagonal~$d_i$ with $s<i<t$ contains the head of
any \r-ribbon of the form $\mu/\\$ or~$\kappa/\mu$. Then one of the following
cases applies:
\statitem both $\xi_0$ and~$\xi_1$ are of the form~$\kappa/\mu$, and
          $\hgt(\xi_1)=\hgt(\xi_0)-1$;
\statitem there is one~$\xi_i$ of each of the forms $\mu/\\$ and~$\kappa/\mu$,
          and $\hgt(\xi_1)=\hgt(\xi_0)$;
\statitem both $\xi_0$ and~$\xi_1$ are of the form~$\mu/\\$, and
          $\hgt(\xi_1)=\hgt(\xi_0)+1$.
\enditem
Moreover the \r-ribbons~$\xi_\ll,\xi_\gg$ of one of the given forms with their
heads on diagonals~$d_k$ with $k$ minimal respectively maximal, are both of
the form~$\kappa/\mu$, while $\hgt(\xi_\ll)=r-1$ and $\hgt(\xi_\gg)=0$.

\proof
Instead of considering sums of~$r-1$ consecutive bits of~$\edge\mu$, consider
sums of~$r$ consecutive bits. As we shift the range of summation upwards so as
to include $\edge\mu_k$ while dropping $\edge\mu_{k-r}$, the sum changes if
and only if there is an \r-ribbon of one of the given forms with its head on
diagonal~$d_k$, with an increase of the sum for ribbons of the form~$\mu/\\$,
and a decrease of the sum for ribbons of the form~$\kappa/\mu$; in either case
the smaller of the two sums gives the height of the ribbon. Considering two
consecutive changes of the sum gives rise to the stated three cases. The final
statement corresponds to the fact that the sum, which takes values
in~$\{0,\ldots,r\}$, tends to~$r$ as $k\to-\infty$, and to~$0$ as
$k\to+\infty$.
\QED

\proclaim Corollary. \interleavecorr
Let $\mu\in\Y$ and $h\in\set{r}$ be fixed. The \r-ribbons of the form
$\xi=\kappa/\mu$ with $\hgt(\xi)=h$ and those of the form $\xi=\mu/\\$ with
$\hgt(\xi)=h$ are perfectly interleaved when ordered from bottom left to top
right (by the diagonal of their head), with ribbons of the first kind at both
ends. In particular, the number of ribbons of the first kind exceeds the
number of those of the second kind by~$1$, exactly.
\QED

The easiest way to understand the corollary is by the same considerations as
in the proof of the proposition: each time the sum over $r$ consecutive bits
of~$\edge\mu$ descends from the range $\{h+1,\ldots,r\}$ to its complement
$\{0,\ldots,h\}$, there is an \r-ribbon of the form $\xi=\kappa/\mu$ with
$\hgt(\xi)=h$; each time it rises back to $\{h+1,\ldots,r\}$ there is an
\r-ribbon of the form $\xi=\mu/\\$ with $\hgt(\xi)=h$. There is also a nice
visual presentation of this proof. Draw the diagram of~$\mu$ and its boundary,
and superimpose a copy of the boundary shifted down by $h+{1\over2}$ units and
leftwards by $r-h-{1\over2}$ units. The boundary and its copy can only cross
in the middle of edge segments, and whenever they do, these two segments
correspond to an \r-ribbon of height~$h$ that can either be added to~$\mu$
(when the copy passes the boundary from the ``inside'', the side of the
diagram, to the outside) or removed from~$\mu$ (when the copy passes back to
the inside). Clearly these types of crossings alternate as one traverses the
boundary, and due to the direction of the shift, one eventually passes from
the inside to the outside. We illustrate this process for $\mu=(6,6,6,4,4,1)$,
$r=4$, and~$h=2$; for those values two ribbons can be added to~$\mu$, and one
removed.
\bigdisplay
 
 \epsfbox{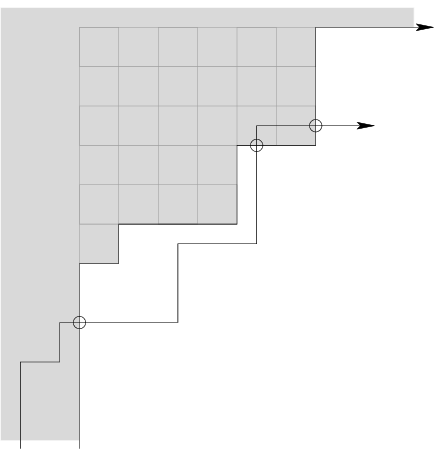}%
\qquad
 \epsfbox{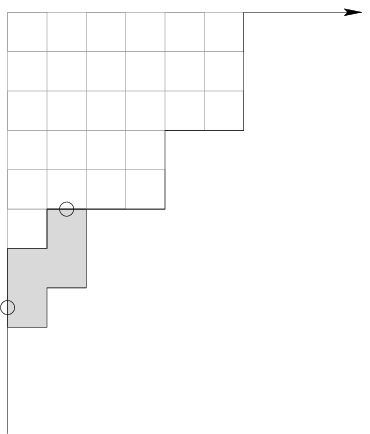}%
\qquad
 \epsfbox{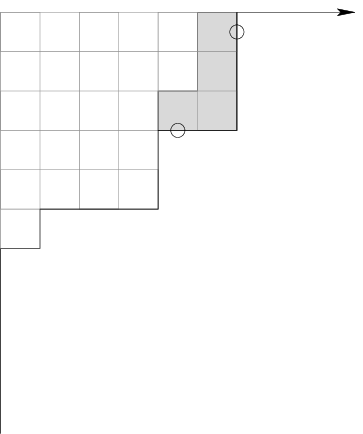}%
\qquad
 \epsfbox{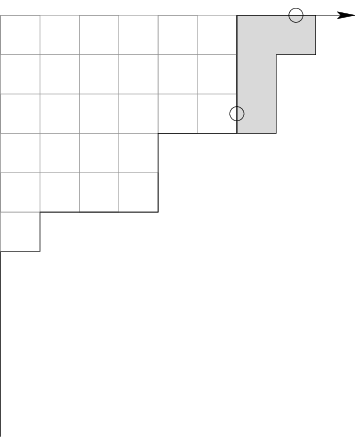}%
$$
It is remarkable that whether we group the \r-ribbons that can be added to or
removed from~$\mu$ by their position class or by their height, they will be
interleaved within each group in exactly the same way. Like for position
classes, this property allows us to define an \r-correspondence
for~$(\Y,\precr)$ in two ways (moving to the top right or to the bottom left).
With the choice of the \r-correspondence moving to the top right (recall that
this direction refers to the extraction process), such an \r-correspondence
was first described in \ref{Shimozono White, propositions~19,20}, and we shall
call it the Shimozono-White \r-correspondence~$\bSW$. For the corresponding
Schensted insertion algorithm, bumping of \r-ribbons is towards the bottom
left, like bumping in the Robinson-Schensted algorithm. Explicitly, if
$\\\precr\mu$ and $\hgt(\mu/\\)=h$, then $\\=\bSW_\mu(\kappa)$ where
$\kappa/\mu$ is an \r-ribbon of height~$h$ that can be added to~$\mu$, and the
first among such ribbons to the bottom left of~$\mu/\\$. Also
$e_h=\bSW_\mu(\kappa)$ where $\kappa/\mu$ is the top-rightmost \r-ribbon of
height~$h$ that can be added to~$\mu$. The corollary guarantees that this
properly describes a unique \r-correspondence. The opposite choice, obtained
by replacing ``bottom left'' by ``top right'' in the description and vice
versa, is also possible; we shall call the \r-correspondence so defined the
transposed Shimozono-White \r-correspondence, and denote it by~$\bWS$ (which
is more convenient than $(\bSW)\tr$, although some might have preferred
$b^{\rm NE}$).

These \r-correspondences, and the Schensted-growths defined by them, behave
well with respect to height. Let $b=\bSW$ or~$b=\bWS$, and consider a square
$\\~\mu\choose\nu~\kappa$ of a Schensted-growth for~$b$, with matrix
entry~$a\in\{0\}\union e_\set{r}$. We assume that $\kappa\notin\{\mu,\nu\}$
(for otherwise one has $\{\\,\kappa\}=\{\mu,\nu\}$ and~$a=0$, which causes our
conclusions below to hold trivially). We have seen that if $\mu\neq\nu$, then
one always has $\hgt(\kappa/\mu)+\hgt(\kappa/\nu)=\hgt(\nu/\\)+\hgt(\mu/\\)$
and~$a=0$. In the case that $\mu=\nu$ and $b_\mu(\kappa)=\\\in\Y$, one has the
stronger set of conditions
$\hgt(\kappa/\mu)=\hgt(\kappa/\nu)=\hgt(\nu/\\)=\hgt(\mu/\\)$ and~$a=0$.
Finally if $\mu=\nu$ and $b_\mu(\kappa)=a\in e_\set{r}$, then $\\=\mu=\nu$ and
$a=e_h$ where $h=\hgt(\kappa/\mu)=\hgt(\kappa/\nu)$. Therefore, with the
conventions that $\n(e_h)=h$ (the ``colour'' of the entry of a coloured
permutation, which we are really interpreting as an integer here) and
$\hgt(\\/\\)=\n(0)=0$, we may conclude that in all cases one has the identity
$\hgt(\kappa/\mu)+\hgt(\kappa/\nu)=\hgt(\nu/\\)+\hgt(\mu/\\)+2\n(a)$. This can
be seen as a ``law of conservation'': if for lattice paths from $(0,n)$ to
$(n,0)$ in the grid of the Schensted-growth one takes the sum of the heights
of all ribbons formed by adjacent shapes on the path, and adds to it twice the
sum of the colours of all matrix entries of the Schensted-growth south-east of
the path (i.e., on the side away from the origin), then this sum will be
invariant under moving the path across a square. Consequently, the sum is the
same for all such paths; in particular (since there are no ribbons for the
path going via the origin) twice the sum of the colours of the entries of the
coloured permutation associated to the Schensted-growth will equal the sum of
the heights of the ribbons in the $P$-symbol and the $Q$-symbol associated to
it. Since half the sum of the heights of the ribbons in a standard
\r-ribbon tableau is called its~``spin'' (see \ref{Carre Leclerc} and
\ref{Lascoux Leclerc Thibon}; the motivation for the terminology is not quite
clear), the computation of Schensted correspondence using the Shimozono-White
\r-correspondence can be called a ``colour-to-spin Schensted algorithm''. We
note in passing that since $P,Q$ may be any pair of standard \r-ribbon
tableaux of the same shape, and the sum of the heights of their ribbons is
always even, the spins of \r-ribbon tableaux of a given shape must either all
lie in~$\N$ or all lie in~$\N+{1\over2}$, a well known fact that is related to
the notion of \r-sign of a shape.

To formulate the enumerative consequences of the existence of such height
respecting correspondences, we give a ``$q$-analogue'' of equation~(\DUcomm)
that holds for the current situation. With $q$ an indeterminate, we consider
endomorphisms $D^r_q,U^r_q$ of the free $\Z[q]$-module $\Z[q]\Y$ on the
set~$\Y$, defined by their action on basis elements
$U^r_q\:\\\mapsto\sum_{\mu\succr\\}q^{\hgt(\mu/\\)}\mu$ and
$D^r_q\:\\\mapsto\sum_{\mu\precr\\}q^{\hgt(\\/\mu)}\mu$. Then the corollary,
together with the conservation of total height observed for ribbons related to
shapes with a common successor and a common predecessor in~$(\Y,\precr)$,
leads to the identity
$$
  D^r_q\after U^r_q=U^r_q\after D^r_q+r_{q^2}\Id,
\ttex{where $r_{q^2}=\sum_{i\in\set{r}}q^{2i}={1-q^{2r}\over1-q^2}\in\Z[q^2]$.}
\label(\qDUcomm)
$$
From it one derives for any \r-core~$\gamma$ (by writing
$\scal<\gamma|(D^r_q)^n((U^r_q)^n(\gamma))>$ and expanding $(D^r_q)^n\after
(U^r_q)^n$ into terms $(U^r_q)^i\after (D^r_q)^i$) the identity:
$$
  \scal<(U^r_q)^n(\gamma)|(U^r_q)^n(\gamma)>=n!(r_{q^2})^n
.\label(\qRSid)
$$
We denote by $\RH(\gamma)$ the connected component of~$(\Y,\leqr)$ containing
the \r-core~$\gamma$ (an \r-rim hook lattice), and by
$F^{\\/\gamma}_{q^{1/2}}\in\Z[q^{1\over2}]$ (using a square root~$q^{1\over2}$
of~$q$) the spin generating series $\sum_Tq^{\spin(T)}$ of the set of standard
\r-ribbon tableaux~$T$ of shape~$\\/\gamma$. We can then write~(\qRSid) in a
form that directly matches the colour-to-spin Schensted correspondence:
$$
  \sum_{\\\in\RH(\gamma)}(F^{\\/\gamma}_{q^{1/2}})^2
 = n!(r_q)^n
 = n! \left(1-q^r\over1-q\right)^n
,\nn
$$
whose right hand side can be interpreted as the sum over all \r-coloured
permutations~$\sigma$ of $q^{\n(\sigma)}$, where $\n(\sigma)$ is the sum of
the colours of the entries of~$\sigma$.

We have finally arrived at the point where we can state our main goal: to find
a Knuth correspondence for~$(\Y,\leqhr,\N^r)$ that satisfies a similar
conservation law involving heights. If one identifies $e_\set{r}$ with the
``standard basis'' of~$\N^r$, this Knuth correspondence should in particular
reduce to a ``colour-to-spin'' Schensted correspondence when all matrix
entries~$a$ satisfy $|a|\leq1$ (i.e., $a\in\{0\}\union e_\set{r}$). To make
our goal precise, one needs to define the spin of a horizontal \r-ribbon
strip, and the contribution to the spin of a matrix entry (for semistandard
\r-ribbon tableaux and complete matrices these quantities will then be defined
by summing over their constituent parts). By definition $\spin(\mu/\\)$ is
equal to the spin of the standardisation of~$\mu/\\$ (cf.~\ref{Lascoux Leclerc
Thibon}); note that it is here that the notion of standardisation of a
horizontal \r-ribbon strip becomes essential, since other interpolations into
a standard \r-ribbon tableau will in general exist, but these will have
different spins (the standardisation in fact achieves the maximal possible
spin, much to the satisfaction of the spin doctors). The contribution of a
matrix entry~$a\in\N^r$ to the spin extends by linearity the fact $e_h$
contributes~$h$, so it will be $\sum_{i\in\set{r}}ia_i$ (one might call this
the ``total colour'' of the entry~$a$, but we are no longer dealing with
coloured permutations, and we wish to drop the colourful terminology). So we
want $(a,\\)=b_{\mu,\nu}(\kappa)$ to be defined in such a way that
$$
  \spin(\kappa/\mu)+\spin(\kappa/\nu)=\spin(\mu/\\)+\spin(\nu/\\)+\n(a),
\ttext{where $\n(a)=\sum_{i\in\set{r}}ia_i$,}
\label(\spinbalance)
$$
always holds. To formulate an enumerative identity that is required for the
existence of such a shape datum, let us introduce $q^{1\over2}$-analogues of
the generating series $U_X$, $D_Y$:
$$
  U^r_{q^{1/2},X}(\\)=\sum_{\mu\geqhr\\}q^{\spin(\mu/\\)}X^{\rwt\mu\\}\mu
,\qquad
  D^r_{q^{1/2},Y}(\\)=\sum_{\mu\leqhr\\}q^{\spin(\\/\mu)}Y^{\rwt\\\mu}\mu
,
\label(\spingenseq)
$$
(one has $U^r_{q^{1/2},X},D^r_{q^{1/2},Y}\in\End(\Z[q^{1\over2}]\Y)[[X,Y]]$),
and a $q$-analogue of $F_{\N^r}$:
$$
  F_{\N^r}(q,T)=\sum_{a\in\N^r}q^{\n(a)}T^{|a|}
  =\prod_{i\in\set{r}}\left(\sum_{a_i\in\N}q^{ia_i}T^{a_i}\right)
  =\prod_{i\in\set{r}}{1\over1-q^iT}
  \in\Z[q][[T]]
.\nn
$$
Then the indicated identity that will be bijectively proved by our main
construction is analogous to~(\sdid), with the expression for $F_{\N^r}(q,XY)$
substituted:
$$
  D^r_{q^{1/2},Y}\after U^r_{q^{1/2},X}
  =(U^r_{q^{1/2},X}\after D^r_{q^{1/2},Y})\prod_{i\in\set{r}}{1\over1-q^iXY}
\,.\label(\qbscid)
$$
Taking the coefficient of~$XY$ gives an identity equivalent to~(\qDUcomm).
What we have seen so far only proves that special case, but (\qbscid) is known
to hold: it can be derived from the commutation relation \ref{Lascoux Leclerc
Thibon,~(21)} proved in~\ref{Kashiwara Miwa Stern}, for operators $B_i$ that
satisfy $\exp\(\sum_{i>0}{B_{-i}\over{i}}X^i\)=U^r_{-q\inv,X}$ and
$\exp\(\sum_{i>0}{B_i\over{i}}Y^i\)=D^r_{-q\inv,Y}$.

We have already hinted at how we shall bijectively prove~(\qbscid), and indeed
indicated the algorithm defining the bijection, but let us here also mention
an obvious idea that does not work. Like for the RSK~correspondence, one may
take standardisations of $\kappa/\mu$ and $\kappa/\nu$, and compute the
Schensted-growth from it for $\bSW$ or $\bWS$, so as to find a shape~$\\$ at
the top left corner, and matrix entries that sum up to a value~$a\in\N^r$. The
fundamental problem with this approach is that in the ribbon setting the
property of being a standardised horizontal strip is not preserved across
Schensted-growths. We have already observed this for the
\r-correspondences derived from the \r-quotient maps, and the situation is
no better for the Shimozono-White \r-correspondence or its transpose.
Informally speaking, when $b_\mu$ is used to find a ribbon $\mu/\\$ matching
the height of~$\kappa/\mu$, that ribbon may have its head on a quite distant
diagonal, causing it to be out of order with respect to other ribbons of the
standardisation. Moreover one cannot ``shuffle the ribbons'' to put them back
into order: even if $\mu/\\$ and $\nu/\\$ should happen to be horizontal
\r-ribbon strips, unless the tableaux obtained from the
Schensted-growth \emph{are} their respective standardisations, their spins
will not have the proper values. Here are very simple examples with $r=2$
where the Schensted-growth does not give a useful result; the first one
uses~$\bSW$, the second one~$\bWS$.
\bigdisplay 
\def\young(#1){\vtop{\kern0pt \smallsquares\Young(#1)}\hfill}
\def~{\vtop{\kern0pt \hbox{$\circ$}}}
\def\:{\lower15pt \hbox{$\longrightarrow$}}
\def\rt#1{\vtop{\kern0pt
                \hbox{$\smallsquares{\rtab2 #1 \edgeseq4,0;1111000 }$}}\hfill}
\matrix { \young(|)   & \young(,|,)  & \young(,|,|,)
  &\: &  \rt{\ribbon1,1:0;1 \ribbon2,0:1;0} \cr
\noalign{\medskip}
  	  \young(,|,) & \young(,|,||)& \young(,|,|,|,)
  &\: &  \rt{\ribbon3,0:0;1 \ribbon3,1:1;1 \edgeseq2,2;11}\cr}
\kern0.15\hsize
\def\:{\lower10pt \hbox{$\longrightarrow$}}
\def\rt#1{\vtop{\kern0pt
                \hbox{$\smallsquares{\rtab2 #1 \edgeseq3,0;11100000 }$}}\hfill}
\matrix { ~         & \young(,)  & \young(,|,)
  &\: &  \rt{\ribbon0,0:0;0 \ribbon1,0:1;0} \cr
        \noalign{\medskip}
  	  \young(|) & \young(,|,)& \young(,,,|,)
  &\: &  \rt{\ribbon1,1:0;1 \ribbon0,2:1;0 \edgeseq2,0;0}\cr}
$$
Again the difficulty is that the ribbon being bumped (the one labelled~$0$ on
the left and the one labelled~$1$ on the right) bypasses another ribbon that
stays on its diagonal. The only way such bypassing could be avoided, is if any
ribbon that risks being overtaken would be bumped instead, which is what
always happens in the case~$r=1$. But for $r>1$, the ribbon bypassed may not
even be candidate for bumping because it is ``at the wrong side of the
boundary'': in the first example in the extraction direction, and in the
second example in the insertion direction, the ribbon being overtaken has
already reached its destination. Such considerations seem to indicate that it
is impossible to define \emph{any} \r-correspondence for~$(\Y,\leqr)$ whose
Schensted-correspondence will preserve standardisations of horizontal
\r-ribbon strips (regardless of concerns about spins), with one notable
exception: the 2-correspondence $\bBVG$ for the Barbasch-Vogan-Garfinkle (or
hyperoctahedral Robinson-Schensted) correspondence described in \ref{van
Leeuwen festschrift, definition~4.2.1} does preserve such standardisations. It
achieves this by choosing the direction of bumping based on the \emph{form} of
the 2-ribbon involved, setting $\bBVG_\mu(\kappa)=\bSW_\mu(\kappa)$ if
$\form(\kappa/\mu)=0$ (a horizontal domino), and
$\bBVG_\mu(\kappa)=\bWS_\mu(\kappa)$ if $\form(\kappa/\mu)=1$ (a vertical
domino); thus the sum of $2$ consecutive bits in the proof of
proposition~\ribbonorderprop\ assumes an extremal value ($0$ or~$2$) at the
start of the search, ensuring that the first change (if any) of its value will
give the domino sought for, without any other dominoes being skipped. The
domino found will also have the same form as the one bumped, which is
essential for bijectivity, and moreover makes the Schensted correspondence
spin preserving. A Knuth correspondence for~$(\Y,\leqh_2,\N^2)$ defined via
standardisation from~$\bBVG$ is described (in different terms) in
\ref{Shimozono White domino}.

In spite of the mentioned difficulty, a method that is an adaptation of the
standardisation idea is given in \ref{Shimozono White}, which handles the
``half semistandard'' case, where $Q$-symbol is still standard, and the matrix
entries are in~$0\union e_\set{r}$. Their solution to the problem is
ingenious: the Schensted-growth is used for trivial steps not
involving~$\bSW$, but if some step would involve~$\bSW$ (there is at most
one), then instead they let the ribbon in question disappear, and show a
ribbon of the same height can be reincarnated \emph{somewhere} among the
images of the other ribbons to produce a proper standardised horizontal strip.
Moreover the same method in the opposite direction will reconstruct the ribbon
at its original location, which is of course essential for bijectivity. That
this is possible is far from trivial, and it implies an enumerative identity
stronger than~(\qDUcomm). That provides evidence that our goal of bijectively
proving~(\qbscid) might be realistic (and it provided our initial motivation
to make an attempt), but the method seems ill suited for adaptation to the
``full semistandard'' case. When the Schensted-growth would invoke~$\bSW$ for
more than one ribbon, the corresponding searches for appropriate replacement
ribbons can interfere with each other, and even in cases where all of them
succeed fairly easily, the relative order of the ribbons found may have
changed, which makes a step-by-step inverse virtually inconceivable.

\newsection A spin preserving shape datum.
\labelsec\mainSec

In this section we shall describe a shape datum for~$(\Y,\leqhr,\N^r)$ that
satisfies (\spinbalance), and which therefore bijectively proves~(\qbscid). An
essential point is that the Schensted correspondence that it generalises is
not the Shimozono-White correspondence but its transpose, and similarly that
for~$r=1$ it reduces to the Burge correspondence rather than to the
RSK~correspondence. We shall first give an informal description, in a form
similar to the algorithm given at the end of~\Sec\factorsec, without yet
proving that this method actually works as it should. After that we shall give
a more formal description, in a form that is more suited to the necessary
proofs, which will then in fact be provided. Notably we shall prove that the
correspondence has an inverse (given by a very similar procedure) and that
equation~(\spinbalance) is indeed always satisfied.

\subsection Informal description, examples.

Let shapes $\mu,\nu$ be given, in the same connected component
of~$(\Y,\leqr)$. Our task is to establish a bijection between on hand the
shapes~$\kappa$ with $\mu\leqhr\kappa\geqhr\nu$ and on the other hand the
pairs of a shape~$\\$ with $\mu\geqhr\\\leqhr\nu$ and an \r-tuple $a\in\N^r$,
such that equations (\rankeq)~and~(\spinbalance) are satisfied in all cases.
Although edge sequences are a central to our considerations, we shall avoid
referring to them in our informal description; however, we shall often compare
\r-ribbons by the diagonal containing the head of the ribbon, so to alleviate
the terminology a bit, we shall say that a diagonal contains a ribbon, and
that the ribbon is on the diagonal, when in fact the diagonal contains the
head of the ribbon. We start with the extraction direction, i.e., the
determination of~$(a,\\)$ given~$\kappa$; we shall assume that the
standardisations of the horizontal \r-ribbon strips $\kappa/\mu$ and
$\kappa/\nu$ have been computed. Our description will be an algorithm that
treats $\\$ and~$a$ as variables, which are initialised as $\\:=\kappa$ and
$a:=0\in\N^r$, and which will contain the desired values at termination.

One starts at the rightmost diagonal~$d_k$ containing a ribbon of at least one
of the standardisations of the strips $\kappa/\mu$ and $\kappa/\nu$; should
both strips be empty, then one terminates immediately. Each time a
diagonal~$d_k$ has been processed as indicated below, one checks if the
diagram of~$\\$ has an empty intersection with any of the diagonals $d_j$ for
$k-r\leq j<k$. If it has, one terminates (there are no \r-ribbons left to
remove from~$\\$); if not, one continues to consider the diagonal~$d_{k-1}$.
No action is required for~$d_k$ unless there exists an \r-ribbon on~$d_k$ that
can be removed from~$\\$. If so, let $\\/\\'$ be that \r-ribbon, and put
$h=\hgt(\\/\\')$. If both the standardisations of the strips $\kappa/\mu$ and
$\kappa/\nu$ have \r-ribbons on~$d_k$, modify~$a$ by setting~$a_h:=a_h+1$. If
at least one of those standardisations has an \r-ribbon on~$d_k$, put
$\\:=\\'$, and processing of~$d_k$ is completed. In the remaining case
(neither of the standardisations has a ribbon on~$d_k$), processing is also
completed (without action) if $a_h=0$; otherwise one sets~$a_h:=a_h-1$ and
$\\:=\\'$, completing the processing for~$d_k$. After the algorithm
terminates, $\kappa/\\$ will be a horizontal \r-ribbon strip whose
standardisation is given by the intermediate values of~$\\$, in other words,
it has ribbons on those diagonals for which the variable~$\\$ was modified.
Then $\mu/\\$ will also be a horizontal \r-ribbon strip, and its
standardisation has ribbons on those diagonals for which~$\\$ was modified and
that do not contain a ribbon of the standardisation of~$\kappa/\mu$; the
situation is similar for $\nu/\\$ and~$\kappa/\nu$.

The reverse procedure is quite similar, but this time the variables $\\$
and~$a$ are initialised from the given shape and matrix entry, and one starts
with $k=-\\\tr_0$, so that the diagonal~$d_k$ contains the leftmost \r-ribbon
that can be added to~$\\$, whose head is the square~$(0,\\\tr_0)$, the first
one below the leftmost column of~$\\$. One proceeds from there for successive
values of~$k$, until $a=0$ holds and there are no more ribbons of the
standardisations of the strips $\mu/\\$ and $\nu/\\$ any diagonal~$d_i$ with
$i\geq k$. Action is needed for~$d_k$ only if an \r-ribbon exists on~$d_k$
that can be added to~$\\$. If so, let $\\'/\\$ be that \r-ribbon, and put
$h=\hgt(\\'/\\)$. If both the standardisations of the strips $\mu/\\$ and
$\nu/\\$ have \r-ribbons on~$d_k$, set~$a_h:=a_h+1$, and then if at least one
of them has an \r-ribbon on~$d_k$ set~$\\:=\\'$, completing the
processing of~$d_k$. In the remaining case processing is also completed
(without action) if $a_h=0$; otherwise the processing of~$d_k$ consists of
setting~$a_h:=a_h-1$ and $\\:=\\'$. This second procedure exactly retraces the
steps of the first one, but this is not as obvious as it might seem at first
glance. Notably, encountering a diagonal where an \r-ribbon of height~$h$ can
be added to~$\\$ triggers no action during the first procedure, but during the
reverse procedure this circumstance \emph{will} cause action, unless $a_h=0$.
We shall see below that in these cases $a_h=0$ always holds, which resolves
the mystery.

Now let us illustrate these computations in an example. In order to exercise
the different cases that can arise, we need rather large diagrams, in
particular horizontally. Our example is for~$r=4$; we take
$\mu=(16,15,15,5,4)$, $\nu=(14,14,14,9,4)$, and
$\kappa=(17,17,16,13,9,5,1,1)$. Then the standardisations of the horizontal
4-ribbon strips involved are as displayed below. We have labelled
individual ribbons to facilitate our discussion of the procedure,
giving the same label to ribbons whose diagram is identical (even if they are
not identical as ribbons, i.e., the step in~$(\Y,\precr)$ in the two
standardisations may differ; this should not cause any confusion).
$$
\cramp e\e9
\def\f{\raise1pt \hbox{$\,f$}}
  \kappa/\mu:~
  {\smallsquares\rtab 4
   \ribbon 2,15:a;101 \ribbon 3,9:c;000 \ribbon 4,6:d;001 \ribbon 4,5:e;100
   \ribbon 5,2:\f;001 \ribbon 7,0:g;110
   \edgeseq5,0;111110000000000000000
   \edgeseq3,13;00
  }
,\quad
  \kappa/\nu:~
  {\smallsquares\rtab 4
   \ribbon 2,15:a;101 \ribbon 2,14:b;110 \ribbon 3,9:c;000 \ribbon 4,5:\e;000
   \ribbon 5,2:\f;001 \ribbon 7,0:g;110
   \edgeseq5,0;1111100000000000000
   \edgeseq3,13;0
  }
$$
The procedure starts with $\\:=\kappa$ and $a:=(0,0,0,0)$, on the diagonal
containing the ribbon~$a$. Since this ribbon has height~$2$ and occurs both in
$\kappa/\mu$ and in~$\kappa/\nu$, we set~$a_2:=1$, and then remove the ribbon
from~$\\$. The next diagonal (to its left) contains~$b$, which only occurs
in~$\kappa/\nu$; therefore $a_2=1$ remains unchanged, and $b$ is removed
from~$\\$, which becomes $(14,14,14,13,9,5,1,1)$. No 4-ribbon can be removed
from~$\\$ on the next two diagonals, but this is possible on the diagonal
after that. This ribbon, which we shall call~$u$, has $\form(u)=(1,0,1)$ and
$\hgt(u)=2$; its diagonal contains no ribbons either in $\kappa/\mu$ or
in~$\kappa/\nu$ and so since $a_2=1$, we decrement $a_2$ back to~$0$ and
remove~$u$ from~$\\$, which now becomes $(14,13,12,12,9,5,1,1)$. On the next
diagonal another 4-ribbon (of form $(1,1,0)$) can be removed, but since there
are still no ribbons either in $\kappa/\mu$ or in~$\kappa/\nu$ on this
diagonal, and by now $a_2=0$, no action is taken here. Skipping another
diagonal, on which no 4-ribbon can be removed from~$\\$, we come to the
diagonal of~$c$. There a 4-ribbon can be removed from~$\\$, but it is not~$c$:
it has form~$(0,0,1)$ and height~$1$ and we shall call it~$c'$. Since~$c$
occurs both in $\kappa/\mu$ and in~$\kappa/\nu$, we set~$a_1:=1$, and
removing~$c'$ from~$\\$ it becomes $(14,13,11,9,9,5,1,1)$. Three diagonals
follow where no 4-ribbons can be removed from~$\\$ (nor could any be added)
and we arrive at the diagonal of~$d$. It occurs only in~$\kappa/\mu$ so we
remove it from~$\\$ and arrive at the diagonal of~$e$ and~$e'$. The 4-ribbon
that can be removed from~$\\$ on this diagonal has the form of~$e$ and
therefore height~$1$; because of the occurrences of $e$~and~$e'$ we raise
$a_1:=2$ and remove~$e$ from~$\\$ which now becomes $(14,13,11,5,5,5,1,1)$.
There follow a diagonal on which a 4-ribbon could be added to~$\\$ and one on
which a 4-ribbon could be removed, but since the latter has height~$2$ while
$a_2=0$ we do nothing here. Then comes~$f$, which can be removed from~$\\$;
since it has height~$1$ and occurs both in $\kappa/\mu$ and in~$\kappa/\nu$ we
further raise $a_1:=3$ and remove~$f$ from~$\\$. The next diagonal allows
removal of a 4-ribbon~$v$ of height~$1$ from~$\\$ but has no ribbons in
$\kappa/\mu$ or in~$\kappa/\nu$; therefore we lower $a_1:=2$ and remove the
ribbon from~$\\$ leaving $(14,13,11,5,1,1,1,1)$. The next two diagonals would
allow adding 4-ribbons to~$\\$ so they are skipped; we come at the diagonal
of~$g$, on which the 4-ribbon~$g'$ that can be removed from~$\\$ is vertical
(height~$3$), so we raise $a_3:=1$ and remove~$g'$ from~$\\$. No squares
of~$\\$ remain on the current diagonal (let alone to its left) so we terminate
with $\\=(14,13,11,5)$ and $a=(0,2,0,1)$.

We summarise the result by displaying the standardisations of $\kappa/\\$,
$\mu/\\$, and $\nu/\\$. For the two latter strips, each of which has only
three ribbons, we have also included (in green when colours are available)
standardisations of their complementary strips in $\kappa/\\$, to emphasise
the relative locations of the ribbons.
$$
\def\f{\raise1pt \hbox{$\,f$}}
\matrix{ \kappa/\\:~ & \green{\kappa/}\mu/\\:~ & \green{\kappa/}\nu/\\:~ \cr
  {\cramp c\c9 \cramp g\g8
   \smallsquares\rtab 4
   \ribbon 2,15:a;101 \ribbon 2,14:b;110 \ribbon 3,12:u;101
   \ribbon 3,9:\c;001 \ribbon 4,6:d;001 \ribbon 4,5:e;100
   \ribbon 5,2:\f;001 \ribbon 5,1:v;100
   \ribbon 7,0:\g;111
   \edgeseq4,0;111100000000000000
  }
&  
  {\cramp u\u9 \cramp v\v8
   \smallsquares\rtab 4
   \green{%
   \ribbon 2,15:a;101 \ribbon 3,9:c;000 \ribbon 4,6:d;001 \ribbon 4,5:e;100
   \ribbon 5,2:\f;001 \ribbon 7,0:g;110
   }%
   \ribbon 2,14:b;110 \ribbon 2,11:\u;001
   \ribbon 4,0:\v;000
   \edgeseq4,0;111100000000000000
  }
&
  {\cramp u\u9 \cramp v\v8 \cramp d\d7 \cramp e\e6
   \smallsquares\rtab 4
   \green{%
   \ribbon 2,15:a;101 \ribbon 2,14:b;110 \ribbon 3,9:c;000 \ribbon 4,5:\e;000
   \ribbon 5,2:\f;001 \ribbon 7,0:g;110
   }%
   \ribbon 2,11:\u;001 \ribbon 3,5:\d;000
   \ribbon 4,0:\v;000 
   \edgeseq4,0;111100000000000000
  }
\cr}
$$
We observe that $|\kappa/\mu|_4-|\nu/\\|_4=6-3=3=|a|$ shows that~(\rankeq)
indeed holds in this case, and similarly
$\spin(\kappa/\mu)+\spin(\kappa/\nu)-\spin(\mu/\\)-\spin(\nu/\\)
={7\over2}+{7\over2}-{3\over2}-{1\over2}=5=\n(a)$ shows that~(\spinbalance)
holds. From the above description, the one for the opposite process (computing
$\kappa$ from the (final) values of~$\\$ and~$a$) can be obtained by a
reversal of its steps, replacing each action by its inverse; we leave it to
the reader to work this out. We do note however that those diagonals where
ribbons could be added to~$\\$ (a fact not really relevant to the initial
computation) now get considered for action, but invariably get skipped
nonetheless because the relevant component $a_h$ is zero. This is for instance
the case just after adding~$g'$ to~$\\$ and decrementing~$a_3$ to~$0$ on the
first diagonal considered: two occasions follow to add ribbons of height~$3$
and~$2$, respectively, but neither is used because now $a_3=0$, and~$a_2=0$.
Similarly, after adding ribbons~$v$ (on the diagonal of~$v'$) and~$f$ (because
this time $a_1=3$ can be decremented) and skipping a diagonal where a ribbon
could be removed, an occasion to add a ribbons of height~$2$ is not used
because one (still) has $a_2=0$; the next occasion involves ribbon~$e$ of
height~$1$, and it is used because $a_1=2$ there.

Let us use this example to make an informal observation: the algorithm turns
out to be surprisingly sensitive the input. Suppose we change the data of the
example slightly be including ribbon~$b$ into~$\kappa/\mu$. Then when this
ribbon is reached, instead of leaving~$a_2=1$, the initial procedure will
increment~$a_2$. One might imagine that this just causes $a_2$ to be
incremented in the final result, but in fact the whole process completely
changes. One may check that the results will be $\\=(14,10,6,1)$ and
$a=(0,0,0,1)$, so rather than an increase of~$a_2$ we see a decrease of~$a_1$.
Indeed the latter component never even gets to be incremented, which is
related to the fact that $\kappa/\\$ now has no ribbons of height less
than~$2$ in its standardisation:
$$
  \kappa/\\:~
  {\smallsquares\rtab 4
   \ribbon 2,15:;101 \ribbon 2,14:;110 \ribbon 3,12:;101 \ribbon3,11:;110
   \ribbon 3,9:;011 \ribbon 4,8:;110 \ribbon 4,6:;011 \ribbon 4,5:;101
   \ribbon 5,3:;011 \ribbon 5,2:;101 \ribbon 5,1:;110
   \ribbon 7,0:;111
   \edgeseq4,0;111100000000000000
  }
$$
Here is how the statistics change: $|\kappa/\mu|_4-|\nu/\\|_4=7-6=1=|a|
=6-5=|\kappa/\nu|_4-|\mu/\\|_4$ and
$\spin(\kappa/\mu)+\spin(\kappa/\nu)-\spin(\mu/\\)-\spin(\nu/\\)
={9\over2}+{7\over2}-{4\over2}-{6\over2}=3=\n(a)$.

\subsection Formal description.

Our formal description will be directly in terms of edge sequences, since all
the conditions and actions used are more easily expressed in terms of these.
At the same time this allows us to generalise by considering arbitrary doubly
infinite sequences of bits; the dependence on the limiting behaviour of these
sequences will be made explicit. We therefore define a \newterm{bit sequence}
to be any function $s\:\Z\to\{0,1\}$; its value at~$i$ will be denoted
by~$s_i$. Moreover, our description will be static, in the sense that there
are no variables that are modified in time. Instead all values used are
represented in a single structure, and the algorithms translate into relations
between components of this structure that allow the whole to be recovered from
partial information, much like growth diagrams can replace the algorithmic
descriptions of Schensted and Knuth correspondences. An additional advantage
of this method of description is that it simultaneously describes the forward
and the reverse algorithms.

A first notion to formalise is that of horizontal \r-ribbon strips. We have
characterised the relation~$\\\leqhr\mu$ by the possibility to change $\\$
into~$\mu$ by adding certain \r-ribbons in a left-to-right fashion. A static
formulation of the relation~$s\leqhr t$ for bit sequences~$s,t$ could be given
by requiring the existence of a set of intermediate bit sequences, and
specifying the relations that force them to correspond to the mentioned type
of transformation. There is however a more economical way to proceed, based on
the observation that adding \r-ribbons induces only a minimal change on the
bit sequence: it will suffice to require the existence of a single
``intermediate'' bit sequence~$w$, with a simple condition for the way it
differs from $s$~and~$t$. At the same time $w$ will allow reading off the
height of the horizontal \r-ribbon strip~$t/s$. During the addition of the
ribbons of a horizontal \r-ribbon strip, a single bit~$s_i$ can get changed at
most twice: once for a ribbon with its head on diagonal~$i$, and once for a
ribbon with its tail on diagonal~$i+1$, and therefore with its head on
diagonal~$i+r$; moreover, should both changes occur, then the former will
precede the latter due to the left-to-right requirement. The bit $w_i$ will
describe the state of~$s_i$ after a possible modification of the first type,
while $t_i$ gives its final state after a possible modification of the second
type. The principal condition on $s$, $w$, and~$t$ relates the two changes
made by adding one ribbon: $t_{i-r}$ differs from~$w_{i-r}$ if and only if
$w_i$ differs from~$s_i$, in which case one must have
$(w_{i-r},t_{i-r},s_i,w_i)=(1,0,0,1)$. The bit sequence to which this ribbon
is added consists of the bits~$t_j$ for $j<i-r$ followed by the bits $w_j$ for
$i-r\leq j<i$ followed by $s_j$ for~$j\geq i$; in particular, the height of
the ribbon is given by $\sum_{j=1}^{r-1}w_{i-j}$. One additional condition is
needed to characterise~$s\leqhr t$, namely that
$(w_{i-r},t_{i-r},s_i,w_i)=(1,0,0,1)$ occurs only for finitely many values
of~$i$ (since horizontal \r-ribbon strips must be finite); this means that
$s_i=w_i=t_i$ for $i\ll0$ and for $i\gg0$. We state formally:

\proclaim Definition.
Let $s,t$ be bit sequences, then $t/s$ is a horizontal \r-ribbon strip
(written $s\leqhr t$) if there exists a bit sequence~$w$ such that
$(w_{i-r},t_{i-r},s_i,w_i)\in
\setof(a,a,b,b):a,b\in\{0,1\}\endset\union\{(1,0,0,1)\}$ for
all $i\in\Z$, while $(w_{i-r},t_{i-r},s_i,w_i)=(1,0,0,1)$ occurs only for
finitely many~$i$. We say that $w$ is the witness for~$s\leqhr t$, that
$t/s$ has an \r-ribbon at position~$i$ if $i\in I=\setof
i\in\Z:(w_{i-r},t_{i-r},s_i,w_i)=(1,0,0,1)\endset$, and that this ribbon has
height~$\sum_{j=1}^{r-1}w_{i-j}$. We put $|t/s|_r=\Card I$ and $\spin(t/s)
={1\over2}\sum_{i\in{I}}\sum_{j=1}^{r-1}w_{i-j}$.

The second part of this definition depends on the uniqueness of a witness
for~$s\leqhr t$, which follows from the finiteness of~$I$ by an easy induction
starting from $\min(I)=\min(\setof i\in\Z:s_i\neq t_i\endset)+r$, or from
$\max(I)=\max(\setof i\in\Z:s_i\neq t_i\endset)$. We shall also need following
more detailed statement.

\proclaim Lemma. \ribbonmerge
Let $s,t,u\:\Z\to\{0,1\}$ be bit sequences such that $s\leqhr t\leqhr u$ and
$s\leqhr u$ hold. Then the set of positions at which $u/s$ has an \r-ribbon is
the disjoint union of the sets $I,J$ of positions at which $t/s$ and $u/t$,
respectively, have \r-ribbons. Moreover the sums of the heights of the ribbons
of~$u/s$ with positions in each of the sets~$I,J$ exceed the sums of the
heights of the ribbons of~$t/s$ and~$u/s$, respectively, in both cases by the
same number $c\in\N$.

\proof
From a relation $p\leqhr q$ it follows that for the first difference $p_i\neq
q_i$ (i.e., with minimal~$i$), if any, one has $p_i=1$ and~$q_i=0$. In the
situation of the lemma this implies that $s=u$ can only happen if also $s=t$
(in which case the statement is trivial), and that otherwise the minimal
position~$i$ at which $u/s$ has an \r-ribbon is the minimal element
of~$I\union J$, and $i$  does not occur in~$I\thru J$. We shall proceed by
induction on $|u/s|_r$. Suppose first that $i\in I$, so that in fact
$i=\min(I)$. Then the \r-ribbons at position~$i$ of $u/s$ and~$t/s$ have the
same height, namely $\sum_{j=1}^{r-1}s_{i-j}$. With $s'$ the bit sequence
obtained by adding this ribbon to~$s$ (i.e., $s'_j=s_j$ except for
$s'_{i-r}=0$ and $s'_i=1$) we shall apply the induction hypothesis
to~$s',t,u$. The bit sequences that are witnesses for $s'\leqhr t$
and~$s'\leqhr u$ differ only from the witnesses for $s\leqhr t$ and~$s\leqhr
u$ at index~$i-r$ where they have a bit~$0$ instead of~$1$ (their bits at
index~$i$ stay~$1$). Then at each position in~$I-\{i\}$ the \r-ribbons of
$t/s'$ and~$u/s'$ have the same heights as the respective ribbons at the same
position in $t/s$ and~$u/s$, and the lemma follows easily for this case.

The case $i\in J$ is sightly more difficult. As before we define~$s'$ by
adding the ribbon at position~$i$ to~$s$, but since $s'\not\!\leqhr t$ we must
add a ribbon to~$t$ as well: we define a bit sequence~$t'$ that differs
from~$t$ only by $t'_{i-r}=0$ and $t'_i=1$. One has $t'\leqhr u$ and $s'\leqhr
u$ for the same reason as one had $s'\leqhr t$ and $s'\leqhr u$ in the
previous case; in particular each \r-ribbon of $u/t'$ or~$u/s'$ has the same
height as the ribbon of $u/t$ or~$u/s$, respectively, at the same position.
The witness~$w'$ for $s'\leqhr t'$ differs from the witness~$w$ for $s\leqhr
t$ in that $(w'_{i-r},w'_i)=(0,1)$ whereas $(w_{i-r},w_i)=(1,0)$ (note that
also $(s_{i-r},s_i)=(t_{i-r},t_i)=(1,0)$ and
$(s'_{i-r},s'_i)=(t'_{i-r},t'_i)=(0,1)$). An \r-ribbon of $t'/s'$ at
position~$i+j\in I$ (with necessarily $j>0$) may have a height greater
by~$1$ than the \r-ribbon of~$t/s$ at the same position, namely if the sum
giving that height involves the term~$w'_i=1$, in other words if $j<r$. So
with $c_0=\Card\setof j\in\set{r}:i+j\in I\endset$ one has
$\spin(t/s)=\spin(t'/s')-{c_0\over2}$. At the same time,
$\hgt(s'/s)=\sum_{j=1}^{r-1}s_{i-j}$ and
$\hgt(t'/t)=\sum_{j=1}^{r-1}t_{i-j}$ are the heights of the initial ribbons
of $u/s$ and of~$u/t$, and in these sums the differences between corresponding
terms occur when $i+j\in I$, in which case $(s_{i-r+j},t_{i-r+j})=(1,0)$;
therefore $\hgt(s'/s)=\hgt(t'/t)+c_0$. We conclude that if the induction
hypothesis holds for $s',t',u$ with a value~$c_1$ for~$c$, then the lemma
holds for $s,t,u$ with $c=c_0+c_1$.
\QED  

\proclaim Definition. \bscdef
Let $l,m,n,k\:\Z\to\{0,1\}$ be bit sequences, and $a_\ll,a_\gg\in\N^r$. Then
we call the data $({l~m\choose n~k},a_\ll,a_\gg)$ a basic square configuration
if the skew shapes $m/l$, $n/l$ $k/m$, $k/n$ and $k/l$ are all horizontal
\r-ribbon strips, and with $w$ the witness for $l\leqhr k$ there exists a map
$a\:\Z\to\N^r$ such that for all $i\in\Z$ the following statements hold, where
$h=\sum_{j=1}^{r-1}w_{i-j}$:
\statitem If $k/l$ has no \r-ribbon at position~$i$, then $a(i)=a(i+1)$.
\statitem If $w_{i-r}\neq w_i$, then $a(i)_h=0$.
\statitem If $k/l$ has an \r-ribbon at position~$i$ (whose height is 
          then~$h$), put $d=1$ if both $m/l$ and $n/l$ have \r-ribbons at
          position~$i$, put $d=-1$ if both $k/m$ and $k/n$ have \r-ribbons at
          position~$i$, and put $d=0$ otherwise. Then $a(i+1)=a(i)+de_h$,
          i.e., $a(i+1)_h=a(i)_h+d$, and $a(i+1)_j=a(i)_j$ for all $j\neq h$.
\statitem If $k/l$ has no \r-ribbons at positions~$j<i$, then $a(i)=a_\ll$.
\statitem If $k/l$ has no \r-ribbons at positions~$j\geq i$, then
          $a(i)=a_\gg$.

We shall show that these conditions encode the rules of the mutually inverse
algorithms described informally above, when $l,m,n,k$ are the respective edge
sequences of $\\,\mu,\nu,\kappa$, when $a_\ll$ is the matrix entry~$a$, and
$a_\gg=0$. The crucial condition~\statitemnr2 will also help us to prove the
fact that the algorithms are indeed each others inverses (a basic square
configuration can be completed in two directions from two different sets of
partial information). It follows from the first condition that~$a(i)$ becomes
stationary for~$i$ at either side of a finite inteval determined by~$k/l$; the
two final conditions merely state that $a_\ll$ and~$a_\gg$ give the limiting
values. When dealing with edge sequences of partitions, one will always
have~$a_\gg=0\in\N^r$. It is nevertheless useful to allow both limiting values
to be arbitrary: firstly because this gives greater generality while
exhibiting the link with the asymptotic behaviour of the bit sequences, and
secondly because it allows us to recognise a formal symmetry between the two
directions of computation (in addition to the more obvious symmetry with
respect to $m$~and~$n$).

\proclaim Proposition. \symprop
Let $\rev{s}$ denote the reverse of a bit sequence~$s$, defined by
$\rev{s}_i=s_{-1-i}$. Then $({l~m\choose n~k},a_\ll,a_\gg)$ is a basic square
configuration if and only if
$({\rev{k}~\rev{n}\choose\rev{m}~\rev{l}},a_\gg,a_\ll)$ is one.

\proof It is easy to see that $s\leqhr u$ if and only if
$\rev{u}\leqhr\rev{s}$, since if $w$ is the witness for the former relation,
then $\rev{w}$ is the witness for the latter. Let $a\:\Z\to\N^r$ be the map
used to establish that $({l~m\choose n~k},a_\ll,a_\gg)$ is a basic square
configuration, then a direct verification then shows that $a'\:i\mapsto
a(r-i)$ establishes the second basic square configuration. In this
verification one uses the fact that when \bscdef\statitemnr2 applies, then so
does \bscdef\statitemnr1, so $a(i)=a(i+1)$. The reverse implication follows by
symmetry.
\QED

In the following theorem, our main theorem, we need information about the
asymptotic behaviour of the bit sequences (since they differ only at finitely
many places, it does not matter which one of them is considered). We extract
this information by the ``$\,\liminf\,$'' operation. The use of this operation
in a combinatorial context is somewhat surprising, but since it is applied to
expressions taking values in a finite set, it merely returns the smallest
value that is assumed infinitely often by the expression. That expression is
the sum of~$r$ consecutive bits, so for the case of edge sequences one has
$h_\ll=r$ below, which will leave $a_\ll$ unrestricted, and $h_\gg=0$ which
will force~$a_\gg=0\in\N^r$. The previous proposition allows us to state the
theorem in just one direction; the ``insertion'' direction slightly
facilitates the formulation.

\proclaim Theorem. \mainthm
Let $r\in\Np$ and let $l,m,n\:\Z\to\{0,1\}$ be bit sequences with
$m\geqhr{l}\leqhr{n}$; let $a_\ll\in\N^r$ be such that with
$h_\ll=\liminf_{i\to-\infty}\sum_{j\in\set{r}}l_{i-j}$ one has $(a_\ll)_h=0$
for $h_\ll\leq h<r$. Then there exist unique $k\:\Z\to\{0,1\}$ and
$a_\gg\in\N^r$ such that $({l~m\choose n~k},a_\ll,a_\gg)$ is a basic square
configuration; moreover, its values~$a,w$ satisfy for all~$i\in\Z$ the
condition $C_i$: $a(i)_h=0$ for $h_i\leq h<r$ where
$h_i=\sum_{j\in\set{r}}w_{i-1-j}$. Finally, one has $(a_\gg)_h=0$ for
$h_\gg\leq h< r$ where $h_\gg=\liminf_{i\to+\infty}\sum_{j\in\set{r}}l_{i-j}$.

\proof
If $m=l=n$ while $(a_\gg)_h=0$ for $h_{\min}\leq h<r$ where
$h_{\min}=\min_{i\in\Z}(\sum_{j\in\set{r}}l_{i-j})$, then one can take $k=l$
and for $a\:\Z\to\N^r$ the constant function with value $a_\ll$ which then is
also the value of~$a_\gg$: in definition~\bscdef, conditions~\statitemnr1,
\statitemnr4, and~\statitemnr5 hold for all~$i$, and one has $w=l$ so that in
the cases where \statitemnr2 applies, its value of~$h$
equals~$h_{\min}\leq\min(h_i,h_{i+1})<r$ whence $a(i)_h=(a_\ll)_h=0$.
That solution for this case is also unique, because if $k/l$ should have any
\r-ribbon, then applying condition~\statitemnr3 for its ribbon at minimal
position~$i$, the conclusion $a(i+1)=a(i)-e_h$ is incompatible with the
fact that for $h=\sum_{j=1}^{r-1}w_{i-j}=\sum_{j=0}^{r-1}l_{i-j}\geq h_{\min}$
one has $a(i)_h=(a_\ll)_h=0$.

In the remaining cases we can set $i_0$ to the smallest value such that either
at least one of $\{m/l,n/l\}$ has an \r-ribbon at position~$i_0$, or
$h=\sum_{j\in\set{r}}l_{i_0-j}$ satisfies $h<r$ and $(a_\gg)_h>0$. By an
argument similar to the one above $k/l$ cannot have any \r-ribbons at
positions $i<i_0$, so that any solution should have $w_i=l_i$, \
$k_{i-r}=l_{i-r}$, and $a(i+1)=a_\ll$ for $i<i_0$; in particular $C_i$ will be
satisfied for~$i\leq i_0$. We now show for any~$i\in\Z$ that, when the values
$a(i)$ and $w_j$, $k_{j-r}$ for $j<i$ are known and $C_i$ holds, one can
uniquely determine $w_i$, $k_{i-r}$ and $a(i+1)$ and deduce $C_{i+1}$; this
suffices to determine a unique solution. Suppose first that
$(w_{i-r},l_i)\neq(1,0)$, so that $k/l$ cannot have an \r-ribbon at
position~$i$, which forces $(k_{i-r},w_i)=(w_{i-r},l_i)$, and $a(i+1)=a(i)$ by
condition~\bscdef\statitemnr1. In this case condition~$C_i$ implies~$C_{i+1}$,
since $h_{i+1}=h_i-w_{i-r}+w_i\geq h_i$. In the more specific case
$(w_{i-r},l_i)=(0,1)$ it does so with a component of~$a(i+1)$ to spare; more
significantly it ensures condition~\bscdef\statitemnr2.

Next assume that $(w_{i-r},l_i)=(1,0)$, and put $h=\sum_{j=1}^{r-1}w_{i-j}$.
Then it will be possible to put $(k_{i-r},w_i)=(0,1)$, so that $k/l$ has an
\r-ribbon at position~$i$, except if \bscdef\statitemnr3 would then cause
$a(i+1)_h$ to become negative, i.e., if neither $m/l$ nor $n/l$ has an
\r-ribbon at position~$i$, and $a(i)_h=0$. Another potential possibility is to
put $(k_{i-r},w_i)=(w_{i-r},l_i)=(1,0)$, so that $k/l$ has no~\r-ribbon at
position~$i$. For this it is necessary that neither $m/l$ nor $n/l$ have such
a ribbon, and \bscdef\statitemnr2 requires moreover that in this case
$a(i)_h=0$; one arrives at the same condition that was excluded for the first
option. As a consequence there remains a unique possibility for the value of
$(k_{i-r},w_i)$ in all cases. As for condition~$C_{i+1}$, in the cases where
$k/l$ has an \r-ribbon at position~$i$ it is equivalent to~$C_i$, since
$w_{i-r}=1=w_i$ implies $h_i=h_{i+1}$ while $h=h_i-1<h_{i+1}$ is the only
index at which $a(i+1)$ may differ from~$a(i)$. In the case where $k/l$ has no
\r-ribbon at position~$i$, one has $w_{i-r}=1$ and $w_i=0$, so
$h_{i+1}=h_i-1=h$, requiring an additional component of $a(i+1)=a(i)$ to
vanish, namely $a(i+1)_h$. But this is equal to the value~$a(i)_h$ that was
required to be zero in order to have this case in the first place, so
$C_{i+1}$ still holds.

It remains to verify some statements after all of $a$ and~$l$ have been
determined. First, let $S\subset\Z$ be the set of positions~$i$ where it was
decided that $k/l$ has an \r-ribbon, then for every $i\in S$ there either
exists an \r-ribbon of $m/l$ or of~$n/l$ (which happens only finitely often),
or one has $|a(i+1)|<|a(i)|$; this forces $S$ to be finite, whence indeed
$l\leqhr k$. The positions~$i$ for which $a(i+1)\neq a(i)$ are a subset
of~$S$, whence $a(i)$ eventually becomes stationary at some value~$a_\gg$ as
$i\to+\infty$; moreover one sees that \bscdef\statitemnr4~and~\statitemnr5
hold. When $i\notin S$ one has $w_i=l_i$, so under `$\liminf$' in the
expression defining~$h_\gg$ one may replace $l$ by~$w$, and the final
statement of the theorem is a consequence of the validity of~$C_i$ for
all~$i\in\Z$.
\QED

Let us recapitulate somewhat less formally the roles played by condition~(2)
of definition~\bscdef\ and condition~$C_i$ of the theorem. The primary
function of~(2) is to remove the ambiguity when one could either add a ribbon
at position~$i$ to~$k/l$, leaving $h_i$ unchanged but decreasing
$a(i)_{h_i-1}$, or not add such a ribbon, decreasing~$h_i$ while leaving
$a(i)$ unchanged. By forbidding the decrease of~$h_i$ to some level~$h$ until
$a(i)_h$ is~$0$, the second option is restricted to the case when the first
option is not available. Due to the symmetry expressed by
proposition~\symprop, a similar ambiguity is removed in the opposite
direction: in the case where $k/l$ has a ribbon at position~$i$ but neither
$k/m$ nor~$k/n$ do, the data~$(a(i+1),h_{i+1})$ can be distinguished from any
that could arise for a case where $h_{i+1}$ was increased from its previous
level $h_i=h$, since the latter requires $a(i)_h=0$. But that requirement
could also block progress, because the case where $h_i$ is increased is one
where there is no choice: should it occur while $a(i)_{h_i}\neq0$, then there
is no way to satisfy definition~\bscdef. This is where condition~$C_i$ comes
in, as is ensures that whenever the current level $h=h_i$ is not~$r$ one has
$a(i)_h=0$. It must require so as well for all higher levels to anticipate
raising of the level, while as we saw condition~(2) takes care of
preserving~$C_i$ when lowering of the level.

By proposition~\symprop\ the final data $a_\gg,k$ of theorem~\mainthm\ can
also be used to uniquely recover $a_\ll$ and~$l$, so that the construction
establishes the following bijective correspondence.

\proclaim Corollary. \maincor
Let $r\in\Np$, let $m,n$ be a pair of bit sequences that differ only at
finitely many positions, and put
$h_\ll=\liminf_{i\to-\infty}\sum_{j\in\set{r}}m_{i-j}$ and
$h_\gg=\liminf_{i\to-\infty}\sum_{j\in\set{r}}m_{i-j}$. Then there is a
bijection between on one side the pairs $(a_\ll,l)$ of a bit sequence~$l$ with
$m\geqhr{l}\leqhr{n}$ and $a_\ll\in\N^r$ with $(a_\ll)_h=0$ for all
$h_\ll\leq{h}<r$, and on the other side the pairs $(a_\gg,k)$ of a bit
sequence~$k$ with $m\leqhr{k}\geqhr{n}$ and $a_\gg\in\N^r$ with $(a_\gg)_h=0$
for all $h_\gg\leq{h}<r$; this bijection is determined by the fact that
$({l~m\choose n~k},a_\ll,a_\gg)$ is a basic square configuration.
\QED

Taking for $m,n$ the edge sequences of partitions $\mu,\nu$ one has $h_\ll=r$
and $h_\gg=0$, the latter of which forces $a_\gg=0\in\N^r$. Then in the
corollary the bit sequences $l,k$ will also be edge sequences of partitions,
say of $\\,\kappa$ respectively, and on has $\mu\geqhr\\\leqhr\nu$ and
$\mu\leqhr\kappa\geqhr\nu$. Thus we arrive at our conclusion:

\proclaim Theorem. \shapedatumthm
Let $r\in\Np$. There exists a shape datum~$(b_{\mu,\nu})_{\mu,\nu\in\Y}$ for
$(\Y,\leqhr,\N^r)$ determined by the requirement that
$(a,\\)=b_{\mu,\nu}(\kappa)$ holds if and only if
$({\edge\\~\edge\mu\choose\edge\nu~\edge\kappa},a,0)$ is a basic square
configuration. Moreover, in that case one has
$\spin(\kappa/\mu)+\spin(\kappa/\nu)=\spin(\mu/\\)+\spin(\nu/\\)+\n(a)$.

\proof
The construction of the bijections~$b_{\mu,\nu}$ is given by
corollary~\maincor; the only things left to check are the condition
$\rwt\kappa{}-\rwt\mu{}-\rwt\nu{}+\rwt\\{}=|a|$ required by the notion of
shape datum, and the one in final statement of the theorem. The former follows
by observing that for any basic square configuration $({l~m\choose
n~k},a_\ll,a_\gg)$ one has $|a_\ll|=|a_\gg|+\sum_{i\in\Z}\(|a(i)|-|a(i+1)|\)$
(the sum is effectively finite), and that $|a(i)|-|a(i+1)|$ is equal (by
\bscdef(1,3) and the first part of lemma~\ribbonmerge) to the contribution of
ribbons at position~$i$ to the value
$\rwt\kappa\mu-\rwt\nu\\=\rwt\kappa\nu-\rwt\mu\\$. For final statement of the
theorem one uses the second part of lemma~\ribbonmerge. For a horizontal
\r-ribbon strip~$t/s$ let us denote by $\heads ts$ the set of positions at
which $t/s$ has a ribbon, and for $i\in\heads ts$ by $\spin_i(t/s)$ the
contribution of that ribbon to $\spin(t/s)$ (i.e.,
${1\over2}\sum_{j=1}^{r-1}w_{i-j}$ where $w$ is the witness for~$s\leqhr{t}$).
Then the lemma says that
$$\eqalignno{
  \sum_{i\in\heads ml}(\spin_i(k/l)-\spin_i(m/l))
&=\sum_{i\in\heads km}(\spin_i(k/l)-\spin_i(k/m))
\cr\noalign{\hbox{and also}}
  \sum_{i\in\heads nl}(\spin_i(k/l)-\spin_i(n/l))
&=\sum_{i\in\heads kn}(\spin_i(k/l)-\spin_i(k/n))
\cr}
$$ 
for ${l~m\choose{n}~k}={\edge\\~\edge\mu\choose\edge\nu~\edge\kappa}$. Then
the expression
$\spin(\kappa/\mu)+\spin(\kappa/\nu)-\spin(\mu/\\)-\spin(\nu/\\)$ that we
should prove equal to~$\n(a)=\n(a_\ll)$ can be rewritten as
$$
  \def\+#1#2{\sum_{i\in I(#1/#2)}\spin_i(k/l)}
  \+km+\+kn-\+ml-\+nl.
$$
It is readily checked that any ribbon considered in \bscdef(3) contributes
$-dh$ to this expression, which can be written as $-\n(de_h)$; according to
\bscdef(3) this is equal to $\n(a(i))-\n(a(i+1))$. Like above one has
$\n(a_\ll)=\n(a_\gg)+\sum_{i\in\Z}\(\n(a(i))-\n(a(i+1))\)$, while
$\n(a_\gg)=0$; this gives the required identity.
\QED

Now by constructing Knuth-growths for this shape datum as described in
\Sec\FominKnuthsec, we obtain the spin preserving Knuth correspondence that we
set out to find; it is symmetric since $b_{\mu,\nu}=b_{\nu,\mu}$ for all
$\mu,\nu\in\Y$. We summarise its main characteristics.

\proclaim Corollary. \spinKnuthcor
Let $r\in\Np$, let $\gamma$ be an \r-core, and $\RH(\gamma)$ the
associated \r-rim hook lattice (the connected component of~$(\Y,\leqr)$
containing~$\gamma$). Then for any $m,n\in\N$ there exists a bijection between
on one side $m\times n$ matrices~$A$ with entries in~$\N^r$, and on the other
side pairs $(P,Q)$ of semistandard \r-ribbon tableaux of equal
shape~$\\/\gamma$ for some $\\\in\RH(\gamma)$, such that one has
$\wt(P)=\(\sum_{i\in\set{m}}|A_{i,j}|\){}_{j\in\set{n}}$ and
$\wt(Q)=\(\sum_{j\in\set{n}}|A_{i,j}|\){}_{i\in\set{m}}$, while also
$\spin(P)+\spin(Q)=\sum_{i,j\in\set{m}\times\set{n}}\n(A_{i,j})$. Moreover the
correspondence is symmetric: if $A$ is mapped to $(P,Q)$, then the transpose
matrix $A\tr$ maps to~$(Q,P)$.
\QED

We remark that theorem~\shapedatumthm\ proves equation~(\qbscid). For
completeness we shall write down the equation the corresponds to the full
Knuth correspondence; it will be the obvious $q^{1\over2}$-analogue
of~(\CauchyKnuthid). So define $G^{(r)}_{\\/\gamma}(q^{1\over2},X)
=\sum_Pq^{\spin(P)}X^{\wt(P)}\in\Z[q^{1\over2}][[X]]$, where the sum is over
all semistandard \r-ribbon tableaux~$P$ of shape~$\\/\gamma$, identifying
tableaux that differ only by stationary steps at the end of the path (staying
at the shape~$\\$). Then one obtains from interpreting
corollary~\spinKnuthcor, or directly by multiplying copies of~(\qbscid):
$$
 \prod_{i,j\in\N}\prod_{k\in\set{r}}{1\over1-q^kX_iY_j}
 =\sum_{\\\in\RH(\gamma)}G^{(r)}_{\\/\gamma}(q^{1\over2},X)
                         G^{(r)}_{\\/\gamma}(q^{1\over2},Y)
.\label(\spinCauchyKnuthid)
$$

Returning to the questions considered in \Sec\enumsec, we observe that
corollary~\maincor\ applied to cases with $n=m$ for which $h_\ll=0=h_\gg$
provides a bijective proof of claim~\polderclaim. For cases with $n=m$ for
which $h_\ll=r$ and~$h_\gg=0$ it similarly proves claim~\partitionclaim, while
for cases with $n=m$ for which $h_\ll=r=h_\gg$ it proves the equation obtained
from that of claim~\alpineclaim\ by multiplying both sides by
$\prod_{k\in\set{r}}{1\over1-XY^k}$.

\subsecno=-1 
\newsection Asymmetric correspondences.
\labelsec\asymSec

In this section we shall define a generalisation to semistandard \r-ribbon
tableau of the asymmetric correspondence defined by Knuth, the one that proves
the identity $\prod_{i,j}(1+X_iY_j)=\sum_\\s_\\(X)s_{\\\tr}(Y)$; it will be
spin preserving in the same sense as our symmetric Knuth correspondence. Much
of the considerations will be similar to those of the symmetric case, but at
the heart of the construction there is a fundamental difference, more
substantial than the difference between Knuth's symmetric and asymmetric
constructions. In fact our construction will be a bit simpler for the
asymmetric case, since it is based on a direct study of all possible values
for~$\\$ and~$\kappa$ for given $\mu,\nu$. 

First we need to give the context in which the construction applies, which
involves less specialised instances of the constructions of \ref{Fomin Schur}
than we have been considering so far. Specifically, we have been supposing
that the horizontal neighbours in the grid, like $\\,\mu$ are related by the
same relation~`$\leqhr$' as the vertical neighbours like $\\,\nu$. There is
however nothing in the construction a Knuth correspondence that requires this,
and the construction applies equally well when two different relations are
used horizontally and vertically. Specifically, we shall consider cases where
the two relations are each others transposes, for instance when requiring
$\\\leqh\mu$ and $\nu\leqh\kappa$, we shall also require $\\\tr\leqh\nu\tr$
and $\mu\tr\leqh\kappa\tr$. Note that the relations `$\prec$' and `$\precr$'
coincide with their transposes, which is why we shall not consider asymmetric
Schensted correspondences (although for appropriate ``dual graded graphs''
they too can be constructed; see \ref{Fomin duality} for examples). To
facilitate notation we shall write $\\\leqv\nu$ instead of $\\\tr\leqh\nu\tr$,
in which case we say $\nu/\\$ is a vertical strip, and similarly $\\\leqvr\nu$
means $\\\tr\leqhr\nu\tr$, and $\nu/\\$ is then said to be a vertical
\r-ribbon strip. Recall that for `$\leqh$' and~`$\leqhr$' we have chosen to
form the symbols for the opposite relations by rotation rather than by
reflection (which convention has the advantage that it can be extended without
ambiguity to relations written vertically or diagonally in diagrams, even if
this does not occur in our paper); by the same token $\nu\geqv\\$ will mean
that $\nu/\\$ is a vertical strip. Standardisation of such strips is defined
so as to commute with transposition: squares are added from top to bottom.

We shall use the relations `$\leqh$' or~`$\leqhr$' between shapes associated
to horizontally adjacent points of our grid, and `$\leqh$' or~`$\leqhr$' for
vertically adjacent grid points; this is easy to remember, but there is of
course no connection between relative positions among squares of individual
shapes and relative positions of those shapes as placed on the grid: growths
with the opposite convention can equally well be defined. The notion of a
shape datum is adapted to the asymmetric context, in the obvious way.

\proclaim Definition. \asddef
Let $\Part$ be a graded set equipped with two relations~`$\leqh$', `$\leqv$',
and $S$ a graded set. An (asymmetric) \newterm{shape datum} for
$(\Part,\leqh,\leqv,S)$ consists of a family
$(b_{\mu,\nu})_{\mu,\nu\in\Part}$ of bijections
$$
\predisplaypenalty=10000 \postdisplaypenalty=10000
  b_{\mu,\nu}\:
 \setof\kappa\in\Part:\mu\leqv\kappa\geqh\nu\endset
\to
 \setof(a,\\)\in S\times\Part:\mu\geqh\\\leqv\nu\endset,
\label(\asymsdeq)
$$
such that $|\kappa|-|\nu|-|\mu|+|\\|=|a|$ holds whenever
$(a,\\)=b_{\mu,\nu}(\kappa)$.

The basic asymmetric shape datum is for $(\Y,\leqh,\leqv,\{0,1\})$, and is
derived from Knuth's asymmetric construction. In fact it comes in two
different flavours, one for each of the natural $1$-correspondences
for~$(\Y,\prec)$. As before they can be defined using standardisations of
$\kappa/\mu$ and $\kappa/\nu$ and computing a partial Schensted-growth across
a rectangle, by taking for~$a$ the sum of the matrix entries found, and
for~$\\$ the shape at the top left corner. There will be at most one matrix
entry~$1$, which occurs at the top right corner (where $\mu$ is placed) in the
case of row insertion, and in the bottom left corner (with~$\nu$) in the case
of column insertion. One cannot really say which of these flavours matches
Knuth's original construction, since that construction performs row insertion
into a row-strict (transpose semistandard) tableau, with the recording tableau
being column-strict (semistandard); with our conventions however the
$P$-symbol will be column-strict and the $Q$-symbol row-strict. Two different
symmetries can be used to map Knuth's construction to our conventions: one may
replace $(P,Q)$ either by~$(P\tr,Q\tr)$ or by~$(Q,P)$; this respectively
results in the column-insertion and the row-insertion flavour, and shows that
the two are essentially equivalent, unlike what we saw for the RSK and Burge
correspondences.

Before considering the ribbon case, we shall reformulate in a more direct way
these asymmetric shape data for $(\Y,\leqh,\leqv,\{0,1\})$. By studying the
sets at both sides of the bijection in~(\asymsdeq), we can give a description
involving hardly any algorithm at all. For given $\mu,\nu\in\Y$, it is clear
that any diagram $\\$ satisfying $\mu\geqh\\\leqv\nu$ will be contained in
$\\'=\mu\thru\nu$, and any $\kappa$ satisfying $\mu\leqv\kappa\geqh\nu$ will
contain $\kappa'=\mu\union\nu$. The conditions $\\'\leqh\mu$ and
$\nu\leqh\kappa'$ are equivalent since the diagrams of $\mu/\\'$ and
$\kappa'/\nu$ are identical, and similarly $\\'\leqv\nu$ is equivalent to
$\mu\leqv\kappa'$; these conditions must be verified for the sets at either
sides of~(\asymsdeq) to be non-empty. Moreover any squares of $\\'/\\$ must be
in distinct rows and in distinct columns, nor can they share a column with any
square of the horizontal strip~$\mu/\\'$, or a row with any square of the
vertical strip~$\nu/\\'$; similar statements hold for $\kappa/\kappa'$. It
follows that there is a set~$S$ of squares that might be removed from~$\\'$,
and removing any subset of~$S$ gives a valid~$\\$, and similarly there is a
set of squares~$T$ that might be added to~$\kappa'$ to get~$\kappa$, and
adding any subset of~$T$ is valid. The squares of $S$ and of~$T$ are perfectly
interleaved in bottom-left to top-right order, with an element of~$T$ at
either end; then after choosing either to match each square of $S$ with the
next or with the previous square of~$T$ (much like the two $1$-correspondences
for~$(\Y,\prec)$), and matching the remaining square of~$T$ with the value of
$a\in\{0,1\}$, one discerns a shape datum is for $(\Y,\leqh,\leqv,\{0,1\})$.
The fact that the squares of $S$ and~$T$ are interleaved can be shown by
procedures that essentially trace bumping sequences in the Schensted-growth
above, but it will be more useful to express $S$ and~$T$ directly in terms of
$\edge\mu$ and $\edge\nu$.

Let us define a doubly infinite sequence $\Edge\\\:\Z\to\N$ for any $\\\in\Y$
by cumulation of~$\edge\\$, setting $\Edge\\_k=\sum_{i\geq k}\edge\\_i$.
This number can be interpreted graphically as the vertical coordinate of the
point where the boundary of~$\\$ crosses the diagonal~$d_k$. Then $\\\leqh\mu$
can be seen to be equivalent to $\Edge\\_k\leq\Edge\mu_k\leq\Edge\\_{k+1}+1$
for all $k\in\Z$, and $\\\leqv\nu$ similarly to
$\Edge\\_k\leq\Edge\nu_k\leq\Edge\\_{k-1}$. This can be expressed more
elegantly by (termwise) comparison of sequences using the shifted sequences
$\Edgem\\,\Edgep\\$ defined by $\Edgem\\_{i-1}=\Edge\\_i=\Edgep\\_{i+1}$:
$$\eqalignno 
{\\\leqh\mu &\iff             \Edgep\mu-1\leq\Edge\\\leq\Edge\mu\leq\Edgem\\+1
            &\nn\cr
 \\\leqv\nu &\iff \phantom{1+{}}\Edgem\nu\leq\Edge\\\leq\Edge\nu\leq\Edgep\\ 
            &\nn\cr
}
$$
where we have added redundant inequalities for symmetry. Now we see that
$\mu\geqh\\\leqv\nu$ is equivalent to
$\Edgep\mu-1\leq\Edge\\\leq\Edge\nu\leq\Edgep\\\leq\Edgep\mu$, which yields a
contradiction unless the sequence $\Delta(\mu,\nu)\defeq\Edgep\mu-\Edge\nu$ has
all its terms in~$\{0,1\}$; the same necessary condition is found for
$\mu\leqv\kappa\geqh\nu$. From the given inequalities it follows moreover that
if $\Delta(\mu,\nu)_k=1$ then $\Edge\\_k=\Edge\nu_k$, while if
$\Delta(\mu,\nu)_{k+1}=0$ then $\Edge\\_k=\Edge\mu_k$, so unless
$\Delta(\mu,\nu)_k=0$ and $\Delta(\mu,\nu)_{k+1}=1$, the term $\Edge\\_k$ is
completely determined by~$\mu$ and~$\nu$. In that remaining case one easily
shows that $\Edge\\_{k-1}=\Edge\mu_{k-1}=\Edge\nu_k$ is exactly one larger
than $\Edge\\_{k+1}=\Edge\nu_{k+1}=\Edge\mu_k-1$, so that $\Edge\\_k$ may
indeed assume either of these values; in other words there is a square of~$S$
on the diagonal~$d_k$. For~$\kappa$, a similar reasoning shows that
$\Edge\kappa_k$ is determined by~$\mu$ and~$\nu$ unless $\Delta(\mu,\nu)_k=1$
and $\Delta(\mu,\nu)_{k+1}=0$, in which case it can assume either of the
values $\Edge\kappa_{k-1}$ and $\Edge\kappa_{k+1}$, so that $T$ has a square
on the diagonal~$d_k$. We conclude that the squares of~$S$ correspond to
occurences in the bit sequence $\Delta(\mu,\nu)$ of~`$01$', while the squares
of~$T$ correspond to occurences of~`$10$'; the announced perfect interleaving
follows by a familiar argument.

That was a bit technical, but it will allow us advance easily to the case of
\r-ribbons. To find a shape datum for $(\Y,\leqhr,\leqvr,\{0,1\}^r)$, it
suffices to consider pairs $\mu,\nu$ whose \r-quotients
$(\mu^i)_{i\in\set{r}},(\nu^i)_{i\in\set{r}}$ are such that each sequence
$\Delta(\mu^i,\nu^i)$ has all its terms in~$\{0,1\}$, so that shapes $\\$ and
$\kappa$ with $\mu\geqhr\\\leqvr\nu$ and $\mu\leqvr\kappa\geqhr\nu$ do exist.
Then the sequence that will take the place of $\Delta(\mu,\nu)$ above is
$\Delta^r(\mu,\nu)$ defined by 
$$
\Delta^r(\mu,\nu)_k=\Edge\mu_{k-r}-\Edge\nu_k.
\nn
$$ 
To visualise this, imagine two points moving simultaneously, the first along
the boundary of~$\nu$, the second along the boundary of~$\mu$, always keeping
$r$~diagonals behind the first (like the fire department ladder truck with
independently steering rear wheels); then $\Delta^r(\mu,\nu)$ tracks the
vertical distance between the two points as a function of time. From the
defining relation $\edge\\_{i+jr}=\edge{\\^i}_j$ of \r-quotients one deduces
$\Edge\\_k=\sum_{i\in\set{r}}\Edge{\\^i}_{k-i\overwithdelims\lceil\rceil r}$,
from which it follows that $\Delta^r(\mu,\nu)_k
=\sum_{i\in\set{r}}\Delta(\mu^i,\nu^i)_{k-i\overwithdelims\lceil\rceil r}$.
This shows that $0\leq\Delta^r(\mu,\nu)_k\leq r$ for all~$k$, and moreover
that
$$
 \edge\mu_{k-r}-\edge\nu_k
=\Delta^r(\mu,\nu)_k-\Delta^r(\mu,\nu)_{k+1}
=\Delta(\mu^i,\nu^i)_j-\Delta(\mu^i,\nu^i)_{j+1}
,\ttex{where $k=i+jr$.}
\label(\deltaeq)
$$
Let $S'$ be the set of indices~$k$ for which the members of this equation have
the value~$-1$, and $T'$ the set of indices~$k$ for which they have the
value~$+1$. We denote by $\\'=\mu\meet\nu$ and $\kappa'=\mu\join\nu$ the meet
and join of the shapes $\mu,\nu$ in $(\Y,\leqr)$; these operations are not the
intersection and union of the shapes, but the components of their \r-quotients
are obtained in that way. Then the set of all shapes~$\\$ with
$\mu\geqhr\\\leqvr\nu$ is in bijection with the set of subsets of~$S'$, the
shape corresponding to $S\subset S'$ being obtained by successively removing
\r-ribbons from~$\\'$ with their heads on the diagonals~$d_k$ with $k\in S$.
Similarly the set of all shapes~$\kappa$ with $\mu\leqvr\kappa\geqhr\nu$ is in
bijection with the set of subsets of~$T'$, the shape corresponding to
$T\subset T'$ being obtained by successively adding \r-ribbons to~$\kappa'$
with their heads on the diagonals~$d_k$ with $k\in T$. By consideration of the
\r-quotients one can see that $\Card{T'}=\Card{S'}+r$ which is sufficient to
obtain the existence of asymmetric shape data for
$(\Y,\leqhr,\leqvr,\{0,1\}^r)$. But we want more than that of course: we want
a shape datum for which equation~(\spinbalance) holds, so that spin shall be
preserved. It is for that purpose that the sequence $\Delta^r(\mu,\nu)$ will
prove its real utility.

We must first define what exactly we mean by the spin of a vertical \r-ribbon
strip. When defining the standardisation of vertical strips, we have required
commutation with transposition, and we do so for the standardisation of
vertical \r-ribbon strips as well, so that their ribbons are added from top
right to bottom left (formally: by decreasing index of the diagonal containing
the head of the ribbon). We note that skew shapes that happen to be
simultaneously a horizontal \r-ribbon strip and a vertical \r-ribbon strip
have two different standardisations (even without considering the possibility
that they are an ordinary horizontal or vertical strip as well!). The context
will always make clear which form of standardisation is intended, but for the
spin, which is defined in terms of the standardisation, we wish to make the
distinction clear in the notation, so we shall write $\spinv(\nu/\\)$ for the
spin of $\nu/\\$ defined using its standardisation as a vertical \r-ribbon
strip; as for horizontal \r-ribbon strip this spin is half the sum of the
heights of the ribbons in the standardisation. Note that this is not the spin
of the horizontal \r-ribbon strip obtained by transposition (this already
fails for single \r-ribbons), although the two quantities are related.

\proclaim Lemma.
Let $\\,\mu,\nu,\kappa\in\Y$ with $\mu\geqhr\\\leqvr\nu$ and
$\mu\leqvr\kappa\geqhr\nu$. Let $S\ssubset\Z$ be the set of indices~$k$ such
that the standardisations both of $\mu/\\$ and of $\nu/\\$ contain an
\r-ribbon with its head on diagonal $d_k$, and let $T\ssubset\Z$ be the set of
indices~$k$ such that the standardisations both of $\kappa/\mu$ and of
$\kappa/\nu$ contain an \r-ribbon with its head on diagonal $d_k$. Then
$$
 \spinv(\kappa/\mu)+\spin(\kappa/\nu)-\spin(\mu/\\)-\spinv(\nu/\\)
=\sum_{k\in T}\Delta^r(\mu,\nu)_{k+1}-\sum_{k\in S}\Delta^r(\mu,\nu)_k
.\label(\spinunbalance)
$$

\proof
We shall use a slight variation of a Schensted-growth to do the accounting for
us. Let $\heads\mu\\\ssubset\Z$ be such that the heads of the ribbons of the
standardisation of $\mu/\\$ occur on diagonals~$d_k$ for $k\in\heads\mu\\$,
and define $\heads\kappa\nu$, $\heads\nu\\$, $\heads\kappa\mu$ similarly; put
$I=\heads\mu\\\union\heads\kappa\nu=\heads\mu\\\union T
=S\union\heads\kappa\nu$, and
$J=\heads\nu\\\union\heads\kappa\mu=\heads\nu\\\union T
=S\union\heads\kappa\mu$. We consider a grid rectangle whose squares are
numbered horizontally by~$I$ increasing from left to right, and vertically
by~$J$ decreasing from top to bottom. Shapes are assigned to all grid points
according to a simple rule: for any top and left justified sub-rectangle
$I'\times J'$ (so $I'\subset I$ is an order ideal and $J'\subset J$ a dual
order ideal), the shape assigned to its bottom right corner is obtained
from~$\\$ by adding one \r-ribbon with its head on~$d_k$ for each $k\in
I'\union J'$ in any order in which this is possible (the results will be the
same). This implies that standardisations of~$\mu/\\$ and of~$\nu/\\$ are
placed on the grid points along the top and left sides, but repeating the same
shape for every step numbered by an element of~$T$, while standardisations
of~$\kappa/\mu$ and of~$\kappa/\nu$ are placed along the right and bottom
sides, with this time steps numbered by an element of~$S$ being trivial.
Moreover, for any grid square numbered $(i,j)$ with $i\neq j$, one has a
configuration as in a Schensted-growth for a square not involving the
\r-correspondence, i.e., opposite sides correspond to ribbons with their head
on the same diagonal, or to no ribbon at all. The remaining squares
numbered~$(k,k)$ are the more interesting ones: if $k\in S$ then the top left
shape differs by an \r-ribbon from the other three which are equal
(representing annihilation of the ribbon rather than bumping), while if $k\in
T$ the same holds for the bottom right shape (creation of an \r-ribbon, as
happens in Schensted-growths when the matrix entry is non-zero).

For a lattice path from the bottom left corner of the rectangle to the top
right corner, we call half the sum of the heights of ribbons encountered along
the path its spin. Then the left hand side of (\spinunbalance) measures the
amount added to the spin when replacing the path passing through the top left
corner by the one passing though the bottom right corner. Moving the path
across one square at a time, we have seen that the spin is unchanged for
squares not numbered~$(k,k)$. We can conclude the proof by showing that the
change of the spin at any square numbered~$(k,k)$ is equal to the contribution
of~$k$ to the right hand side of~(\spinunbalance). That change is equal to
$\pm\hgt(\kappa'/\\')$, where $\\'$ and $\kappa'$ are the shapes assigned to
the top left respectively bottom right corners of the square, and the sign is
negative if $k\in S$ and positive if $k\in T$. By definition
$\hgt(\kappa'/\\')=\Edge{\\'}_{k-r+1}-\Edge{\\'}_k$, which we can write as
$\Edge{\\'}_{k-r}-\Edge{\\'}_{k+1}-1$ since
$(\edge{\\'}_{k-r},\edge{\\'}_k)=(1,0)$. Now the value of
$\Edge{\\'}_{k-r}$ is unaffected by addition to or removal from~$\\'$ of
\r-ribbons with their head on a diagonal~$d_i$ with $i\geq k$, so
$\Edge\cdot_{k-r}$ is constant on the four corners of our square, and on all
grid points above and to the right of them; in particular
$\Edge{\\'}_{k-r}=\Edge\mu_{k-r}$. It follows by a similar argument that
$\Edge{\\'}_{k+1}=\Edge\nu_{k+1}$, so one obtains
$\hgt(\kappa'/\\')=\Edge\mu_{k-r}-\Edge\nu_{k+1}-1$. Now since $S\subset S'$
and $T\subset T'$ as defined below~(\deltaeq), the value of
$(\edge\mu_{k-r},\edge\nu_k)$ is $(0,1)$ when $k\in S$, and it is $(1,0)$ when
$k\in T$, which allows us to write $\hgt(\kappa'/\\')=\Delta^r(\mu,\nu)_k$ in
the former case and $\hgt(\kappa'/\\')=\Delta^r(\mu,\nu)_{k+1}$ in the latter,
giving the desired contribution to the right hand side of~(\spinunbalance).
\QED

This lemma allows us to define a spin preserving asymmetric shape datum $\bA$
for $(\Y,\leqhr,\leqvr,\{0,1\}^r)$ in the way the reader may have guessed we
wanted to define it. We match elements of~$S'$ corresponding to a rise of
$\Delta^r(\mu,\nu)$ from~$h$ to~$h+1$ with elements of~$T'$ corresponding to a
fall of $\Delta^r(\mu,\nu)$ from~$h+1$ to~$h$, while matching the remaining
element of~$T'$ of that form with the component~$a_h$ of $a\in\{0,1\}^r$; this
is possible since $\lim_{k\to-\infty}\Delta^r(\mu,\nu)_k=r$ and
$\lim_{k\to+\infty}\Delta^r(\mu,\nu)_k=0$. Then any $\kappa$ in the domain
of~$\bA_{\mu,\nu}$, in other words with $\mu\leqvr\kappa\geqhr\nu$, determines
a subset $T\subset T'$ which gives rise to a subset $S\subset S'$ determining
$\\$ with $\mu\geqhr\\\leqvr\nu$, and to a value $a\in\{0,1\}^r$. Then
matching elements of $S$ and~$T$ cancel each others contribution to the right
hand side of~(\spinunbalance), and the contribution of the remaining elements
of~$T$ amounts to $\n(a)$, whence one obtains the asymmetric counterpart of
equation~(\spinbalance). We must still specify the precise matching of
elements of~$S'$ and of~$T'$; of the two equally natural possibilities we
choose the one which is most similar to the symmetric case: we match each rise
of $\Delta^r(\mu,\nu)$ to the \emph{next} descent back to the same level,
leaving the very first descent from~$h+1$ to~$h$ to match~$a_h$.

If one represents $\mu$ and~$\nu$ by their edge sequences and the \r-ribbon
strips $\mu/\\$, $\nu/\\$, $\kappa/\mu$, $\kappa/\nu$ by the sets
$\heads\mu\\,\ldots,\heads\kappa\nu$ giving the diagonals containing the heads
of the ribbons in their standardisations, then there is a simple procedural
description of this asymmetric shape datum. We shall describe it in the
``insertion'' direction, calculating $\kappa=(\bA_{\mu,\nu})\inv(a,\\)$ by
specifying the elements of $\heads\kappa\mu$ and~$\heads\kappa\nu$. The
procedure uses $a$ as an initialised \r-bit variable; a diagonal index~$k$
that traverses a sufficiently large interval of~$\Z$ in increasing order,
while an integer variable $h$ taking values $0\leq h\leq r$ keeps track of the
current value of $\Delta^r(\mu,\nu)_k$. The initial value of~$k$ is taken
sufficiently small so that $k\leq i$ for all
$i\in\heads\mu\\\union\heads\nu\\$ and $(\mu_{i-r},\nu_i)=(1,1)$ for
all~$i<k$; correspondingly $h$ is initialised to~$r$. For each~$k$ the
following cases are distinguished. If $(\mu_{k-r},\nu_k)=(1,0)$ then
necessarily $h>0$ and $k\notin\heads\mu\\\union\heads\nu\\$, and one starts by
setting $h:=h-1$; if now $a_h=1$ then it is established that
$k\in\heads\kappa\mu$ and that $k\in\heads\kappa\nu$, and otherwise that
$k\notin\heads\kappa\mu$ and that $k\notin\heads\kappa\nu$; finally one sets
$a_h:=0$. If $(\mu_{k-r},\nu_k)=(0,1)$ then necessarily $h<r$ and $a_h=0$; in
this case it is established that $k\notin\heads\kappa\mu$ and
$k\notin\heads\kappa\nu$, one sets $a_h:=1$ if
$k\in\heads\mu\\\thru\heads\nu\\$ (otherwise
$k\notin\heads\mu\\\union\heads\nu\\$ and $a_h$ stays~$0$), and finally one
sets $h:=h+1$. In the remaining cases ($\mu_{k-r}=\nu_k$), both $h$ and~$a$ are
unchanged; it will be established that $k\in\heads\kappa\mu$ if and only
if $k\in\heads\nu\\$, and that $k\in\heads\kappa\nu$ if and only if
$k\in\heads\mu\\$ (at most one of these conditions can hold). Traversal may
terminate when $k>i$ for all $i\in\heads\mu\\\union\heads\nu\\$,
and $h=0$.

This procedure has some similarity to the one given for the symmetric shape
datum, in particular if there the level~$h$ is also be maintained in an
incremental fashion. Like we saw for that procedure, $a_i=0$ here holds for
all~$h\leq i<r$ after each step, forcing $a=0$ at the end. It also enables a
step-by-step inverse procedure, but here that is for a rather simple reason:
using the bit $a_h$ to record the presence of $k$ in $\heads\mu\\$ and
in~$\heads\nu\\$ can only be undone if $a_h$ previously had a known state
(cleared). More generally the current procedure is a lot simpler, because $h$
evolves independently of the placement of ribbons.

Let us give an example of this asymmetric shape datum. We take $r=5$, \ 
$\\=(10,10,10,10,4,3,3,1)$, \
$a=(1,0,0,1,1)$, and $\mu,\nu$ are such that the standardisations
of $\mu/\\$ and $\nu/\\$ are
$$
  \mu/\\:~
  {\smallsquares\rtab 5
   \edgeseq15,0;11
   \ribbon 12,0:;1111 \ribbon 10,1:;1110 \ribbon 8,2:;0111
   \ribbon 8,4:;1111 \ribbon 7,5:;1110 \ribbon 6,6:;1010 \ribbon6,7:;0101
   \ribbon 3,10:;1110
   \edgeseq8,0;111111110000000000
   \edgeseq0,12;00000
  }
,\qquad
  \nu/\\:~
  {\smallsquares\rtab 5
   \edgeseq15,0;11111
   \ribbon 9,0:;1000 \ribbon 7,1:;0010 \ribbon 5,3:;0100
   \edgeseq4,7;000111
   \ribbon 0,10:;0000
   \edgeseq8,0;111111110000000000
   \edgeseq0,15;00
  }
.
$$
One can read off ${\edge\mu\choose\edge\nu} =\let~\,{
\phantom{11111}~11011~01100~01010~00101~10111~01000~00
\choose  11111 ~01000~11011~00100~01110~00001~00\phantom{000~00}}$ at
indices $-15\leq i<17$, which we displayed shifted so as to align
$\edge\mu_{i-5}$ with $\edge\nu_i$; from this one can easily determine the
terms of~$\Delta^5(\mu,\nu)$ at indices $-10\leq k\leq17$, namely
$\let~\,~5,4,4,4,3,~2,3,3,2,3,~4,4,3,4,3,~3,3,4,4,5,~4,3,3,2,1,~1,1,0$. It
then follows that $S'=\{-5,-2,-1,2,6,8\}$ and
$T'=\{-10,-7,-6,-3,1,3,9,10,12,13,16\}$, and the injection $S'\to T'$ is given
by $-5\mapsto-3$, \ $-2\mapsto12$, \ $-1\mapsto1$, \ $2\mapsto3$, \
$6\mapsto10$, $8\mapsto9$, and the remaining elements of~$T'$ that match
$\li(a4..0)$ are $-10,-7,-6,13,16$.
One has
$\heads\mu\\=\{-8,-5,-2,0,2,4,5,11\}$ and $\heads\nu\\=\{-5,-2,2,14\}$, whose
intersection gives $S=\{-5,-2,2\}$; from this and the fact that the non-zero
bits $a_i$ are $a_4,a_3,a_0$, one obtains $T=\{-10,-7,-3,3,12,16\}$, and it
can be concluded that $\heads\kappa\mu=\{-10,-7,-3,3,12,14,16\}$ and
$\heads\kappa\nu=\{-10,-8,-7,-3,0,3,4,5,11,12,16\}$. Thus one finds
$\kappa=(\bA_{\mu,\nu})\inv(a,\\)=(17,16,15,12,10,10,10,10,6,6,4,3,3,3,1)$,
with the following 5-ribbon strips
$$
  \kappa/\mu:~
  {\smallsquares\rtab 5
   \ribbon 14,0:;1001 \ribbon 12,1:;1010 \ribbon 9,2:;0001
   \ribbon 7,6:;0001
   \edgeseq6,10;110
   \ribbon 3,11:;1000 \ribbon 1,11:;0000 \ribbon0,12:;0000
   \edgeseq13,0;1111111111111000000000000
  }
,\qquad
  \kappa/\nu:~
  {\smallsquares\rtab 5
   \ribbon 14,0:;1111 \ribbon 13,1:;1111 \ribbon 13,2:;1111
   \ribbon 10,3:;1011 \ribbon 9,5:;1111
   \ribbon 7,6:;1101 \ribbon 7,7:;1011 \ribbon 7,8:;0111
   \ribbon 3,10:;1100 \ribbon 3,11:;1001 \ribbon 2,14:;1010
   \edgeseq10,0;1111111111000000000000000
  }
.
$$
One checks from these results that
$\rwt\kappa\mu-\rwt\nu\\=7-4=|a|=11-8=\rwt\kappa\nu-\rwt\mu\\$ and that
$\spinv(\kappa/\mu)+\spin(\kappa/\nu)-\spin(\mu/\\)-\spinv(\nu/\\)
={7\over2}+{34\over2}-{24\over2}-{3\over2}=7=\n(a)$ as it should be. The
mentioned procedural description of the shape datum can be seen to do
essentially the same computation, but in a more orderly left-to-right fashion.
The reader may check that in fact everything can be done entirely in terms of
edge sequences and in a single pass from left to right, finding $\heads\mu\\$
and $\heads\nu\\$ on the fly by transforming a copy of~$\edge\\$ into
$\edge\mu$ and a copy of~$\edge\nu$ into~$\edge\\$; the information obtained
about $\heads\kappa\nu$ can be used to simultaneously transform another copy
of~$\edge\nu$ into~$\edge\kappa$. Such a description does not appear to given
much new insight into correspondence however. We shall content ourselves here
with summarising the result found.

\proclaim Theorem. 
For every $r>0$ there exists an asymmetric shape datum $\bA$
for~$(\Y,\leqhr,\leqvr,\{0,1\}^r)$ such that whenever
$(a,\\)=\bA_{\mu,\nu}(\kappa)$, then
$\spinv(\kappa/\mu)+\spin(\kappa/\nu)=\spin(\mu/\\)+\spinv(\nu/\\)+\n(a)$.
Hence
$$
  \widetilde D^r_{q^{1/2},Y}\after U^r_{q^{1/2},X}
  =(U^r_{q^{1/2},X}\after \widetilde D_{r,Y})\prod_{i\in\set{r}}(1+q^iXY)
,\nn
$$
where $U^r_{q^{1/2},X}$ is the spin generating series for horizontal \r-ribbon
strips defined in equation~(\spingenseq), while $\widetilde
D^r_{q^{1/2},Y}=\sum_{\mu\leqvr\\}q^{\spinv(\\/\mu)}Y^{\rwt\\\mu}\mu$ is the
similar spin generating series for vertical \r-ribbon strips.
\QED

\proclaim Corollary.
For $r\in\Np$ there exists a spin preserving asymmetric Knuth correspondence
for \r-ribbon tableaux: for any \r-core $\gamma$, and any $m,n\in\N$, there is
a bijection between $m\times n$ matrices~$A$ with entries in~$\{0,1\}^r$, and
pairs consisting of a semistandard \r-ribbon tableau~$P$ and a transpose
semistandard \r-ribbon tableau~$Q$ of equal shape with
$\wt(P)=\(\sum_{i\in\set{m}}|A_{i,j}|\){}_{j\in\set{n}}$ and
$\wt(Q)=\(\sum_{j\in\set{n}}|A_{i,j}|\){}_{i\in\set{m}}$, which is such that
in addition one has
$\spin(P)+\spinv(Q)=\sum_{i,j\in\set{m}\times\set{n}}\n(A_{i,j})$.
\QED

This corollary gives a counterpart of (\spinCauchyKnuthid), which uses
besides $G^{(r)}_{\\/\gamma}(q^{1\over2},X)$
the spin generating series $\widetilde G^{(r)}_{\\/\gamma}(q^{1\over2},Y)
=\sum_Qq^{\spinv(Q)}Y^{\wt(Q)}$ of transpose
semistandard \r-ribbon tableaux of shape~$\\/\gamma$; it reads
$$
 \prod_{i,j\in\N}\prod_{k\in\set{r}}(1+q^kX_iY_j)
 =\sum_{\\\in\RH(\gamma)}G^{(r)}_{\\/\gamma}(q^{1\over2},X)
              \widetilde G^{(r)}_{\\/\gamma}(q^{1\over2},Y)
.\nn
$$

\newsection Conclusion.

We have defined constructions showing that symmetric and asymmetric Knuth
correspondences for semistandard \r-ribbon tableaux, of which trivial examples
can be obtained using the \r-quotient map, can be refined so as to satisfy the
additional requirement that the spins of tableaux be respected. More
interesting than the mere existence of such correspondences is the fact that
this requirement forces defining the basic ingredient (the shape datum) for
the construction by a radically different method than the traditional one of
applying a Schensted correspondence to standardisations. Instead, the method
of tracing level changes of a function derived from edge sequences, which is
at the heart of the Shimozono-White \r-correspondence, is extended to deal
with many ribbons at the same time. This shows the importance of edge
sequences for understanding of \r-ribbon tableaux, even when the \r-quotient
map cannot be used.

\Supplement References.

\input \jobname.ref

\iftoc
  \closeout\tocfile \vfil\eject \def\tocitem#1=#2\onpage#3. {\par\line{\hbox
  to 1.5pc{\bf #1\hss}\quad#2\dotfill#3}\ignorespaces} \input \jobname.toc

\fi

\bye